\documentclass[11pt]{amsart}
\usepackage{latexsym}
\usepackage{amscd}
\usepackage[all]{xy}
\usepackage[mathscr]{euscript}

\normalsize

\RequirePackage{xspace}

\newcommand{\thought}[1]{}
\renewcommand{\thought}[1]{ \textbf{[#1]}}

\usepackage{enumerate}        
\newenvironment{roenumerate}{\begin{enumerate}[\upshape (i)]}{\end{enumerate}}

\newcommand\nc {\newcommand}
\newcommand\rnc{\renewcommand}
\nc\script{\mathscr}

\newcount\blopone
\newcount\xone
\newcount\xtwo
\newcount\ytwo

\newtheorem{theorem}{Theorem}[section]
\newtheorem{prop}[theorem]{Proposition}
\newtheorem{equivform}[theorem]{Equivalent Formulation}

\newtheorem{refinement}[theorem]{Refinement}
\newtheorem{summary}[theorem]{Summary}
\newtheorem{importnota}[theorem]{Important Notation}
\newtheorem{prblm}[theorem]{Problem}
\newtheorem{notation}[theorem]{Notation}
\newtheorem{defin}[theorem]{Definition}
\newtheorem{caution}[theorem]{Caution}
\newtheorem{remark}[theorem]{Remark}
\newtheorem{reminder}[theorem]{Reminder}
\newtheorem{lemma}[theorem]{Lemma}
\newtheorem{construction}[theorem]{Construction}
\newtheorem{corollary}[theorem]{Corollary}
\newtheorem{example}[theorem]{Example}
\newtheorem{conclusion}[theorem]{Conclusion}
\newtheorem{triviality}[theorem]{Triviality}
\newtheorem{proto}[theorem]{Prototype Quasifibration}
\newtheorem{cauex}[theorem]{Cautionary Example}
\newtheorem{hypo}[theorem]{Hypothesis}
\newtheorem{subth}{ }[theorem]
\newtheorem{case}{Case}[theorem]
\newtheorem{ssubth}{ }[subth]

\nc\tri[1]{\begin{triviality}
\label{#1}}
\nc\eqf[1]{\begin{equivform}
\label{#1}\begin{em}}
\nc\cas[1]{\begin{case}
\label{#1}
\begin{em}}
\nc\rfn[1]{\begin{refinement}
\label{#1}}
\nc\prt[1]{\begin{proto}
\label{#1}}
\nc\lem[1]{\begin{lemma}
\label{#1}}
\nc\pro[1]{\begin{prop}
\label{#1}}
\nc\thm[1]{\begin{theorem}
\label{#1}}
\nc\cor[1]{\begin{corollary}
\label{#1}}
\nc\dfn[1]{\begin{defin}
\begin{em}
\label{#1}}
\nc\sthm[1]{\begin{subth}
\label{#1}}
\nc\exm[1]{\begin{example}
\label{#1}
\begin{em}}
\nc\plm[1]{\begin{prblm}
\label{#1}
\begin{em}}
\nc\rmk[1]{\begin{remark}
\label{#1}
\begin{em}}
\nc\rmd[1]{\begin{reminder}
\label{#1}
\begin{em}}
\nc\ntn[1]{\begin{notation}
\label{#1}
\begin{em}}
\nc\smr[1]{\begin{summary}
\label{#1}
\begin{em}}
\nc\cau[1]{\begin{caution}
\label{#1}
\begin{em}}
\nc\hyp[1]{\begin{hypo}
\label{#1}\begin{em}}
\nc\imn[1]{\begin{importnota}
\label{#1}
\begin{em}}
\nc\cax[1]{\begin{cauex}
\label{#1}
\begin{em}}
\nc\con[1]{\begin{construction}
\label{#1}
\begin{em}}
\nc\ssthm[1]{\begin{ssubth}
\label{#1}
\begin{em}}
\nc\cnc[1]{\begin{conclusion}
\label{#1}
\begin{em}}

\nc\elem{\end{lemma}}
\nc\eeqf{\blimy \end{em}
\end{equivform}}
\nc\erfn{\end{refinement}}
\nc\eprt{\end{proto}}
\nc\ethm{\end{theorem}}
\nc\ecor{\end{corollary}}
\nc\edfn{\blimy \end{em}
\end{defin}}
\nc\esthm{\end{subth}}
\nc\epro{\end{prop}}
\nc\etri{\end{triviality}}
\nc\eexm{\blimy \end{em}
\end{example}}
\nc\ermk{\blimy \end{em}
\end{remark}}
\nc\ermd{\blimy \end{em}
\end{reminder}}
\nc\eplm{\end{em}
\end{prblm}}
\nc\ecas{\end{em}
\end{case}}
\nc\ecau{\end{em}
\end{caution}}
\nc\ecax{\end{em}
\end{cauex}}
\nc\eimn{\end{em}
\end{importnota}}
\nc\entn{\blimy \end{em}
\end{notation}}
\nc\econ{\end{em}
\end{construction}}
\nc\esmr{\end{em}
\end{summary}}
\nc\ehyp{\blimy \end{em}
\end{hypo}}
\nc\ecnc{\end{em}
\end{conclusion}}
\nc\essthm{\end{em}
\end{ssubth}}

\nc\blimy{{\begin{flushright} $\Box$\par
\end{flushright}}\vskip 2mm }

\newtheorem{hypothesis}[theorem]{Hypothesis}
\nc\sst{\scriptstyle}
\newcommand{\comment}[1]{}
\newcommand{\ri}{\longrightarrow}

\newcommand{\zz}{{\mathbb Z}}

\nc\bR{{\mathbf R}}
\nc\bS{{\mathbf S}}
\nc\bT{{\mathbf T}}
\nc\bU{{\mathbf U}}
\nc\z{\zeta}
\nc\bc{{\mathbb{BC}}}
\nc\ct{{\script T}}
\nc\cs{{\script S}}
\nc\car{{\script R}}
\nc\ca{{\script A}}
\nc\cb{{\script B}}
\nc\cc{{\script C}}
\nc\cd{{\script D}}
\nc\ce{{\script E}}
\nc\ci{{\script I}}
\nc\co{{\script O}}
\nc\bZ{{\mathbb Z}}
\nc\bd{\begin{description}}
\nc\ed{\end{description}}
\nc\ctob{{\script C}at\big(\ci^{op},\ca\big)}
\nc\clim{{\ds\mathop{\rm lim}_{\ds\longleftarrow}}}
\nc\climi{\clim^{\!i}\,}
\nc\climn{\clim^{\!n}\,}
\nc\colim{{\ds\mathop{\rm colim}_{\ds\la}}}
\nc\oa{\overline{\ca}}
\nc\s{\sigma}
\nc\ta{\tau}
\nc\os{\overline\sigma}
\nc\ot{\overline\tau}
\nc\T{\Sigma}
\nc\de[1]{{\mathop{\rm deg(#1)}}}
\nc\Ad[1]{\mathop{\rm Ad}(#1)}
\nc\ad[1]{\mathop{\rm ad}(#1)}
\nc\Tor{\text{\rm Tor}}

\def\der #1 {D\left(#1\right)}
\nc\prf{\begin{proof}}
\nc\eprf{\end{proof}}
\nc\ds{\displaystyle}

\nc\ab{{\script A}b}

\nc\csab{{\script C}at\big(\cs^{op},\ab\big)}
\nc\ctab{{\script C}at\Big({\{\ct^\alpha\}}^{op},\ab\Big)}
\nc\csex{{\script E}x\big(\cs^{op},\ab\big)}
\nc\ctex{{\script E}x\Big({\{\ct^\alpha\}}^{op},\ab\Big)}
\nc\sub{\qquad\subset\qquad}
\nc\ctr[1]{{\left.\ct\left(-,#1\right)\right|}_{\cs}}
\nc\ctrf[2]{{\left.\ct\left(#1,#2\right)\right|}_{\cs}}
\nc\Ctr[1]{{\left.\ct\left(-,#1\right)\right|}_{\ct^\alpha}}
\nc\Ctrf[2]{{\left.\ct\left(#1,#2\right)\right|}_{\ct^\alpha}}

\nc\la{\longrightarrow}
\nc\oti{{^L\otimes_R^{}}}
\nc\rs{\s^{-1}R}
\nc\br{{\{\s^{-1}R\}}}

\nc\nin{\noindent}
\nc\cad[1]{\text{card}(#1)}
\nc\eq{\quad=\quad}
\nc\BA{\begin{array}{c}}
\nc\EA{\end{array}}
\nc\kth{{\it K}--theory}

\nc\barr{
\[
\begin{array}{cccccccccccccccc}
}
\nc\earr{
\end{array}
\]
}
\nc\as[1]{{\langle S\rangle}^{#1}}
\nc\sh{\hbox{\it shift}}

\nc\yy[1]{{\left.\ct\left(-,#1\right)\right|}_{\ct^c}}
\nc\vrep[2]{{\left.\ct\left(#1,#2\right)\right|}_{\ct^\alpha}}
\nc\da{\downarrow}
\nc\Hom{{\mathop{\rm Hom}}}
\nc\End{{\mathop{\rm End}}}
\nc\Ext{{\mathop{\rm Ext}}}
\nc\mr\Modtc
\nc\PExt{{\mathop{\rm PExt}}}
\nc\bA{{\mathbf A}}
\nc\bB{{\mathbf B}}
\nc\bC{{\mathbf C}}
\nc\bD{{\mathbf D}}

\nc\ra\ri
\nc\y[1]{\mathbf{y}#1}
\nc\x[1]{\mathbf{z}#1}
\nc\Mod[1]{\ensuremath{\mathop{\textup{Mod-}#1}}\xspace}
\nc\Md {\ensuremath{\mathop{\textup{Mod}}}}
\rnc\mod[1]{\ensuremath{\mathop{\textup{mod-}#1}}\xspace}
\nc\Modtc{\Mod{\ct^c}}
\nc\pgldim[1]{\mathop{\rm pgldim}\,#1}

\begin{document}

\author{Amnon Neeman}
\address{A.N. : Center for Mathematics and its Applications \newline
\indent School of Mathematical Sciences\newline
\indent John Dedman Building\newline
\indent The Australian National University\newline
\indent Canberra, ACT 0200\newline
\indent AUSTRALIA}
\email{Amnon.Neeman@anu.edu.au}
\author{Andrew Ranicki}
\address{A.R. : Department of Mathematics and Statistics \newline
\indent University of Edinburgh \newline
\indent James Clerk Maxwell Building \newline
\indent King's Buildings \newline
\indent Mayfield Road \newline
\indent Edinburgh EH9 3JZ\newline
\indent SCOTLAND, UK}
\email{aar@maths.ed.ac.uk}

\title[Noncommutative localization I]
{Noncommutative localization and chain complexes\\
I. Algebraic $K$- and $L$-theory}
\date{18 September, 2001}

\begin{abstract}\hskip2mm
The noncommutative (Cohn) localization $\sigma^{-1}R$ of a ring $R$ is
defined for any collection $\sigma$ of morphisms of f.g.  projective
left $R$-modules. We exhibit $\sigma^{-1}R$ as the endomorphism ring 
of $R$ in an appropriate triangulated category.  We use this
expression to prove that if
$\text{Tor}^R_i(\sigma^{-1}R,\sigma^{-1}R)=0$ for $i \geq 1$
then every bounded f.g.  projective $\sigma^{-1}R$-module chain complex $D$ 
with $[D] \in \text{im}(K_0(R)\ri K_0(\sigma^{-1}R))$ is chain equivalent 
to $\sigma^{-1}C$ for a bounded f.g.  projective $R$-module chain complex $C$, 
and that there is a localization exact sequence in higher algebraic $K$-theory
$$\dots \ri K_n(R) \ri K_n(\sigma^{-1}R) \ri K_n(R,\sigma) \ri
K_{n-1}(R) \ri \dots~,$$
extending to the left the sequence obtained for $n\leq1$ by Schofield.
For a noncommutative localization $\sigma^{-1}R$ of a
ring with involution $R$ there are analogous results for algebraic
$L$-theory, extending the results of Vogel from quadratic to symmetric
$L$-theory.
\end{abstract}


\maketitle

\tableofcontents

\section*{Introduction}
\label{S0}

This is the first of a series of papers on the algebraic $K$- and
$L$-theory of noncommutative localizations of rings. We adopt
throughout the following convention. Suppose $R$ is an associative
ring. Unless otherwise specified,
by $R$-module we shall mean left $R$-module.
\medskip 

Let $\sigma\subset R$ be a multiplicative set of elements in the centre
$Z(R)$ of the ring $R$.  It is very classical to define $\rs$ as the
ring of fractions $r/s$, with $r\in R$ and $s\in\s$.  The ring $\rs$ is
called the {\em commutative localization} of $R$ with respect to the
multiplicative set $\s$.  Note that the rings $R$ and $\rs$ are not
assumed commutative; the only commutativity hypothesis is that
$\s\subset Z(R)$.  The localization exact sequences relating the
algebraic $K$- and $L$-groups of $R$ and $\sigma^{-1}R$ are basic
computational tools.  \medskip

More recently, it has turned out to be useful to generalise the notion
of rings of quotients, in which much more general $\sigma$'s are
allowed.  From now on, the elements of $\sigma$ will be maps $s:P\ri
Q$, with $P$ and $Q$ f.g.  projective $R$-modules, and localization
will invert these maps.  The classical case of a multiplicative set
$\s\subset R$ is just the special case where $P=Q=R$, in other words
$P$ and $Q$ are free $R$ modules of rank 1.  The morphisms $s:R\ri R$
to be inverted are given by right multiplication by $s\in\s$.  If all
$s\in\s$ lie in the centre of $R$, we are in the traditional situation. 
\medskip

Noncommutative localization is characterized by the following
universal property. A ring homomorphism $R\ri S$ is called {\it
$\sigma$-inverting\/} if $1\otimes s:S\otimes_R P\ri S\otimes_R Q$ is
an $S$-module isomorphism for every $s:P \ri Q$ in $\s$.  The category
of $\sigma$-inverting ring homomorphisms $R \to S$ has an initial object,
denoted $R \ri \sigma^{-1}R$.  This means that any $\sigma$-inverting
ring homomorphism $R \ri S$ factors uniquely as $R \ri \sigma^{-1}R \ri S$. 
The ring $\sigma^{-1}R$ is called a {\it noncommutative localization\/}
or a {\it universal localization\/} of $R$ inverting $\sigma$. 
\medskip

Noncommutative localization was pioneered by Ore \cite{Ore} and Cohn
\cite{Cohn2}, in order to study embeddings of noncommutative rings in
skewfields.  See Ranicki \cite{Ranicki1998} for some of the
applications of the algebraic $K$- and $L$-theory of noncommutative
localization to the topology of codimension 2 submanifolds, such as
knots.  
\medskip

In Part I of the paper we study the algebraic $K$- and $L$-theory of a
noncommutative localization $\sigma^{-1}R$ by means of triangulated
categories, generalizing the work of Vogel \cite{Vogel1982} and
Schofield \cite{Schofield}.  In Part~II we shall obtain a chain complex
interpretation of the normal form of Gerasimov \cite{Gerasimov} and
Malcolmson \cite{Malcolmson} for elements of $\sigma^{-1}R$.  
\medskip

The ring homomorphism $R \ri \sigma^{-1}R$ gives $\sigma^{-1}R$ the
structure of a right $R$-module in the usual manner, and we have a
functor
$$\begin{array}{l}
\sigma^{-1}~=~\sigma^{-1}R\otimes_R-~
:~\{\hbox{$R$-modules}\} \ri 
\{\hbox{$\sigma^{-1}R$-modules}\}~;\\[1ex]
\hphantom{\sigma^{-1}~=~\sigma^{-1}R\otimes_R-~:~\{\hbox{$R$-modules}\} 
\ri}
M \mapsto \sigma^{-1}M~=~\br\otimes_RM~.
\end{array}$$
A $\sigma^{-1}R$-module is {\it induced f.g. projective} 
if it is of the form $\sigma^{-1}P$ for a f.g.  projective $R$-module $P$.  
\medskip

The {\it chain complex lifting problem} is to decide if a bounded chain
complex $D$ of induced f.g.  projective $\sigma^{-1}R$-modules is chain
equivalent to $\sigma^{-1}C$ for a bounded chain complex $C$ of f.g. 
projective $R$-modules.  The problem has a trivial affirmative solution
for a commutative localization, by the clearing of denominators, when
$D$ is actually isomorphic to $\sigma^{-1}C$.  In Part I of the paper
we apply triangulated categories to study the problem for a
noncommutative localization $\sigma^{-1}R$.  
\medskip

A systematic solution of the chain complex lifting problem leads to the
extensions to the noncommutative case of the localization exact
sequences for the algebraic $K$- and $L$-groups of a commutative
localization.  It turns out that the problem has a systematic solution
if and only if $\sigma^{-1}R$ is `stably flat over $R$', and that
there are homological obstructions to stable flatness in general.  Here
is what we mean by stably flat over $R$.  
\medskip

\dfn{Stableflat}
Let $R\ri S$ be a ring homomorphism. The ring $S$ is called
{\em stably flat over $R$} if~: 
\begin{roenumerate}
\item the multiplication map $\mu:S\otimes_R^{}S\ri S$ is an
  isomorphism, 
\item $\Tor^R_i(S,S)=0$ for all $i\geq1$.
\end{roenumerate}
\edfn

\rmk{Stableflat1}
In the case of a noncommutative localization $R\ri\rs$, it is always
true that $\mu:S\otimes_R^{}S\ri S$ is an isomorphism, and that 
$\Tor^R_1(\rs,\rs)=0$.  A proof may be found on page~58 of
Schofield~\cite{Schofield}, or also in Corollary~\ref{C8.4} of this
article.  In general $\Tor^R_i(\rs,\rs)\neq0$ for $i\geq2$.
\ermk

In fact, Schofield has constructed examples of noncommutative
localizations $\sigma^{-1}R$ which are not stably flat over $R$, with
$\Tor_2^R(\rs,\rs)\neq 0$.  These examples pinpoint the difference
between commutative and noncommutative localization.  In commutative
localization $\sigma^{-1}R$ is always a flat $R$-module (Example
\ref{comm1} below); Schofield's examples show that in noncommutative
localization the ring $\rs$ need not even be stably flat over $R$.  We
shall describe these examples in Part~II.

\exm{comm1} {\rm
Let $\sigma \subset R$ be a multiplicatively closed subset 
with $1 \in \sigma$, such that the elements $s \in \sigma$
are central in $R$ or more generally satisfy the {\it Ore conditions}~:
\begin{itemize}
\item[(i)] for all $r \in R$, $s \in \sigma$ there exist 
$q \in R$, $t \in \sigma$ such that $rt=sq\in R$,
\item[(ii)] if $r \in R$, $s \in \sigma$ are such that $sr=0 \in R$
then $rt=0 \in R$ for some $t \in\sigma$.
\end{itemize}
The {\it Ore localization} $\sigma^{-1}R$ is the ring of fractions 
$r/s$ ($r \in R, s \in \sigma)$,
which are the equivalence classes of pairs $(r,s)$ with
\begin{itemize}
\item[]
$(r,s) \sim (q,t)$ if there exist $u,v \in R$ such that 
$ru=qv \in R$, $su=tv \in \sigma$.
\end{itemize}
An Ore localization $\sigma^{-1}R$ is flat (Stenstr\"om  
\cite{Stenstrom}, Prop. II.3.5) and hence stably flat over $R$.}
\eexm

\exm{E0.2} Recall that a ring $R$ is called {\em hereditary} if all
$R$-modules have projective dimension $\leq1$.  A noncommutative
localization of a hereditary ring $R$ is a hereditary ring
$\sigma^{-1}R$ (Bergman and Dicks \cite{Bergman-Dicks}).  In this case,
the vanishing of $\Tor^R_i(\rs,\rs)$ with $i\geq2$ follows just from
the fact that every $R$-module has projective dimension $\leq1$.  Thus
$\rs$ is stably flat over $R$.  
\eexm

\exm{E0.3} For $\mu \geq 1$ let $F_{\mu}$ be the free group on $\mu$
generators.  Given a commutative ring $k$ let $R=k[F_\mu]$, and let $\sigma$ be
the set of all $R$-module morphisms $s:R^n \ri R^n$ ($n \geq 1$)
inducing $k$-module isomorphisms $\epsilon(s):k^n \ri k^n$ via the
augmentation $\epsilon:R \ri k$.  The noncommutative localization
$\sigma^{-1}R$ is flat for $\mu=1$ (when it is commutative), 
but is not flat for $\mu \geq 2$. If $k$ is a principal ideal domain
Farber and Vogel \cite{FarberVogel}
identified $\sigma^{-1}R$ with the ring of rational functions in $\mu$
non-commuting variables, and proved that $\sigma^{-1}R$ is stably 
flat over $R$. (If $k$ is a field then $k[F_\mu]$ is hereditary
by Cohn \cite{Cohn1}, and the stable flatness is given by Example
\ref{E0.2}).
\eexm

Now that we have seen some examples of stably flat localizations
$R\ri \rs$, it is time for a comment. The main results of
the article are about stably flat localizations. We do have
some weak results that hold without the hypothesis of stable
flatness; but the powerful theorems assume that $\rs$ is stably
flat over $R$. It becomes interesting to find equivalent
formulations of the hypothesis that $\rs$ is stably flat over $R$.

\eqf{EQ0.05}
Let $R \ri S$ be a ring homomorphism. The ring $S$ is stably flat over
$R$ if and only if
\begin{roenumerate}
\item the multiplication map $\mu:S\otimes_R^{}S\ri S$ is an isomorphism,
\item for all $S$-modules $M$ and all $i\geq1$, we have $\Tor^R_i(S,M)=0$.
\end{roenumerate}
The if implication is trivial; if $\Tor^R_i(S,M)$ vanishes for all
$S$-modules $M$, then it vanishes in particular for $M=S$.  The only if
implication may be found in Lemma~\ref{L8.5}.  
\eeqf

\rmk{Rtermin}
(i) Equivalent Formulation~\ref{EQ0.05} should explain the terminology. 
$S$ is flat as a right $R$-module if and only if $\Tor^R_i(S,M)$
vanishes for all $R$-modules $M$, while $S$ is stably flat over $R$ if
and only if $\Tor^R_i(S,M)$ vanishes for all $S$-modules $M$.  In
\ref{EQ0.1} we shall see an equivalent formulation of stable flatness,
this time in terms of Vogel's chain complex $E(C)$.  It seems natural
to postpone this equivalent formulation until we are ready to discuss
Vogel's construction.  \\
(ii) A ring homomorphism $R \ri S$ such that ${\rm Tor}_i^R(S,M)=0$ for all
$S$-modules $M$ and all $i \geq 1$ is called a {\it lifting} by Dicks 
\cite{Dicks} (p. 565).
\ermk

Let $\cc$ be a triangulated category. We say that $\cc$ satisfies
[TR5] if arbitrary coproducts exist in $\cc$.
An object $c$ in a triangulated category $\cc$ satisfying [TR5]
is called {\em compact} if
$$\cc\left(c,\coprod_{\lambda \in \Lambda}t_{\lambda}\right)~=~
\bigoplus\limits_{\lambda \in \Lambda} \cc(c,t_{\lambda})$$
for every collection $\{t_{\lambda}\,\vert\,\lambda \in \Lambda\}$ of
objects in $\cc$.  Let $\cc^c \subset \cc$ be the full subcategory of
all the compact objects in $\cc$.  For the derived category $D(R)$ of
unbounded $R$-module chain complexes, the compact category
$D(R)^c=D^c(R) \subset D(R)$ has for its objects the bounded f.g. 
projective $R$-module chain complexes, and any objects isomorphic to
these.  Let $D(R,\sigma) \subset D(R)$, $D^c(R,\sigma) \subset D^c(R)$
be the subcategories generated by $\sigma$.  The objects of
$D^c(R,\sigma)$ are the bounded f.g.  projective $R$-module chain
complexes $C$ such that $H_*(\sigma^{-1}C)=0$ (Proposition \ref{P5.3}). 
We shall be working with triangulated categories and functors
$$\xymatrix{
\ca^c=D^c(R,\sigma) \ar[r] \ar[d] & \cb^c=D^c(R) 
\ar[r]^-{\displaystyle{\pi}}\ar[d] &
\cc^c \ar[r]^-{\displaystyle{T}}  \ar[d] & 
\cd^c=D^c(\sigma^{-1}R)\ar[d] \\
\ca=D(R,\sigma) \ar[r]
& \cb=D(R) \ar[r]^-{\displaystyle{\pi}}
 & \cc=D(R)/D(R,\sigma) 
\ar[r]^-{\displaystyle{T}} & \cd=D(\sigma^{-1}R)}$$
The unnamed maps are inclusions.  The natural map $\cb^c/\ca^c\ri
\cc^c$ is an idempotent completion (see \ref{T3.8.4}).  The functor
$\pi:\cb \ri \cc=\cb/\ca$ is the projection to the quotient.  The
functor
$$T\pi~:~\cb \ri \cd~;~
X \mapsto \{\sigma^{-1}R\}^L\!\otimes_RX~=~\sigma^{-1}P$$
is given by the tensor ${~}^L\otimes_R$ in the derived category,
constructed using a sufficiently nice projective $R$-module chain
complex $P$ with a homology equivalence $P \ri X$.  And the functor $T$
is determined as the unique factorization of $T\pi$ through $\pi$.  We
shall be particularly concerned with the extent to which
$T=\sigma^{-1}:\cc^c \ri \cd^c$ is an equivalence of triangulated
categories.
\medskip

Here is our main result (Theorem \ref{Tfinal})~:
\medskip

\thm{0.1} The following conditions on a noncommutative localization
$R \ri \sigma^{-1}R$ are equivalent~:
\begin{itemize}
\item[(i)] $\sigma^{-1}R$ is stably flat over $R$,
\item[(ii)] the functor $T:\cc^c \ri \cd^c$ is an equivalence of
categories.
\end{itemize} 
\hfill $\Box$ \ethm

Our main Theorem \ref{0.1} has an immediate consequence (Theorem
\ref{Plifting})~:

\thm{0.2}
If $\sigma^{-1}R$ is stably flat over $R$ then the chain
complex lifting problem can always be solved~: for every bounded chain
complex $D$ of induced f.g.  projective $\sigma^{-1}R$-modules there
exists a bounded f.g.  projective $R$-module chain complex $C$ with a
chain equivalence $\sigma^{-1}C \simeq D$.  
\hfill$\Box$ 
\ethm

Without the stable flatness hypothesis there are Toda bracket
obstructions to lifting chain complexes; we explain this briefly
in Section~\ref{lift2}. A much fuller treatment will come in
Part~II of this series.
\medskip

Let $H(R,\sigma)$ be the exact category of $\sigma$-torsion $R$-modules
$T$ of projective dimension 1, i.e.  the $R$-modules with a f.g. 
projective $R$-module resolution
$$\xymatrix{0 \ar[r]& P \ar[r]^{\displaystyle{s}} & Q \ar[r] & T 
\ar[r]&0}$$
such that $\sigma^{-1}s:\sigma^{-1}P \ri \sigma^{-1}Q$ 
is an isomorphism (e.g.  if
$s \in \sigma$).  We shall only be dealing with $H(R,\sigma)$ in the
special case of a noncommutative localization $\sigma^{-1}R$ when each
morphism $s:P \ri Q$ in $\sigma$ is injective. This happens, for 
example, if  $R \ri
\sigma^{-1}R$ is injective (see Proposition~\ref{inj2} for details).
  \medskip

The algebraic $K$-theory localization exact sequence for an injective
Ore localization $R \ri \sigma^{-1}R$
$$\dots \ri K_n(R) \ri K_n(\sigma^{-1}R) \ri K_n(R,\sigma) \ri 
K_{n-1}(R) \ri \dots$$
was obtained by Bass \cite{Bass} for $n=1$ and Gersten
\cite{Gersten74}, Quillen \cite{Quillen} and Grayson \cite{Grayson76},
\cite{Grayson80} for $n \geq 2$, with $K_*(R,\sigma)=K_{*-1}(H(R,\sigma))$. 
Schofield \cite{Schofield} established the algebraic $K$-theory
localization exact sequence in the classical dimensions 0,1
$$K_1(R) \ri K_1(\sigma^{-1}R) \ri K_1(R,\sigma) \ri K_0(R) \ri
K_0(\sigma^{-1}R)$$
for any injective noncommutative localization  $R\ri\sigma^{-1}R$.
It is easy to extend the exact sequence to the right, using the
lower $K$-groups $K_{-*}$ of \cite{Bass}.
\medskip

Waldhausen \cite{Waldhausen85} identified the algebraic $K$-groups of
$D^c(R)$ with the Quillen \cite{Quillen} algebraic $K$-groups of $R$
$$K_*(D^c(R))~=~K_*(R)~.$$
By the localization theorem of \cite{Waldhausen85}, for any $R,\sigma$
there is defined a long exact sequence in algebraic $K$-theory
$$\dots \ri K_n(R) \ri 
K_n(D^c(R)/D^c(R,\sigma)) \ri K_n(R,\sigma) \ri K_{n-1}(R) \ri 
\dots~.$$
with 
$$K_*(R,\sigma)~=~K_{*-1}(D^c(R,\sigma))~~({\rm definition})~.$$
The idempotent completion $D^c(R)/D^c(R,\sigma)\ri\cc^c$ 
induces an isomorphism on higher \kth, so that
$$K_*(D^c(R)/D^c(R,\sigma))~=~K_*(\cc^c)~.$$ 
The map $T:\cc^c\ri D^c(\rs)$ induces a map on Waldhausen's $K$--theory
$$K(T)~:~K(\cc^c)\ri K(D^c(\rs))~=~K(\rs)~.$$
In Sections~\ref{S8} and
\ref{S9} we show that this map is an isomorphism on $K_0$ and $K_1$,
without any hypothesis on $\s$.  If we assume that $\rs$ is stably
flat over $R$, then Theorem~\ref{Tfinal} tells us that
$T:\cc^c\ri\cd^c=D^c(\rs)$ is an equivalence of categories, and as a
corollary we deduce Theorem~\ref{Kexact} :

\thm{0.3} If $\sigma^{-1}R$ is stably flat over $R$ the functor 
$T:\cc^c \ri \cd^c=D^c(\sigma^{-1}R)$ induces isomorphisms
$$T~:~K_*(\cc^c)~=~K_*(D^c(R)/D^c(R,\sigma)) \ri 
K_*(\cd^c)~=~K_*(\sigma^{-1}R)$$
and there is a localization exact sequence in algebraic $K$-theory 
$$\dots \ri K_n(R) \ri K_n(\sigma^{-1}R) \ri K_n(R,\sigma) \ri 
K_{n-1}(R) \ri \dots~.\eqno{\Box}$$
\ethm

In Theorem  \ref{torsionK} we shall prove~:

\thm{0.4} If each morphism in $\sigma$ is injective the Waldhausen
$K$-groups of $D^c(R,\sigma)$ are just the Quillen $K$-groups of
$H(R,\sigma)$
$$K_*(R,\sigma)~=~K_{*-1}(D^c(R,\sigma))~=~K_{*-1}(H(R,\sigma))~.$$
If in addition $\sigma^{-1}R$ is stably flat over $R$ 
there is defined a localization exact sequence in algebraic $K$-theory
$$\dots \ri K_n(R) \ri K_n(\sigma^{-1}R) \ri K_{n-1}(H(R,\sigma))
\ri K_{n-1}(R) \ri \dots~.\eqno{\Box}$$
\ethm

However, in general the morphisms $s:P \ri Q$, $R \ri \sigma^{-1}R$ are
not injective and, except in section~\ref{Storsion} and parts of
section~\ref{Ltheory}, we do not assume this to be the case.  
\medskip

In section \ref{Ltheory} we consider the $L$-theory of noncommutative
localizations, obtaining the following results (Theorems 
\ref{L2}, \ref{L3}, \ref{Lfin})~:

\thm{0.5} Let $R \ri  \sigma^{-1}R$ be a noncommutative localization 
of a ring with involution $R$, such that $\sigma$ is invariant under
the involution.\\
{\rm (i)} There is a localization exact sequence of quadratic 
$L$-groups
$$\xymatrix{\dots \ar[r] & L_n(R) \ar[r]&
L_n^I(\sigma^{-1}R) \ar[r]^-{\partial} &
L_n(R,\sigma)\ar[r] & L_{n-1}(R) \ar[r] & \dots}$$
with $I={\rm im}(K_0(R) \ri K_0(\sigma^{-1}R))$, and
$L_n(R,\sigma)$ the cobordism group of $\sigma^{-1}R$-contractible
$(n-1)$-dimensional quadratic Poincar\'e complexes over $R$.\\
{\rm (ii)} If $\sigma^{-1}R$ is stably flat over $R$
there is a localization exact sequence of symmetric $L$-groups
$$\xymatrix{\dots \ar[r] & L^n(R) \ar[r]&
L^n_I(\sigma^{-1}R) \ar[r]^-{\partial} &
L^n(R,\sigma)\ar[r] & L^{n-1}(R) \ar[r] & \dots}$$
with $L^n(R,\sigma)$ the cobordism group of 
$\sigma^{-1}R$-contractible
$(n-1)$-dimensional symmetric Poincar\'e complexes over $R$.\\
{\rm (iii)} If $R \ri \sigma^{-1}R$ is injective then $L^n(R,\sigma)$ 
(resp. $L_n(R,\sigma)$) is the cobordism group of  $n$-dimensional 
symmetric 
(resp. quadratic) Poincar\'e complexes of $\sigma$-torsion
$R$-modules of projective dimension 1.
\hfill$\Box$
\ethm

The $L$-theory exact sequences of Theorem \ref{0.5} for an injective
Ore localization $R \ri \sigma^{-1}R$ were obtained in Ranicki
\cite{Ranicki1981}.  The quadratic $L$-theory exact sequence of
\ref{0.5} (i) for arbitrary $R \ri \sigma^{-1}R$ was obtained by Vogel
\cite{Vogel1980}, \cite{Vogel1982}.  The symmetric $L$-theory exact
sequence of \ref{0.5} (ii) is new.
\medskip

We shall now give yet another equivalent formulation of stable
flatness.  However, in Part I, we will use Definition~\ref{Stableflat}
as the working definition.  
\medskip

Recall that a right $R$-module $S$ is flat if and only if, 
for any chain complex $C$ of left $R$-modules
$$H_*(S\oti C)~=~S\otimes_R^{}H_*(C)$$
where $S^L$ is a projective right $R$-module resolution of $S$.

\eqf{EQ0.1} For any noncommutative localization $\sigma^{-1}R$
and any chain complex of $R$-modules $C$, we define a contravariant
functor 
$$\begin{array}{c}
[[-,C]]~:~\{\hbox{$R$-module chain complexes}\}
\ri \{\hbox{$\bZ$-modules}\}~;\\[1ex]
A \mapsto [[A,C]]~=~\mathop{\varinjlim}\limits_{(B,\beta)}[A,B]
\end{array}$$
with $[A,B]$ the $\bZ$-module of chain homotopy classes of chain maps
$A \ri B$, and the direct limit taken over all the $R$-module chain
complexes $B$ with a chain map $\beta:C \ri B$ such that the mapping
cone of $\beta$ is quasi-isomorphic to a bounded complex of f.g. 
projective $R$-modules, and $H_*(S\oti C) \cong H_*(S\oti B)$ with
$S=\sigma^{-1}R$.  (A chain map is a quasi-isomorphism if it induces
isomorphisms in homology).  This functor is representable.  That is
$$[[A,C]]~=~[A,E(C)]$$
for some $R$-module chain complex $E(C)=\mathop{\varinjlim} B$.  Such a
complex $E(C)$ was first constructed by Vogel in his
paper~\cite{Vogel1982}.  There is a map of $R$-module chain complexes
$C\ri E(C)$ inducing $S$-module isomorphisms
$$H_*\big(S\oti C\big)~\cong~H_*\big(S\oti E(C)\big)$$
and such that for each $i \in \zz$
$$H_i(E(C))~=~[[\Sigma^iR,C]]~=~\mathop{\varinjlim}\limits_{(B,\beta)}H_i(B)$$
is an $S$-module.  Furthermore, $H_0(E(R))=S$.  If $H_i(C)$ are
$S$-modules for all $i\in\zz$, then the map $C\ri E(C)$ is a
quasi-isomorphism. There is a map of $R$-module chain complexes $E(C)\ri S\oti C$ inducing
$S$-module morphisms
$$H_*(E(C)) \ri H_*(S\oti C)~\cong~H_*\big(S\oti E(C)\big)~.$$
Here is the equivalent formulation~: the noncomutative localization
$S=\rs$ is stably flat over $R$ if and only if $E(R) \ri S\oti R=S$ is
a quasi-isomorphism. The proof of this assertion
is the equivalence of (i) and (ii) in Theorem~\ref{Tfinal}, coupled
with the fact that $E(R)=G\pi R$.
 If $S$ is stably flat then $E(C) \ri S\oti C$ is
a quasi-isomorphism for any $R$-module chain complex $C$. 
The proof of this statement comes from the fact that
the full subcategory on which the map is an isomorphism is
triangulated and closed under coproducts, coupled with 
Lemma~\ref{L3.1.7}.
 In Part~II
of this series we plan to explore in far greater depth the relation
between our work and the earlier work in the subject, especially the
important contributions of Vogel.  
\eeqf

So far, we have used chain complexes and homology.  This makes the
notation above consistent with the the usage standard in the
$L$--theory literature.  For triangulated categories, it is more usual
to work with cochain complexes and cohomology, and in Part I we shall
be working with these (except in the $L$-theory section \ref{Ltheory}). 
The sole object is to make it easier for the reader to check our
references to the literature.  
\medskip

The basic tool in the proof of our main Theorem \ref{0.1} is the fact 
that
$$\pi~:~\cb~=~D(R) \ri \cc~=~D(R)/D(R,\sigma)$$
admits a Bousfield localization. That is, $\pi$ has a right adjoint 
$G:\cc  \ri \cb$,  meaning 
$$\cb(x,Gy)~=~\cc(\pi x,y)$$
for any objects $x \in \cb$, $y \in \cc$.  The $R$-module cochain
complex $G\pi R$ (which is essentially the same as Vogel's
$E(R)$; see Equivalent Formulation \ref{EQ0.1} above) has the following 
properties~:
\begin{itemize}
\item[(i)] $H^{-i}(G\pi R)$ is the group of morphisms $\Sigma^i R \ri
R$ in $\cc$.  Such morphisms are equivalence classes of pairs of
$R$-module chain maps $(\alpha:R \ri Y,\beta:\Sigma^i R \ri Y)$ 
with $Y \in D^c(R)$, $C(\alpha) \in D^c(R,\sigma)$.
The fact that $Y$ may
be taken to be in $D^c(R)\subset D(R)$ 
is not supposed to be trivial: see \ref{T3.8.4}.
\item[(ii)] The groups $H^*(G\pi R)$, ${\rm
Tor}_*^R(\sigma^{-1}R,\sigma^{-1}R)$ are 
$\sigma^{-1}R$--$\sigma^{-1}R$
bimodules, and there is a cohomology spectral sequence
$$E^{i,j}_2~=~
{\rm Tor}^R_{-i}(\sigma^{-1}R,H^j(G\pi R)) \Longrightarrow 
H^*(\{\sigma^{-1}R\}^L\!\otimes_RG\pi R)~=~\sigma^{-1}R~,$$
with $\sigma^{-1}R$ concentrated in degree 0. In particular
$$H^0(G\pi R)~=~\sigma^{-1}R~.$$
\item[(iii)] $H^i(G\pi R)=0$ for $i\neq0$ if and only if $\sigma^{-1}R$
is stably flat over $R$.
\end{itemize}
\medskip

\section{Noncommutative localization of a ring}
\label{S1}

Given a ring $R$ and a collection $\sigma$ of morphisms
$s:P \ri Q$ of f.g. projective $R$-modules we recall the
universal property of the noncommutative localization $\sigma^{-1}R$,
and the original construction.

\dfn{D1} 
(i) A ring homomorphism $R \ri S$ is {\it $\sigma$-inverting} if for
every $s:P\ri Q$ in $\sigma$ the induced $S$-module morphism $1 \otimes
s:S\otimes_RP \ri S\otimes_RQ$ is an isomorphism.\\
(ii) A ring homomorphism $R \ri S$ is {\it universally $\sigma$-inverting}
if it is $\sigma$-inverting, and any other $\sigma$-inverting morphism 
$R \ri S'$ has a unique factorization $R \ri S \ri S'$.
\edfn

Any two universally $\sigma$-inverting ring homomorphisms $R \ri S$,
$R \ri S'$ are related by a canonical isomorphism $S \ri S'$ such
that $R\ri S'$ is the composite $R \ri S \ri S'$.

\thm{Cohn} {\rm  (Cohn \cite{Cohn2})}\\
For any $R$, $\sigma$ there exists a universally $\sigma$-inverting ring
morphism $R \ri \sigma^{-1}R$.\hfill$\Box$
\ethm 

As in the introduction, for any $R$-module $M$
we define $\s^{-1}M=\br\otimes_R^{}M$.
Because $R\ri\rs$ is $\s$-inverting, every $s:P \ri Q$ in $\sigma$
induces an isomorphism $s:\sigma^{-1}P \ri\sigma^{-1}Q$.
\medskip

The original construction of $\sigma^{-1}R$ in \cite{Cohn2} was
for a set $\sigma$ of $R$-module morphisms $s:R^n \ri R^n$, i.e.  for a
set of square matrices $s$, with $\sigma^{-1}R$ obtained from $R$ by
adjoining one generator for each component of a formal inverse $s^{-1}$
and the relations given by
$$s s^{-1}~=~s^{-1} s~=~I~.$$
Gerasimov \cite{Gerasimov}, Malcolmson \cite{Malcolmson} 
and Schofield \cite{Schofield} constructed $\sigma^{-1}R$ for any $\sigma$
as the ring of equivalence classes of triples of morphisms of f.g. 
projective $R$-module morphisms 
$$\big(\, f:P \ri R\,,\, s:P \ri Q\,,\,g:R \ri Q\,\big)$$ 
with $s \in \sigma$. The triple $(f,s,g)$ represents 
$fs^{-1}g \in \sigma^{-1}R$.

\section{Bousfield localization in triangulated categories}
\label{S3}

In section \ref{S4} we shall express a noncommutative localization
$\sigma^{-1}R$ of a ring $R$ as the endomorphism of $R$ in  the 
triangulated category $\cc^c=(D(R)/D(R,\sigma))^c$
$$\sigma^{-1}R~=~{\rm End}_{\cc^c}(R)~.$$ 
The main tool is Bousfield localization, which we
review in this section.  We give careful statements of what is known,
and refer elsewhere for the proofs.

\dfn{D3.1.1}
Let $\cb$ be a triangulated category.  A {\em triangulated subcategory}
is a non-empty full subcategory $\ca\subset\cb$ closed under suspension
and triangles; that is, given a distinguished triangle in $\cb$
$$
\CD
X @>>> Y @>>> Z @>>> \Sigma X~,
\endCD
$$
if $X$ and $Y$ lie in $\ca$ then so do all their suspensions, and so
does $Z$.  
\edfn

\dfn{D3.1.2}
Let $\cb$ be a triangulated category.  A triangulated subcategory
$\ca\subset\cb$ is called {\em thick} (or {\em \'epaisse}) if it
contains all direct summands of its objects.  
\edfn

And now we get to the first theorem~:

\thm{T3.1.3} {\bf (Verdier localization).\ }
Let $\cb$ be a triangulated category, $\ca\subset\cb$
a triangulated subcategory. There is a quotient triangulated
category
$\cc=\cb/\ca$, and a natural triangulated functor
$\cb\ri\cc$. The composite functor
$$
\CD
\ca @>>> \cb @>>> \cc ~=~\cb/\ca 
\endCD
$$
takes every object in $\ca$ to an object isomorphic in $\cc$
to zero. The functor $\cb\ri \cb/\ca$ is universal among
the functors $\cb\ri \cd$ taking the objects of $\ca$
to objects isomorphic to zero. Furthermore, if $\ca\subset\cb$
is thick, then all the objects in $\cb$ whose images in $\cc=\cb/\ca$
are isomorphic to zero lie in $\ca$. 
\ethm

\prf
The theorem is due to Verdier, and may be found in his 
thesis~\cite{Verdier96}.
For a very complete and detailed proof, see Theorem~2.1.8,
on page~74 of \cite{Neeman99}. The full proof occupies 
pages~75-99 op. cit.
\eprf

Suppose $\cc=\cb/\ca$ is as in Theorem~\ref{T3.1.3}. It may
happen that the functor $\pi:\cb\ri\cc$ has a right adjoint,
that is a functor $G:\cc\ri\cb$ such that for every object $x$ in $\cb$
and every object $y$ in $\cc$
$$\cb(x,Gy)~=~\cc(\pi x,y)~.$$
We shall be working in the following situation :

\dfn{D3.1.4}
Let $\cb$ be a triangulated category, and let $\ca\subset\cb$
be a thick subcategory. Let $\cc=\cb/\ca$ be the quotient of
Theorem~\ref{T3.1.3}. We say that a {\em Bousfield localization
functor exists for the pair $\ca\subset\cb$} if the natural
functor $\pi:\cb\ri\cc$ has a right adjoint $G:\cc\ri\cb$.
\edfn  

\rmk{R3.1.4.5}
The adjoint $G$, if it exists, must be a triangulated functor.
Adjoints of triangulated functors are always triangulated. See
Lemma~5.3.6
of \cite{Neeman99}.
\ermk

\nin
Let us next summarise some useful facts about Bousfield localizations.
This is not an exhaustive list; later on in the article
we shall cite more properties. What comes here is a handy 
list of basic, core properties.

\thm{T3.1.5}
Let $\cb$ be a triangulated category, and let $\ca\subset\cb$
be a thick subcategory. Suppose a Bousfield localization
functor exists for the pair $\ca\subset\cb$. That is,
the functor $\pi:\cb\ri\cc=\cb/\ca$ has a right adjoint
$G:\cc\ri\cb$. Then the following
statements are true.
\sthm{T3.1.5.1}
Coproducts in $\cb$ of objects of $\ca$ must lie
in $\ca$.
\esthm
\sthm{T3.1.5.2}
The functor $G$ is fully faithful. 
\esthm
\sthm{T3.1.5.3}
For any two objects $x$ and $y$ in $\cb$, we have an isomorphism
$$\cb(G\pi x,G\pi y) \ri \cb(x,G\pi y)~;~\alpha \mapsto \alpha \eta_x~,$$
where $\eta$ is the unit of adjunction, and $\alpha\eta_x^{}$ is
the composite
$$\alpha \eta_x~:~x \xrightarrow[]{\eta_x} G\pi x \xrightarrow[]{\alpha}
G\pi y~.$$
The assertion is really that any $\beta:x\ri G\pi y$ factors
uniquely as $\alpha\eta_x^{}$.
\esthm
\sthm{T3.1.5.4}
If $b\in\cb$ lies in the image of $G:\cc\ri\cb$, then 
$\eta_t^{}:b\ri G\pi b$ is an isomorphism.
\esthm
\sthm{T3.1.5.5}
Let $b\in\cb$ be any object. The unit of adjunction gives
a map $\eta_b^{}:b\ri G\pi b$. Complete it to a triangle
$$
\CD
a @>>> b @>\eta_b^{}>> G\pi b @>>> \Sigma a~.
\endCD
$$
Then the object $a$ lies in $\ca\subset\cb$.
\esthm
\ethm

\prf
To prove \ref{T3.1.5.1}, just observe that the functor $\pi:\cb\ri\cc$
has a right adjoint, and therefore takes coproducts to coproducts.  Now
$\pi$ takes objects of $\ca$ to zero, and hence takes any coproduct of
objects in $\ca$ to zero.  But because $\ca$ is thick,
Theorem~\ref{T3.1.3} tells us that any object of $\cb$ whose image in
$\cc$ vanishes must lie in $\ca$.  Hence any $\cb$-coproduct of objects
of $\ca$ lies in $\ca$.

Now for \ref{T3.1.5.2}.  In the proof of Lemma~9.1.7 of \cite{Neeman99}
we see that, for all $x\in\cb$, the map $\varepsilon_{\pi x}^{}$ is an
isomorphism.  Any object of $\cc$ is of the form $\pi x$, hence
$\varepsilon:\pi G\Longrightarrow 1$ is an isomorphism.  But then
Lemma~A.2.9 of \cite{Neeman99} tells us that $G$ is fully faithful.

Next comes \ref{T3.1.5.3}.
We have
$$
\begin{array}{rclcl}
\cb(x,G\pi y)&=&\cc(\pi x,\pi y)&\qquad&\text{by adjunction,}\\
&=&\cb(G\pi x,G\pi y)&\qquad&\text{because $G$ is fully faithful.}
\end{array}
$$

\nin
Finally, \ref{T3.1.5.4} may be found in Lemma~9.1.7 of 
\cite{Neeman99}, while \ref{T3.1.5.5} may be found in
Proposition~9.1.8 loc. cit.
\eprf

It might be useful to illustrate these properties in the case of most
interest to us, with $\cb=D(R)$ the (unbounded) derived category of
chain complexes of left $R$-modules.  When we refer to the object $R\in
D(R)$, we mean the chain complex which is $R$ in degree 0, and zero
elsewhere.  Let us make an observation.

\pro{L3.1.6}
Let $\cb=D(R)$ be 
the derived category of a ring $R$. Let $\ca\subset\cb=D(R)$
be a thick subcategory. Suppose a Bousfield localization
functor exists for the pair $\ca\subset D(R)$. That is,
the functor $\pi:\cb\ri\cc=\cb/\ca$ has a right adjoint
$G:\cc\ri\cb$. Then $G\pi R$ is a chain complex in $D(R)$.
We assert that $S=H^0(G\pi R)$ 
has a natural structure of an algebra over $R$, i.e. there
exists a ring homomorphism $R \ri S$.
\epro

\prf
We note first that $H^0(X)=\cb(R,X)$ for every object $X\in\cb=D(R)$.
Applying this to $X=G\pi R$, we have
$$
\begin{array}{rclcl}
H^0(G\pi R)&=&\cb(R,G\pi R)&\qquad&\text{by the above,}\\
&=&\cb(G\pi R,G\pi R)&\qquad&\text{by \ref{T3.1.5.3}.}
\end{array}
$$
Now $S=\cb(G\pi R,G\pi R)$ is a ring, being the endomorphism
ring of an object in an additive category. The fact that $G\pi$
is an additive functor gives us a ring homomorphism
$$
\CD
\cb(R,R) @>>> \cb(G\pi R,G\pi R)~=~S~.
\endCD
$$
Let us agree that we shall view $R$ and $G\pi R$ as right 
modules for, respectively, $\cb(R,R)$ and  $\cb(G\pi R,G\pi R)$.
Then $\cb(R,R)=R$ with the usual right action, and
we have a homomorphism of rings
$$
\CD
R @>>> \cb(G\pi R,G\pi R)~=~S~.
\endCD
$$
\eprf

\lem{L3.1.6.4}
Let the situation be as in Proposition~\ref{L3.1.6}.
We remind the reader: $\cb=D(R)$ is 
the derived category of a ring $R$. Let $\ca\subset\cb=D(R)$
be a thick subcategory. Suppose a Bousfield localization
functor exists for the pair $\ca\subset D(R)$. That is,
the functor $\pi:\cb\ri\cc=\cb/\ca$ has a right adjoint
$G:\cc\ri\cb$. The unit of adjunction $\eta_R^{}:R\ri G\pi R$
induces a map 
$$H^0(\eta_R)~:~R~=~H^0(R) \ri H^0(G\pi R)~.$$
We assert that it agrees with the ring homomorphism
of Proposition~\ref{L3.1.6}.
\elem

\prf
Recall that $R$ can be viewed as $\Hom_R^{}(R,R)$, acting on
the right. For any $r\in R$, right multiplication by $r$ is a 
left-module homomorphism. We denote this homomorphism as
$r:R\ri R$.

The naturality of $\eta$ gives a commutative square
$$
\CD
R @>\eta_R^{}>> G\pi R \\
@VVrV       @VV{G\pi r}V \\
R @>\eta_R^{}>> G\pi R
\endCD
$$
Taking $H^0$, we have a commutative square
$$
\CD
R @>H^0(\eta_R^{})>> H^0(G\pi R) \\
@VVrV       @VV{G\pi r}V \\
R @>H^0(\eta_R^{})>>H^0( G\pi R)
\endCD
$$
Since this commutes for any $r\in R$, we deduce that 
$H^0(\eta_R^{}):R\ri H^0(G\pi R)$ is a homomorphism 
of right $R$-modules. Here, $R$ is a right $R$-modules
in the obvious way, while $S= H^0(G\pi R)$ is a right
module over the ring $S=\cb(G\pi R,G\pi R)$, and the
ring homomorphism $R\ri S$ of Proposition~\ref{L3.1.6} turns
it into a right $R$-module. To prove that the $R$-module
homomorphism $H^0(\eta_R^{}):R\ri H^0(G\pi R)$ agrees
with the ring homomorphism $R\ri S$, we need only check that
$1\in R$ maps to $1\in S$.

But this is easy. For any element $s\in S$,
the identifications
$$
s\in S~=~H^0(G\pi R)~=~\cb(R,G\pi R)~=~\cb(G\pi R,G\pi R)
$$
are quite explicit. Given $\alpha\in \cb(G\pi R,G\pi R)$,
that is a morphism $\alpha:G\pi R\ri G\pi R$, the isomorphism
$$
\cb(R,G\pi R)~=~\cb(G\pi R,G\pi R)
$$
of \ref{T3.1.5.3} sends it to the composite
$$
\CD
R @>\eta_R^{}>> G\pi R @>\alpha>> G\pi R~.
\endCD
$$
In particular, if $\alpha=1_S^{}$, then it is
sent to $\eta_R^{}$. The isomorphism
$$
H^0(G\pi R)~=~\cb(R,G\pi R)
$$
sends $\eta_R^{}\in \cb(R,G\pi R)$ to the image
of $1_R^{}\in R=H^0(R)$ under the map
$$
\CD
H^0(\eta_R^{}):R @>>> H^0(G\pi R)~. 
\endCD
$$
We conclude that under the natural isomorphisms, 
$1_S^{}=H^0(\eta_R^{})(1_R^{})$.
This proves that $H^0(\eta_R^{})$
takes $1\in R$ to $1\in S$.
\eprf

\section{Noncommutative localization using triangulated categories}
\label{S4}

In Lemma~\ref{L3.1.6.4} we learned that, given a suitable subcategory
$\ca\subset\cb=D(R)$ for which the projection $\pi:\cb \ri \cc=\cb/\ca$
has a right adjoint $G:\cc \ri \cb$, there is a canonical ring
homomorphism $R \ri H^0(G\pi R)$. In this section we shall show that
given a set $\sigma$ as in section \ref{S1} and an appropriate
choice of $\ca$, this is the noncommutative localization
$R \ri H^0(G\pi R)=\sigma^{-1}R$ considered in section \ref{S1}.
\medskip 

The first step is to give sufficient conditions for the Bousfield
localization $G$ to exist.  We shall state a general existence theorem,
and then narrow it to the case of interest.  To state the general
theorem, we need to remind the reader of compact objects in
triangulated categories satisfying [TR5].  We recall the definitions.

\dfn{D3.1}
A triangulated category $\cc$ is said to satisfy [TR5] if arbitrary
coproducts exist in $\cc$.
\edfn

\exm{E3.2}
If $R$ is an associative ring and $D(R)$ is its unbounded derived
category, then $D(R)$ satisfies [TR5]. If a triangulated category
satisfies [TR5], then coproducts of distinguished triangles are
distinguished; see Proposition~1.2.1 and Remark~1.2.2 of \cite{Neeman99}.
\eexm

\dfn{D3.3}
Let $\cc$ be a triangulated category satisfying [TR5].
An object $c\in\cc$ is called {\em compact} if the functor
$\cc(c,-)$ respects coproducts. Equivalently, $c$ is compact if
every map
$$
\CD
c @>>> \ds\coprod_{\lambda\in\Lambda}t_\lambda^{}
\endCD
$$
factors through a finite coproduct.
\edfn

\rmk{R3.4}
In the derived category $D(R)$ of Example~\ref{E3.2}, the object $R\in
D(R)$ is compact.  By $R\in D(R)$ we mean the chain complex which is
zero in all dimensions but 0, and is $R$ in dimension 0.  The proof
that $R$ is compact is the following.  For any object $X\in D(R)$,
$$
D(R)\big(R,X\big)~=~H^0(X)~.
$$
Since $H^0(-)$ is a functor commuting with coproducts, the compactness
follows.  All suspensions of a compact object are compact.  
\ermk

It is useful to recall the following fact.

\lem{L3.1.7}
With the notation as in Remark~\ref{R3.4}, the category $D(R)$ contains
an object $R$ (the complex with $R$ in degree zero, and zero in all
other degrees).  Any triangulated subcategory of $D(R)$, closed under
coproducts and containing $R\in D(R)$, is all of $D(R)$.  
\elem

\prf
By Remark~\ref{R3.4}, the object $R$ is compact.  In the terminology of
\cite{Neeman99}, $R$ is an $\aleph_0$-compact object in $D(R)$.  Also,
if for every $n\in\zz$ we have
$$D(R)\big(\Sigma^nR,X\big)~=~D(R)\big(R,\Sigma^{-n}X\big)~=~H^{-n}(X)~=~0~,$$
then $X$ is acyclic and vanishes in $D(R)$.  This makes the set
$\{\Sigma^nR,n\in\zz\}$ an $\aleph_0$-compact generating set for
$D(R)$, as in Definition~8.1.6 of \cite{Neeman99}.  But then
Theorem~8.3.3 loc.  cit.  tells us that
$$\langle\{\Sigma^nR,n\in\zz\}\rangle~=~D(R)~,$$ 
that is any triangulated subcategory of $D(R)$, closed 
under coproducts and containing $R\in D(R)$, is all of $D(R)$.
\eprf

\dfn{D3.5}
Let $\cc$ be a triangulated category satisfying [TR5].
The subcategory $\cc^c\subset\cc$ is defined to be the
full subcategory of all compact objects in $\cc$.
\edfn

\rmk{R3.6}
The subcategory $\cc^c\subset\cc$ is always a thick subcategory.
\ermk

\exm{E3.7}
If $R$ is an associative ring and $\cb=D(R)$ is its unbounded derived
category, then $\cb^c$ turns out to be the full subcategory
of all objects isomorphic in $D(R)=\cb$ to bounded complexes
of f.g. projective $R$-modules. In the notation
of the representation theorists, what we have been denoting
$D(R)$ would be written $\cb=D(R\text{-Mod})$, and $\cb^c$ is
$D^b(R-\text{proj})$. The proof that $D^b(R-\text{proj})$
is precisely $\cb^c$ goes as follows. It is clear that 
$D^b(R-\text{proj})$ is a triangulated subcategory of $\cb$,
and Proposition~3.4 in \cite{Bokstedt-Neeman93} proves that it
is also closed under direct summands. That is,
$D^b(R-\text{proj})$ is a thick subcategory of $\cb$.
Example~1.13 of \cite{Neeman92A} shows that every object in
$D^b(R-\text{proj})$ is compact; that is,
$$
D^b(R-\text{proj})\subset\cb^c~.
$$
But $R$ is an object of $D^b(R-\text{proj})$, and by
Lemma~\ref{L3.1.7} any triangulated category of $\cb$
closed under coproducts and containing $R$ is all
of $\cb$. The fact that the inclusion
$$
D^b(R-\text{proj})\subset\cb^c
$$
is an equality now follows from Lemma~4.4.5 of \cite{Neeman99},
more precisely from the case where $\alpha=\beta=\aleph_0$
in that Lemma.

In the rest of this article, to keep the notation as simple
as possible, we write $\cb=D(R)$ for $D(R\text{-Mod})$,
and $\cb^c=D^c(R)$ for $D^b(R-\text{proj})$.
\eexm

Next we state the main existence theorem for Bousfield 
localizations. The theorem gives sufficient conditions
for the existence of a right adjoint $G$ to the
functor $\pi:\cb\ri\cc=\cb/\ca$. It also tells us
how the subcategories of compact objects, that is
$\ca^c$, $\cb^c$ and $\cc^c$, are related to each 
other.

\thm{T3.8}
Let $\cb$ be a triangulated category satisfying [TR5].
Suppose there are two sets of 
compact objects $A,B\subset\cb^c$
so that 
\begin{roenumerate}
\item Every triangulated subcategory of $\cb$, closed under
coproducts containing $B\subset\cb^c$, must equal all of $\cb$.
\item
The smallest triangulated subcategory of $\cb$, closed
under coproducts and containing $A\subset\cb^c$, will be denoted $\ca$.
\end{roenumerate}
We have an inclusion $\ca\subset\cb$. Let $\cc=\cb/\ca$
be the quotient. The following is true
\sthm{T3.8.1}
A Bousfield localization functor exists for the pair 
$\ca\subset\cb$. If $G$ is right adjoint to the
natural projection $\pi:\cb\ri\cc=\cb/\ca$, then
$G$ respects coproducts.
\esthm
\sthm{T3.8.2}
The categories $\ca^c$, $\cb^c$ and $\cc^c$ are all essentially
small.
\esthm
\sthm{T3.8.3}
$\ca^c=\ca\cap\cb^c$. Furthermore, $\ca^c$ can also be
described as the smallest thick subcategory of $\cb^c$
containing $A$.
\esthm
\sthm{T3.8.4}
The natural map $\cb^c/\ca^c\ri\cc$ is fully faithful, and factors
through $\cc^c\subset\cc$. The functor
$$\pi^c~:~\cb^c/\ca^c \ri \cc^c$$
is almost an equivalence of categories. It is fully faithful by
the above, and every object of $\cc^c$ is a direct summand of
an object in the image of $\pi^c$.
\esthm
\ethm

\prf
For the proof of \ref{T3.8.1}, see for example Lemma~1.7,
Remark~1.8 and Proposition~1.9 in \cite{Neeman92A}.
The proof of \ref{T3.8.4} may be found in Theorem~2.1
loc. cit. 

For \ref{T3.8.3}, observe that the inclusion $\ca\cap\cb^c\subset
\ca^c$ is obvious. An object of $\ca$ which is compact as
an object in the larger $\cb$ must be compact in $\ca$. 
The reverse inclusion, $\ca^c\subset\ca\cap\cb^c$,
may be found in Lemma~2.2 of \cite{Neeman92A},
where it is also proved that $\ca^c$ is the smallest
thick subcategory containing $A$.

Finally, for \ref{T3.8.2} we need to show that 
$\ca^c$, $\cb^c$ and $\cc^c$ are all essentially small.
For $\ca^c$, this is obvious by \ref{T3.8.2}; after
all, $\ca^c$ is obtained from the set $A$ of objects
by completing with respect to triangles and direct summands.
For $\cb$, we deduce it from the previous case by choosing
$A=B$, and therefore $\ca=\cb$ and $\ca^c=\cb^c$.
Now we know that $\cb^c$ and $\ca^c$ are essentially small. Hence 
so is $\cb^c/\ca^c$, and $\cc^c$ is obtained from it
by splitting idempotents. Therefore $\cc^c$ is also 
essentially small.
\eprf

\exm{E3.9}
The case of most interest to us here is where $\cb=D(R)$ is
the derived category of the ring $R$; we remind the reader,
this means the derived category of all chain
complexes of left $R$-modules. By Lemma~\ref{L3.1.7}, we
know that we can let the set $B$ be $B=\{R\}$. Every
triangulated subcategory of $\cb$, closed under
coproducts and containing $B\subset\cb^c$, is all of $\cb$.

What Theorem~\ref{T3.8} tells us is that we may choose
any set of objects $A\subset\cb^c$, and all the statements
of the theorem hold.
Let us begin by specialising to the case
we wish to study.
\eexm

\ntn{N5.1}
Suppose that
we are given a ring $R$. Let $\cb=D(R)$ be the derived category
of chain complexes of left $R$-modules. Let $\s$ be a set of objects 
in $\cb^c=D^c(R)$ of the form
$$\dots \ri0\ri c^\ell\ri c^{\ell+1}\ri0\ri\dots $$
That is, there exists an integer $\ell$ so that the
only non-zero terms in the complex occur in dimensions
$\ell$ and $\ell+1$. And the $R$-modules $c^\ell$ and 
$c^{\ell+1}$ are both f.g. and projective.
Note that we shall freely confuse the set of 
complexes $\s$ with the set of maps
$$c^\ell\ri c^{\ell+1}~.$$
Both sets will be called $\s$. We henceforth apply Theorem~\ref{T3.8}
to $\cb=D(R)$ with $A=\s$ and $B=\{R\}$. Let $\ca=D(R,\sigma)$
be the smallest triangulated subcategory of $\cb$ which is closed
under coproducts and contains $A=\s \subset\cb^c$. We shall
write the full subcategory of $\cb^c=D^c(R)$ consisting of the compact
objects in $\ca$ as
$$\ca^c~=~D^c(R,\sigma)~.$$
With $\cc=\cb/\ca$ as in  Theorem~\ref{T3.8}, the functor $\pi:\cb\ri\cc$ has
a right adjoint $G:\cc\ri\cb$.
\entn

It is now time to recall another well-known fact about
Bousfield localization. We begin with a definition.

\dfn{D100.1}
Let $\cb$ be a triangulated category satisfying [TR5].
Let $\s$ be a set of objects in $\cb$. An object $X$
is called {\em $\s$-local} if, for every object $s\in\s$
and every integer $n\in\zz$, $\cb(s,\Sigma^n X)=0$.
\edfn

\rmk{R100.1.5}
The full subcategory of all $\s$-local objects is triangulated.
\ermk

\lem{L100.2}
Let the notation be as in \ref{N5.1}. That is, 
$\cb=D(R)$, $\s$ is a set of objects  in $\cb^c$,
and $\ca=D(R,\sigma)$ is the smallest triangulated subcategory of $\cb$
containing $\s$ and closed under coproducts. Suppose 
$x\in\cb$ is $\s$-local. Then the unit of adjunction
$$
\CD
\eta_x~:~x @>>> G\pi x
\endCD
$$
is an isomorphism.
\elem

\prf
Fix a $\s$-local object $x\in\cb$.
Consider the full subcategory $\ct\subset\cb$ of all objects $t\in\cb$
so that, for every $n\in\zz$,
$$
\cb(t,\Sigma^n x) ~=~\cb(\Sigma^{-n}t,x)~=~0~.
$$
Because $x$ is $\s$-local, $\ct$ must contain $\s$. But
from its definition, it is clear that $\ct$ is closed under
coproducts and triangles. Hence it must contain $\ca$.
We conclude that, for every object $a\in\ca$ and
any integer $n\in\zz$, $\cb(\T^n a,x)=0$.

This means that, in the notation of Definition~9.1.10
of \cite{Neeman99}, $x$ lies in $^\perp\ca$.  
By Corollary~9.1.9 loc. cit, the unit of adjunction
$$
\CD
\eta_x:x @>>> G\pi x
\endCD
$$
is an isomorphism.
\eprf

\rmk{t-structure} The next lemmas will make use of the standard
$t$-structure on $\cb=D(R)$.  The reader will find an excellent
exposition of this topic in Chapter 1 of \cite{BeiBerDel82}.  We recall
that the full subcategory $\cb^{\leq n}$ is defined by
$$\text{\rm Ob}(\cb^{\leq n})~=~\{X\in\text{\rm Ob}(\cb)\mid
H^r(X)=0\,\,\hbox{\rm for all}~r>n\}~.$$
\ermk

\lem{L6.0.1}
As in Notation~\ref{N5.1}, let $\cb=D(R)$, let $\s$ be a set of objects 
in $\cb^c=D^c(R)$ of the form
$$\dots \ri0\ri c^\ell\ri c^{\ell+1}\ri0\ri\dots $$
If $x$ is a $\s$-local object, then so are its
$t$-structure truncations $x^{\leq n}$ and $x^{\geq n}$.
\elem

\prf
Pick a $\s$-local object $x$ and an integer $n\in\zz$;
without loss of generality we may assume $n=0$.
Let us begin by proving that $x^{\geq 0}$ is $\s$-local.
Take any $s\in\s$, that is a chain complex
$$\dots \ri0\ri c^\ell\ri c^{\ell+1}\ri0\ri\dots $$
and any integer $m\in\zz$. By Definition~\ref{D100.1},
it suffices to show that any map $\Sigma^{-m}s\ri x^{\geq 0}$
vanishes. That is, we must prove that any map from
$$\dots \ri0\ri c^{\ell+m}\ri c^{\ell+m+1}\ri0\ri\dots $$
to $x^{\geq 0}$ must vanish in $\cb$.  If $\ell+m+1<0$, then there is
no problem.  We have $\Sigma^{-m}s\in\cb^{<0}$ while $x^{\geq
n}\in\cb^{\geq 0}$, and hence any map $\Sigma^{-m}s\ri x^{\geq 0}$ must
vanish.  Assume therefore that $\ell+m+1\geq 0$.

We have a triangle
$$
\CD
x^{<0} @>>> x @>>> x^{\geq 0} @>w>> \Sigma x^{<0}~.
\endCD
$$
The composite 
$$\CD
\Sigma^{-m}s @>>>  x^{\geq 0}@>w>> \Sigma x^{<0}
\endCD$$
is a map from a bounded above complex of projectives
$\Sigma^{-m}s$ to some object in $\cb=D(R)$, and hence
it is represented by a chain map. But the chain
complex $\Sigma^{-m}s$ lives in degrees $\ell+m$ and
$\ell+m+1$, both of which are $\geq-1$, while the complex
$\Sigma x^{<0}$ lies in $\cb^{\leq-2}$. Hence the map vanishes.
 From the triangle we deduce that the map
$\Sigma^{-m}s\ri x^{\geq 0}$ must factor as
$$
\CD
\Sigma^{-m}s @>>> x @>>> x^{\geq 0}~.
\endCD
$$
But $x$ is $\s$-local by hypothesis, and hence 
any map $\Sigma^{-m}s \ri x$ vanishes.

This proves that $x^{\geq0}$ is $\s$-local. In the triangle
 $$
\CD
x^{<0} @>>> x @>>> x^{\geq 0} @>w>> \Sigma x^{<0}~.
\endCD
$$
we now know that both $x$ and $x^{\geq0}$ are $\s$-local,
and Remark~\ref{R100.1.5} permits us to deduce that
$x^{<0}$ is also $\s$-local.
\eprf

\lem{L6.1}
As in Notation~\ref{N5.1}, let $\cb=D(R)$, let $\s$ be a set of objects 
in $\cb^c=D^c(R)$ of the form
$$\dots \ri0\ri c^\ell\ri c^{\ell+1}\ri0\ri\dots$$
and let $\ca=D(R,\sigma)$ be the smallest triangulated subcategory of $\cb$
containing 
$\sigma$ and closed under coproducts. We assert that if $x\in\cb^{\leq n}$, 
then so is $G\pi x$.
\elem

\prf
We may assume without loss that $n=0$. Pick any 
$x\in\cb^{\leq 0}$. For any $s\in\sigma$ and any
$n\in\zz$, we have by adjunction
$$\cb(\Sigma^{-m}s,G\pi x)~=~\cc(\pi\Sigma^{-m}s,\pi x)$$
which vanishes since $\pi\Sigma^{-m}s=0$. Hence
$G\pi x$ is $\s$-local, and by Lemma~\ref{L6.0.1}
so is ${\{G\pi x\}}^{\leq 0}$.

Now the unit of adjunction $\eta_x:x\ri G\pi x$
is a map from an object $x\in\cb^{\leq 0}$, and
therefore factors (uniquely) as
$$
\CD
x @>\alpha>> {\{G\pi x\}}^{\leq0} @>f>> G\pi x~.
\endCD
$$
Since ${\{G\pi x\}}^{\leq 0}$ is $\s$-local, Lemma~\ref{L100.2}
coupled with \ref{T3.1.5.3} tell us that the map $\alpha$
factors (uniquely) as
$$
\CD
x @>\eta_x>> G\pi x @>g>> {\{G\pi x\}}^{\leq0}~.
\endCD
$$
The composite
$$
\CD
x @>\eta_x>> G\pi x @>g>> {\{G\pi x\}}^{\leq0} @>f>> G\pi x
\endCD
$$
is $\eta_x$, and the uniqueness statement in
\ref{T3.1.5.3} says that $fg:G\pi x\ri G\pi x$ is the identity.
The identity factors through an object in $\cb^{\leq0}$.

But then the identity must induce the zero map on all 
$H^n(G\pi x)$ with $n>0$. That is, for $n>0$ we have 
$H^n(G\pi x)=0$. In other words, $G\pi x\in\cb^{\leq0}$.
\eprf

\rmk{RT} The {\it derived tensor product} is defined for any $R-R$
bimodule $S$ to be the triangulated, coproduct-preserving functor
$$D(R)\ri D(R)~;~X\mapsto S{^L\otimes_R}X~=~S\otimes_RP$$ with $P$ a
sufficiently nice projective $R$-module chain complex quasi-isomorphic
to $X$.  The existence of this functor may be found, for example, in
Theorem~2.14 of B\"okstedt and Neeman \cite{Bokstedt-Neeman93}. 
The following are sufficiently nice~:
\begin{itemize}
\item[(a)] bounded above chain complexes of projectives,
\item[(b)] coproducts of the above,
\item[(c)] mapping cones of chain maps between the above.
\end{itemize}
\indent In this article, we consider tensor products both in the
category of modules and in the derived category.  We try to be careful
to distinguish them in the notation.  \ermk

\pro{L8.1}
As in Notation~\ref{N5.1}, let $\cb=D(R)$, let $\s$ be a set of objects 
in $\cb^c=D^c(R)$ of the form
$$\dots \ri0\ri c^\ell\ri c^{\ell+1}\ri0\ri\dots $$
and let $\ca=\cd(R,\sigma) \subset \cb$ be the smallest triangulated
subcategory containing $\sigma$ and closed under coproducts.

By Lemma~\ref{L3.1.6.4}, we have a ring homomorphism
$H^0(\eta_R^{}):R\ri S$, with $S=H^0(G\pi R)$.
Let $g:R \ri T$ be a ring homomorphism such that for
every $s\in \s\subset\cb^c$ we have
$T{^L\otimes_R^{}}s=0$. Then $g$ factors uniquely through the ring homomorphism
$H^0(\eta_R^{}):R\ri S$, as below~:
$$\xymatrix@C+5pt@R+15pt{& S \ar@{.>}[dr]^{\exists ! \, f} & \\
R \ar[ru]^{H^0(\eta_R)} \ar[rr]^{g} && T}$$
\epro

\prf
For any object $X\in D(R)=\cb$, we have
$$
D(R)\big(X,T\big)~=~D(T)\big(T{^L\otimes_R^{}}X,T
\big)~.
$$
By hypothesis, $T{^L\otimes_R^{}}s$ vanishes for all $s\in \s$.
Since
$$D(R)\big(\Sigma^r s,T\big)~=~
D(T)\big(T{^L\otimes_R^{}}\Sigma^r s,T\big)~=~0$$
for all $r\in\zz$ and all $s\in \s$, it follows that $T\in D(R)$ is
$\s$-local, and by Lemma~\ref{L100.2} the unit of adjunction
$$\eta_{T}^{}:T\ri G\pi T$$ is an isomorphism.  But then
$$\begin{array}{rclcl}
\cb(R,T)&=&\cb(R,G\pi T)&\qquad&
\text{because $\eta_{T}^{}$ is an isomorphism}\\
 &=& \cb(G\pi R,G\pi T) &\qquad&
\text{by \ref{T3.1.5.3}.}
\end{array}$$
This concretely translates to saying that the map $R\ri T=G\pi T$
factors uniquely through $\eta_R^{}:R\ri G\pi R$. Applying
$H^0$ to this, we have a factorization
$$\CD
R @>H^0(\eta_R^{})>> H^0(G\pi R)=S @>f>> T~.
\endCD$$
To see that the factorization is unique note that, by
Lemma~\ref{L6.1},
we know that $G\pi R\in\cb^{\leq0}$. The homomorphism 
$R\ri S$ therefore factors, in the derived category $D(R)=\cb$,
as the composite
$$\CD
R@>\eta_R^{}>> G\pi R@>>> {\{G\pi R\}}^{\geq0}~=~H^0(G\pi R)~=~S~.
\endCD$$
Given any morphism of $R$-modules $f':S\ri T$
we can form the composite
$$
\CD
R@>\eta_R^{}>> G\pi R@>>> {\{G\pi R\}}^{\geq0}~=~H^0(G\pi R)
@>f'>> T~.
\endCD
$$
If this composite agrees with $g$, then by the uniqueness
of the factorization through $R\ri G\pi R$  above,
the composites
$$
\CD
G\pi R@>>> {\{G\pi R\}}^{\geq0}
@>f>> T \\
G\pi R@>>> {\{G\pi R\}}^{\geq0}
@>f'>> T
\endCD
$$
must agree.
But $T\in\cb^{\geq0}$.
Therefore any map $G\pi R\ri T$ factors uniquely 
through the $t$-structure truncation
$$
\CD
G\pi R@>>> {\{G\pi R\}}^{\geq0}
@>>> T~.
\endCD
$$
Hence $f=f'$.
We have proved that there exists a unique map $f$ of left $R$-modules
making the triangle commute
$$\xymatrix@C+5pt@R+15pt{& S \ar@{.>}[dr]^{\exists ! \, f} & \\
R \ar[ru]^{H^0(\eta_R)} \ar[rr]^{g} && T}$$
\indent Now consider the functor 
$$D(R) \ri D(T)~;~X \mapsto T\oti X~.$$ 
By the hypothesis on $T$, this functor
annihilates all complexes $s\in\s$. Hence it annihilates
the subcategory $\ca=D(R,\sigma)\subset\cb=D(R)$. The functor factors
as
$$
\CD
D(R) @>\pi>> \cc~=~D(R)/D(R,\sigma) @>>> D(T)
\endCD
$$
Hence we have ring homomorphisms
$$
\CD
\End_{D(R)}(R) @>\alpha>> \End_{\cc}(\pi R) @>\beta>> \End_{D(T)}(T)~.
\endCD
$$
In Lemma~\ref{L3.1.6.4}, we checked that $\alpha$ agrees with the
homomorphism $H^0(\eta_R^{}):R\ri H^0(G\pi R)=S$, while from
the definition of the functor $D(R)\ri D(T)$ by tensor products,
the composite $\beta\alpha$ is nothing other than the given map
$g:R\ri T$. It follows that there is a commutative diagram
of ring homomorphisms
$$\xymatrix@C+5pt@R+15pt{& S \ar@{.>}[dr]^{\exists \, f'} & \\
R \ar[ru]^{H^0(\eta_R)} \ar[rr]^{g} && T}$$
But we already know that $f'$ is unique, even
as a map of $R$-modules.
There is a unique
ring homomorphism rendering commutative the triangle.
\eprf

\lem{L100.3}
The ring $S=H^0(G\pi R)$ satisfies the property that,
for all $s:P \ri Q$ in $\s$ the induced $S$-module morphism
$1\otimes s:S\otimes_RP \ri S\otimes_RQ$ is an isomorphism.
\elem

\prf
By Lemma~\ref{L6.1}, $G\pi R\in \cb^{\leq0}$. Hence
$S=H^0(G\pi R)={\{G\pi R\}}^{\geq0}$. By Lemma~\ref{L6.0.1},
${\{G\pi R\}}^{\geq0}$ is $\s$-local.

But then, for all $\{c^\ell\ri c^{\ell+1}\}\in\s$, the map
$$
\CD
\Hom_R(c^{\ell+1},S) @>>> \Hom_R(c^\ell, S)
\endCD
$$
is an isomorphism. Applying the functor $\Hom_S(-,S)$ 
to this isomorphism, and recalling that
$$\Hom_S\big[\Hom_R(A,S)\,,\,S\big]~=~S\otimes_R^{} A~,$$
we deduce that
$S\otimes_R^{} c^\ell\ri S\otimes_R^{} c^{\ell+1}$ is an isomorphism.
\eprf

\thm{T100.4}
With the notation as above, the ring homomorphism $R \ri S=H^0(G\pi R)$
satisfies the universal property of $R\ri \s^{-1}R$. 
\ethm

\prf
Lemma~\ref{L100.3} and Proposition~\ref{L8.1}.
\eprf

\ntn{N100.5}
We have now proved that, with $R$ a ring, $\cb=D(R)$
its derived category, and $\s\subset\cb^c$, 
$\ca=D(R,\sigma)\subset\cb$
and $\cc=\cb/\ca$ as in Notation~\ref{N5.1}, the ring homomorphism
$\eta_R^{}:R\ri H^0(G\pi R)$ satisfies Cohn's universal
property.  From now on, we shall freely confuse
 $H^0(G\pi R)=\s^{-1}R$.
\entn

The object $G\pi R$ in $D(R)$ is isomorphic to a chain complex 
$E(R)$ of free $R$-modules bounded above by 0.
We shall prove that $\s^{-1}E(R)$ is isomorphic to $\s^{-1}R$.
In fact, Vogel~\cite{Vogel1982} gave a direct construction for $E(R)$
as the direct limit of successive mapping cones.
This construction is very reminiscent of Bousfield's 
original proof of the existence of a Bousfield
localization; see Bousfield's~\cite{Bousfield75}
and \cite{Bousfield79A}. By now, of course, there
are other proofs of the existence of $G$.

\rmk{Rsymdual}
The universal property of $\rs$ is self-dual in the following
sense. We are given a set $\s$ of morphisms of f.g. projective 
left $R$--modules. If $s:P\ri Q$ is an element of $\s$,
we can look at the dual map $s^*:Q^*\ri P^*$. Here, $X^*=\Hom_R^{}(X,R)$.
To say that
\[\CD
1\otimes s:S\otimes_R^{}P @>>> S\otimes_R^{}Q
\endCD\]
is an isomorphism is equivalent to saying that 
\[\CD
s^*\otimes 1:Q^*\otimes_R^{}S @>>> P^*\otimes_R^{}S
\endCD\]
is an isomorphism. In other words, we can proceed to
do the entire construction of $\rs$ in terms of right
$R$--modules, just by replacing $\s$ by the set $\s^*$
of maps $s^*:Q^*\ri P^*$. Every theorem we prove has
a dual version for right modules.
\ermk

\pro{P8.3}
Let $M$ be a right $R$--module so that, for all $s:P\ri Q$
in $\s$, $1\otimes_R^{} s:M\otimes_R^{} P\ri M\otimes_R^{} Q$
is an isomorphism. Then the 
map $1{^L\otimes_R^{}}\,\eta_R^{}:M\ri
 M\,{^L\otimes_R^{}}\,G\pi R$
is an isomorphism.
\epro

\prf
Form the triangle
$$
\CD
a @>>> R @>\eta_R^{}>> G\pi R @>>> \Sigma a~.
\endCD
$$
By \ref{T3.1.5.5}, $a\in\ca$.
Now the set of objects $x\in\cb=D(R)$ such that
$M\oti x$ vanishes is a triangulated category
containing $\s$, and closed under coproducts.
It must contain all of $\ca$. Therefore $M\,{^L\otimes_R^{}}\,a=0$.
The triangle
$$
\CD
M\,{^L\otimes_R^{}}\,a @>>> M 
@>M\,{^L\otimes_R^{}}\,\eta_R^{}>> 
M\,{^L\otimes_R^{}}\,G\pi R @>>> M\,{^L\otimes_R^{}}\,\Sigma a
\endCD
$$
tells us that $M \,{^L\otimes_R^{}}\,\eta_R^{}$ must be an isomorphism.
\eprf

\rmk{Rdually}
Dually as in Remark~\ref{Rsymdual}, given any left module
$M$ so that, for any $s^*:Q^*\ri P^*$ in $\s^*$,
the map $s^*\otimes 1$ is an isomorphism, then 
$\eta'_R\oti 1:M\ri
 G'\pi' R\oti M$ is an isomorphism. 
[Here, $G'$, $\pi'$ and $\eta'$ are the duals of $G$, $\pi$
and $\eta$ in the derived category of right $R$-modules].
Note that
\[
P^*\oti M~=~\Hom_R^{}(P,R)\otimes_R^{} M~=~\Hom_R^{}(P,M)~. 
\]
We deduce that $s^*\otimes 1_M^{}$ is an isomorphism
if and only if $\cb(\T^n s,M)=0$ for all $n\in\zz$.
The $M$'s for which $s^*\otimes 1_M^{}$ are isomorphisms
whenever $s\in\s$ are precisely the $\s$--local objects,
as in 
Definition~\ref{D100.1}. Summarising:
whenever
$M$ is a left $R$-module, which as an object of $\cb=D(R)$
happens to be $\s$--local,
then the dual of Proposition~\ref{P8.3}
asserts that $\eta'_R\oti 1:M\ri
 G'\pi' R\oti M$ is an isomorphism.  
\ermk

\pro{P8.3.4}
Let $M$ be any right $R$-module so that, for all $s:P\ri Q$
in $\s$, $1\otimes_R^{} s:M\otimes_R^{} P\ri M\otimes_R^{} Q$
is an isomorphism. Then the map $1\otimes H^0(\eta_R^{}):M\ri
M\otimes_R^{}\br$ is an isomorphism. Furthermore,
$$\text{\rm Tor}^R_1(M,\s^{-1}R)~=~0~.$$
\epro

\prf
There is a spectral sequence computing the cohomology 
of $M{^L\otimes_R^{}}G\pi R$. The $E_2$  term is
$$E_2^{i,j}~=~\text{\rm Tor}^R_{-i}\big(M,H^j(G\pi R)\big)~.$$
By Lemma~\ref{L6.0.1}, $H^j(G\pi R)=0$ for all $j>0$.
And $\text{\rm Tor}^R_{-i}$ vanishes for $i>0$. This spectral
sequence is concentrated in negative degrees.
Because all the differentials in and out of the following
terms vanish, we conclude that,
in the spectral sequence, $E_2^{0,0}=E_\infty^{0,0}$ and
$E_2^{-1,0}=E_\infty^{-1,0}$.

By Proposition~\ref{P8.3},
we have that $M{^L\otimes_R^{}}G\pi R=M$. The above spectral
sequence converges to $H^{i+j}(M)$, which is $M$ if $i+j=0$
and zero otherwise. We immediately conclude that 
$$E_2^{0,0}~=~M\otimes_R^{}\s^{-1}R~=~M$$
and
$$E_2^{-1,0}~=~\text{\rm Tor}^R_{1}(M,\s^{-1}R)~=~0~.$$
\eprf

\cor{C8.4} 
The natural multiplication map $\s^{-1}R\otimes_R^{}\s^{-1}R\ri \s^{-1}R$ is
an isomorphism, and 
$$\text{\rm Tor}^R_1(\s^{-1}R,\s^{-1}R)~=~0~.$$
\ecor
\prf
In Proposition~\ref{P8.3.4}, put $M=\rs$. We immediately deduce
that $\text{\rm Tor}^R_1(\s^{-1}R,\s^{-1}R)~=~0~.$ We also have 
that the map
$$\CD
1\otimes H^0(\eta_R^{}):\rs @>>>
\br\otimes_R^{}\br\endCD$$
is an isomorphism. But the composite
$$\CD
\rs @>1\otimes H^0(\eta_R^{})>>
\br\otimes_R^{}\br @>\text{multiplication}>> \rs
\endCD$$
is clearly the identity, and hence the multiplication map must
be the two--sided inverse of the invertible map
$1\otimes H^0(\eta_R^{})$.
\eprf

\rmk{Rschof} The result $\text{\rm Tor}^R_1(\s^{-1}R,\s^{-1}R)=0$ is due to
Schofield \cite{Schofield}, p.58.\ermk

\cor{C8.4.1}
Suppose $M$ is a left $R$-module, which as an object of 
$\cb=D(R)$ is $\s$--local (see Definition~\ref{D100.1}).
Then the $R$-module structure of $M$ extends, uniquely,
to a $\rs$-module structure.
\ecor

\prf
The uniqueness of the extension is clear: to say that
$M$ is a left $R$-module is to give a homomorphism
\[
R\ri\Hom_\zz^{}(M,M),
\]
and the universal property of the noncommutative localization 
tells us that the factorization of this ring homomorphism
through $\rs$ is certainly unique, if it exists.

It remains to prove existence. By the dual of 
Proposition~\ref{P8.3.4}, the map
\[
\CD
M=R\otimes_R^{}M @>H^0(\eta_R^{})\otimes_R^{} 1>> \br\otimes_R^{}M
\endCD
\]
is an isomorphism. Thus $M$ is isomorphic, as a left $R$--module,
to $\br\otimes_R^{}M$, and $\br\otimes_R^{}M$ is certainly a module
over
$\rs$.
\eprf

In the
remainder of this section, we want to relate
$H^n(G\pi R)$ for $n<0$ to
$\text{Tor}^R_*(\sigma^{-1}R,\sigma^{-1}R)$.
First we prove a lemma.

\lem{L8.5}
Suppose $R \ri S$ is a ring homomorphism such that the multiplication map
$S\otimes_R^{}S\ri S$ is an isomorphism.  Suppose also that for some $n \geq 1$
$$\text{\rm Tor}^R_i(S,S)~=~0~(1\leq i\leq n)~.$$
Then for every $S$-module $M$, we have

\sthm{L8.5.1}
The multiplication map $S\otimes_R^{}M\ri M$ is an isomorphism.
\esthm
\sthm{L8.5.2}
$\text{\rm Tor}^R_i(S,M)=0$ for all $1\leq i\leq n$.
\esthm
\elem

\prf
Choose a resolution of $M$ by free $S$-modules
$$
\CD
\cdots @>>> Q^{-2} @>>> Q^{-1} @>>> Q^0 @>>> M @>>> 0~,
\endCD
$$
and a resolution of $S$ by free $R$-modules
$$
\CD
\cdots @>>> P^{-2} @>>> P^{-1} @>>> P^0 @>>> S @>>> 0~.
\endCD
$$
The tensor product $P\otimes Q$ gives a double
complex whose cohomology computes $\text{\rm Tor}^R_{-i-j}(S,M)$.
But there is a spectral sequence for it, whose $E_1$ term is
$$E_1^{i,j}~=~\text{\rm Tor}^R_{-j}(S,Q^i)$$ 
Now $E_1^{i,0}=S\otimes_R^{}Q^i=Q^i$, since $Q^i$ is free
and, by hypothesis,  $S\otimes_R^{}S\ri S$ is an
isomorphism. In $E_2$, we have
$$E_2^{i,0}~=~\left\{
\begin{array}{ccl}
M &\qquad & \text{if $i=0$}\\
0 &\qquad & \text{otherwise.}
\end{array}
\right.$$
But by hypothesis, we also have $\text{\rm Tor}^R_{-j}(S,S)=0$, for all
$1\leq -j\leq n$, and since $Q^i$ are free, this gives $\text{\rm
Tor}^R_{-j}(S,Q^i)=0$, for all $i$ and for all $1\leq -j\leq n$.  In
other words, $E_1^{i,j}=0$ if $1\leq -j\leq n$, and hence $E_2^{i,j}=0$
if either $j=0$, $i\neq0$, or if $1\leq -j\leq n$.  The assertions of
the Lemma immediately
follow.
\eprf

\cor{C8.6}
Let the notation be as in Notation~\ref{N5.1}.
Suppose $\text{\rm Tor}^R_{i}(\s^{-1}R,\s^{-1}R)=0$, for all $1\leq i\leq n$.
Then for all $1\leq i\leq n-1$ we have $H^{-i}(G\pi R)=0$,
and
$$\text{\rm Tor}^R_{n+1}(\s^{-1}R,\s^{-1}R)~=~H^{-n}(G\pi R)~.$$
\ecor

\prf
The proof is a slightly more sophisticated computation 
with the same spectral
sequence we saw in Proposition~\ref{P8.3.4}. Recall that we have
a spectral sequence whose $E_2$ term is
$$
E_2^{i,j}~=~\text{\rm Tor}^R_{-i}\big(\s^{-1}R,H^j(G\pi R)\big)~,
$$
which converges to $H^{i+j}(\s^{-1}R)$. Now we know that $G\pi R$
is $\s$--local. By Lemma~\ref{L6.0.1} so
are its $t$--structure truncations
$H^j(G\pi R)={\left\{{[G\pi R]}^{\leq j}\right\}}^{\geq j}$.
Corollary~\ref{C8.4.1} now tells us that $H^j(G\pi R)$ are
all left $\rs$--modules.
Lemma~\ref{L8.5} now applies, and we deduce that
if $1\leq -i\leq n$ then $E_2^{i,j}=0$. This forces
the differential
$$E_2^{-i-1,0}\ri E_2^{0,-i}$$
to be an isomorphism, for all $1\leq i\leq n$. For 
$1\leq i\leq n-1$ we read off that $H^{-i}(G\pi R)=0$.
For $i=n$, we deduce that
$$\text{\rm Tor}^R_{n+1}(\s^{-1}R,\s^{-1}R)~=~H^{-n}(G\pi R)~.$$
\eprf

\section{A bound on the length of complexes in $D^c(R,\sigma)$}
\label{Sbound}

\ntn{NSbound.1}
Our notation stays as in Notation~\ref{N5.1}. $D(R)=\cb$
is the derived category of a ring $R$, and we are given a
set $\s$ of maps of f.g. projective $R$-modules. The 
category $\ca=D(R,\sigma)$ is the smallest triangulated subcategory of $\cb$
containing $\s$ and closed under coproducts.
Let $\cc=\cb/\ca$. Let $\pi:\cb\ri\cc$ be the projection,
$G:\cc\ri\cb$ the fully faithful right adjoint. Identify 
$R \ri \s^{-1}R$ with $R \ri H^0(G\pi R)$. 
\entn

For any f.g. projective $R$-modules $M,N$ a morphism $\pi M \to
\pi N$ in $\cc$ is an equivalence class of diagrams
$$\xymatrix@C+15pt{M  \ar[r]^-{\beta} & Y & N \ar[l]_-{\alpha}}$$
with $\alpha$ a morphism  in $\cb$ which becomes an isomorphism in
$\cc$. Later in the article, we will need to have bounds
on the length of $Y$.
In this section we carry out the
preparatory technical work.

\dfn{DSbound.2}
The full subcategory of all objects in $\cb=D(R)$ which vanish outside
the range $[m,n]$ will be denoted $\cb[m,n]$. 
We allow $m$ or $n$ to be infinite; the categories
$\cb[m.\infty)$ and  $\cb(-\infty,n]$
have the obvious definitions.
\edfn

\rmk{RSbound.3}
The reader should note that the categories 
$\cb[n,\infty)$ and  $\cb(-\infty,n]$ should not be
confused with $\cb^{\geq n}$ and $\cb^{\leq n}$.
It is true that every object in $\cb^{\leq n}$ is
{\em isomorphic in $\cb$} to a chain complex
\[
\cdots \ri X^m\ri X^{m+1}\ri\cdots
\ri X^{n-1}\ri X^n\ri0\ri0\ri\cdots
\]
An isomorphism in $\cb=D(R)$ is after all just a 
homology isomorphism.
For any object in $\cb^{\leq n}$, there is an
object in $\cb(-\infty,n]$ homology isomorphic to it. But for
once we want to have a name for the complexes which 
are actually supported on the interval $[m,n]$, not 
just isomorphic in $\cb$ to such objects.
\ermk

\dfn{Dbound.4}
The category $\cs$ will be the smallest full
subcategory of $\cb$ such that
\sthm{Dbound.4.1}
Every suspension of every object in $\s$ lies in 
$\cs$. That is, $\cs$ contains all the complexes
\[
\cdots\ri0\ri c^{\ell+n}\ri c^{\ell+n+1}\ri0\ri\cdots
\]
\esthm
\sthm{Dbound.4.2}
Given any chain map of objects in $\cs$
\[
\CD
\cdots @>\partial>> X^{i-1} @>\partial>> X^i 
@>\partial>>X^{i+1} @>\partial>>\cdots \\
@. @Vf_{i-1}^{}VV @Vf_{i}^{}VV @Vf_{i+1}^{}VV @. \\
\cdots @>\partial>> Y^{i-1} @>\partial>> Y^i 
@>\partial>>Y^{i+1} @>\partial>>\cdots 
\endCD
\]
then the mapping cone
\[
\CD
\cdots @>>> X^{i}\oplus Y^{i-1}
@>{\left(
\begin{array}{rl}
-\partial & 0 \\
f_{i}^{} & \partial
\end{array}
\right)}>> X^{i+1}\oplus Y^{i}
@>{\left(
\begin{array}{rl}
-\partial\,\,\,\, & 0 \\
f_{i+1}^{} & \partial
\end{array}
\right)}>>
X^{i+2}\oplus Y^{i+1} @>>>\cdots
\endCD
\]
also lies in $\cs$.
\esthm
\nin
As in Remark~\ref{RSbound.3}, we mean equality
of chain complexes, not homotopy equivalence.
\edfn

\dfn{Dbound.5}
The subcategories $\cs[m,n]$ are defined as the intersection
\[
\cs[m,n]~=~\cs\cap\cb[m,n]~.
\]
As in Definition~\ref{DSbound.2}, we allow $m$ and $n$ to be infinite.
\edfn

\lem{Lbound.6}
Suppose $n\in\zz$ is an integer.
Then every object $Z\in\cs$ can be expressed as
a mapping cone on a chain map $Z_1\ri Z_2$, as below
\[
\CD
\cdots @>\partial>> Z_1^{n-1} @>\partial>> Z_1^n @>\partial>> Z_1^{n+1} 
@>>>0 @>>>\cdots \\
@. @VVV @Vf_{n}^{}VV @Vf_{n+1}^{}VV @VVV @. \\
\cdots @>>> 0 @>>> Z_2^n
@>\partial>>Z_2^{n+1} @>\partial>>Z_2^{n+2} @>\partial>>\cdots 
\endCD
\]
that is, $Z_1\in\cs(-\infty,n+1]$ and
$Z_2\in\cs[n,\infty)$.
\elem

\prf
Let $\ct$ be the full subcategory of $\cs$ containing the 
objects for which the assertion of the lemma holds. That
is, an object $Z\in\cs$ belongs
to $\ct$ if and only if, for every $n\in\zz$,
there exist $Z_1\in\cs(-\infty,n+1]$ and
$Z_2\in\cs[n,\infty)$
and a chain map $Z_1\ri Z_2$ so that $Z$ is equal to 
the mapping cone. It suffices to prove that $\ct=\cs$,
for which we need only show that any suspension of
an object of $\s$ lies in $\ct$, and that mapping cones 
on maps in $\ct$ lie in $\ct$.

Assume therefore that we are given a complex $s$ below
  \[
\cdots\ri0\ri c^{\ell}\ri c^{\ell+1} \ri0 \ri\cdots
\]
which is some suspension of an object in $\s$.
Choose any $n\in\zz$. If $n\leq\ell$, then $s\in\cs[n,\infty)$,
and $s$ is the mapping cone of the chain map $0\ri s$.
If $n\geq \ell+1$, then $\Sigma^{-1}s\in\cs(-\infty,n+1]$
and $s$ is isomorphic to the mapping cone on the chain
map $\Sigma^{-1}s\ri0$. Either way, $s\in\ct$.

Next suppose we are given two object $X$ and $Y$ in $\ct$,
and a chain map $f:X\ri Y$. Let $Z$ be the mapping cone
of $f$. We need to show that $Z$
is in $\ct$. For every integer $n\in\zz$, we need to
express $Z$ as a mapping cone on a map of objects
$Z_1\ri Z_2$, with $Z_1\in\cs(-\infty,n+1]$ and
$Z_2\in\cs[n,\infty)$. Without loss of generality, assume
$n=0$.

Because $X\in\ct$, we may express it as a mapping cone on
a map $X_1\ri X_2$, with $X_1\in\cs(-\infty,2]$ and
$X_2\in\cs[1,\infty)$.  Because $Y\in\ct$, we may express it as the
mapping cone on a map $Y_1\ri Y_2$, with $Y_1\in\cs(-\infty,1]$ and
$Y_2\in\cs[0,\infty)$. We have a diagram, where the rows
are short exact sequences of chain complexes
\[
\CD
0 @>>> X_2 @>>> X  @>>> \Sigma X_1 @>>> 0 \\
@.     @.     @VfVV      @.            @. \\
0 @>>> Y_2 @>>> Y  @>>> \Sigma Y_1 @>>> 0
\endCD
\]
The composite
\[
\CD
X_2 @>>> X  @.  \\
  @.     @VfVV      @.    \\
@. Y  @>>> \Sigma Y_1 
\endCD
\]
is a chain map from $X_2\in\cs[1,\infty)$ to
$\Sigma Y_1\in\cs(-\infty,0]$, and therefore must vanish.
It follows that we may complete to a commutative 
diagram of chain complexes
\[
\CD
0 @>>> X_2 @>>> X  @>>> \Sigma X_1 @>>> 0 \\
@.     @Vf_1VV     @VfVV      @V\Sigma f_2VV           @. \\
0 @>>> Y_2 @>>> Y  @>>> \Sigma Y_1 @>>> 0
\endCD
\]
Let $Z_1$ be the mapping cone on $f_1:X_1\ri Y_1$, and
let $Z_2$ be the mapping cone on $f_2:X_2\ri Y_2$. 
Then $Z_1\in\cs(-\infty,1]$ while $Z_2\in\cs[0,\infty)$.
Furthermore, $Z$, which is the mapping cone on $f:X\ri Y$,
can also be expressed as a mapping cone on
a map $Z_1\ri Z_2$.
\eprf

\lem{Lbound.6.5}
Let $\tilde\cs$ be the full subcategory of all objects
in $\cb$ isomorphic to objects in $\cs$. That is, any
object of $\cb=D(R)$ isomorphic to a chain complex 
in $\cs$ lies in $\tilde\cs$. The subcategory 
$\tilde\cs\subset\cb$ is triangulated.
\elem

\prf
The point is that the objects of $\cs$ are bounded chain
complexes of projectives. Let $f:X\ri Y$ be a morphism in $D(R)$,
between objects in $\cs$. Because $X$ is a bounded-above
complex of projectives, there is a chain map representing
the morphism. The mapping cone on this chain map
completes $f:X\ri Y$ to a triangle, and lies in $\cs$.
Up to isomorphism in $\cb=D(R)$, the third edge is unique.
Therefore in any triangle
\[
\CD
X @>f>> Y @>>> Z @>>> \Sigma X
\endCD
\]
we have $Z\in\tilde\cs$.
\eprf

\lem{Lbound.7}
The category $\cs$ is contained in $\ca^c$. Furthermore,
the inclusion is nearly an equality.
Every object in $\ca^c$ is
a direct summand
of an object isomorphic in $D(R)$ to an object in $\cs$.
\elem

\prf
The inclusion $\cs\subset\ca^c$ is easy. The category $\ca^c$
contains $\s$ and is closed under mapping cones, and $\cs$ is
the smallest such.

Next observe that, by \ref{T3.8.3}, the category $\ca^c$ 
is the smallest thick subcategory of $\cb$ containing
$\s$, and hence $\ca^c$ is the smallest thick
subcategory containing the triangulated subcategory $\tilde\cs$
of Lemma~\ref{Lbound.6.5}.
Now Corollary~4.5.12 of \cite{Neeman99} tells us that, for every
object $X\in\ca^c$, the object $X\oplus \Sigma X$ lies in
$\tilde\cs$. In particular
$X\in\ca^c$ is a direct summand of an object isomorphic
in $\cb$ to an object in $\cs$.
\eprf

\lem{Lbound7.5}
The natural map $\cb^c/{\tilde\cs}^c\ri \cc$ is fully faithful.
\elem

\prf
By \ref{T3.8.4}, the functor
\[
\cb^c/\ca^c\ri\cc^c\subset\cc
\]
is fully faithful. 
Now by Lemma~\ref{Lbound.7}, 
the thick closure of $\tilde\cs$ is all of $\ca^c$,
and hence $\cb^c/\ca^c=\cb^c/{\tilde\cs}^c$.
\eprf

\rmk{Rbound.200}
One way we shall use the results of this section is as
follows. Suppose $x\ri y$ is a morphism in $\cb^c$, which
maps to zero in $\cc=\cb/\ca$. By Lemma~\ref{Lbound7.5},
the map vanishes already in $\cb^c/\tilde\cs$.
By Lemma~2.1.16 of \cite{Neeman99}, $x\ri y$ must factor
as $x\ri s\ri y$ with $s\in\cs$. The results of this 
section will enable us to replace $s$ by shorter complexes.
\ermk

In section~\ref{Ltheory} we will need the following
result

\pro{Pbound.8}
Let $M$ and $N$ be any f.g. projective $R$--modules,
which we view as objects in the derived category $\cb=D(R)$,
concentrated in degree $0$. Then any map in $\cc(\pi M,\pi N)$ can be
represented as $\alpha^{-1}\beta$, for some $\alpha$, $\beta$ as below
$$M \xrightarrow[]{\beta} Y \xleftarrow[]{\alpha} N$$
The map $\alpha:N\ri Y$ fits in a triangle
\[
\CD
X @>>> N @>\alpha>> Y @>>> \Sigma X
\endCD
\]
and $X$ may be chosen to lie in $\cs[0,1]$.
\epro

\prf
Both $M$ and $N$ are assumed to be objects
of $\cb^c$, and by Lemma~\ref{Lbound7.5}
the map 
$\cb/{\tilde\cs}^c\ri\cc$ if fully faithful. Hence
\[
\cc(\pi M,\pi N)=
\{\cb^c/{\tilde\cs}^c\}(M,N).
\]
Any morphism may be represented by a diagram
$$M \xrightarrow[]{\beta} Y \xleftarrow[]{\alpha} N$$
so that in the triangle
\[
\CD
X @>>> N @>\alpha>> Y @>>> \Sigma X
\endCD
\]
$X$ may be chosen to lie in $\cs$.

By Lemma~\ref{Lbound.6}, there exists a triangle in
$\ca^c$
\[
\CD
X' @>>> X @>>> X'' @>>> \Sigma X'
\endCD
\]
with $X'\in\cs[1,\infty)$ and $X''\in\cs(-\infty,1]$.
The composite $X'\ri X\ri N$ is a map from
$X' \in\cs[1,\infty)$ to $N\in\cs[0,0]$, which must
vanish. Hence we have that $X\ri N$ factors as
$X\ri X''\ri N$. We complete to a morphism of triangles
\[
\CD
X @>>> N @>\alpha>> Y @>>> \Sigma X \\
@VVV  @V1VV     @V\gamma VV         @VVV \\
X'' @>>> N @>\gamma\alpha>> Y'' @>>> \Sigma X''
\endCD
\]
and another representative of our morphism is the diagram
$$M \xrightarrow[]{\gamma\beta} Y'' \xleftarrow[]{\gamma\alpha} N$$
We may, on replacing $Y$ by $Y''$, assume $X\in\cs(-\infty,1]$.

Applying Lemma~\ref{Lbound.6}
again, we have that any 
$X\in\cs(-\infty,1]$ admits a triangle
\[
\CD
X' @>>> X @>>> X'' @>>> \Sigma X'
\endCD
\]
with $X'\in\cs[0,1]$ and  $X''\in\cs(-\infty,0]$. Form the
octahedron
\[
\CD
X' @>>> N @>\alpha'>> Y' @>>> \Sigma X' \\
@VVV  @V1VV     @V\gamma VV         @VVV \\
X @>>> N @>\alpha>> Y @>>> \Sigma X \\
@. @. @VVV @VVV \\
@. @.   \Sigma X'' @>1>> \Sigma X'' 
\endCD
\]
The composite $M\ri Y\ri \Sigma X''$ is a map
from $M\in\cs[0,0]$ to $\Sigma X''\in\cs(\infty,-1]$,
and must vanish. The map $\beta:M\ri Y$ therefore factors as
$M\stackrel{\beta'}\ri Y'
\stackrel\gamma\ri Y$, and our morphism in $\cc$ has a
representative
$$M \xrightarrow[]{\beta'} Y' \xleftarrow[]{\alpha'} N$$
so that in the triangle
\[
\CD
X' @>>> N @>\alpha'>> Y' @>>> \Sigma X'
\endCD
\]
$X'$ may be chosen to lie in $\cs[0,1]$.
\eprf

\section{The functor $T:(D(R)/D(R,\sigma))^c\to D^c(\sigma^{-1}R)$}
\label{S5}

Let $R\ri S$ be a ring homomorphism.  There is a triangulated
functor $D(R)\ri D(S)$, taking $X\in D(R)$ to $S\oti X$.  In this
section, we shall study this functor in the case where $R$ is a ring,
and $S=\rs$ is the noncommutative localization of $R$ inverting $\sigma$.

\pro{P5.1}
The functor 
$$\cb~=~D(R) \ri \cd~=~D(\sigma^{-1}R)~;~X \mapsto \br\oti X$$
factors uniquely (up to canonical natural isomorphism)
through $\cb\ri\cc=D(R)/D(R,\sigma)$. The unique factorization will 
be written as
$$\CD 
\cb @>\pi>> \cc @>T>> \cd~.
\endCD
$$
The functor $T:\cc\ri \cd$ takes $\cc^c\subset\cc$
to $\cd^c=D^c(\rs)\subset \cd=D(\rs)$.
\epro

\prf
Let $c^\ell\ri c^{\ell+1}$ be a map in $\s$. Tensoring
with $\rs$ takes it to an isomorphism. Hence tensoring
with $\rs$ takes the chain complex
$$
\cdots\ri 0 \ri c^\ell\ri c^{\ell+1}\ri0 \ri\cdots
$$
to an acyclic complex.  Therefore the functor 
$X\mapsto\br\oti X:\cb \ri \cd$ kills all
the objects in $\s$.  Since derived tensor product preserves triangles
and coproducts, the subcategory of $\cb$ annihilated by 
$X\mapsto\br\oti X$ must be
closed under triangles and coproducts, and therefore contains all of
$\ca=D(R,\sigma)$.  By the universal property of the Verdier quotient
$\cc=\cb/\ca$, there is a unique factorization
$$\CD
\cb @>\pi>> \cc @>T>> \cd~.
\endCD$$
It remains to show that $T$ takes $\cc^c\subset\cc$ to $\cd^c\subset \cd$.\\
\indent It is clear that the map $T\pi:\cb\ri \cd$ takes a bounded
complex of f.g. projective $R$-modules to a bounded complex of 
f.g. projective
$\rs$-modules; the map just tensors with $\rs$.  In other words, the
functor $T\pi$ obviously takes $\cb^c$ to $\cd^c \subset \cd$.  
By \ref{T3.8.4}, every object in $\cc^c$ is a direct summand of an object
in the image of $\pi : \cb^c \ri \cc^c$.  Therefore the functor $T$
takes any object in $\cc^c$ to a direct summand of an object in
$\cd^c$.  
But by Proposition~3.4 of \cite{Bokstedt-Neeman93}, any
direct summand of an object in $\cd^c$ lies in $\cd^c$.
\eprf

\pro{P5.2}
For any two projective $R$-modules $P$ and $Q$, one has
$$\cc(\pi P,\pi Q)~=~\cd(T\pi P,T\pi Q)~=~\Hom_{\rs}^{}(\s^{-1}P,\s^{-1}Q)~.$$
\epro
\prf
The identity $\cd(T\pi P,T\pi Q)~=~\Hom_{\rs}^{}(\s^{-1}P,\s^{-1}Q)$
is just by definition. We have defined
$T\pi P=\br\times_R^{}P=\s^{-1}P$.

There is a natural map, induced by the functor $T$,
\[
\CD
\cc(\pi P,\pi Q) @>>> \cd(T\pi P,T\pi Q).
\endCD
\]
We need to prove it an isomorphism.
The case where $P=Q=R$ is easy; we have
$$
\begin{array}{rclcl}
\cc(\pi R, \pi R)
&~=~&\cb(R, G\pi R)&\qquad&\text{by adjunction}\\[1ex]
&~=~&H^0(G\pi Q)\\[1ex]
&~=~&\rs&\qquad&\text{by Theorem~\ref{T100.4}}\\[1ex]
&~=~&\cd(T\pi R,T\pi R)~. & &
\end{array}
$$
But the collection of all $P$ and $Q$ for which
the map $T:\cc(\pi P,\pi Q)\ri \cd(T\pi P,T\pi Q)$
is an isomorphism is clearly closed under direct
sums and direct summands, and hence contains all
projective modules.
\eprf

In Section~\ref{S8}, we shall need to know that the
only object of $\cc^c$ annihilated by $T$ is the zero object.
The next two propositions prove this.

\pro{P5.3}
For any object $X$ in
$\cb^c=D^c(R)$, we have the implication
$$
\big\{\br{^L\otimes_R^{}}X=0\big\}\,\,\Longrightarrow\,\,
\{X\in \ca^c\}.
$$
\epro

\prf
Take any $X\in\cb^c$.
Since $X\in\cb^c=D^c(R)$, we know that $X$ is isomorphic
to
a bounded complex of f.g. projective $R$-modules.
Up to suspension, $X$ may be written as a complex
$$
\ri 0\ri X^0 \ri X^1 \ri \cdots \ri
X^{n-1} \ri X^n \ri 0 \ri
$$
If $\br{^L\otimes_R^{}}X=0$, then the complex
$$
\ri 0\ri \sigma^{-1}X^0 \ri  \cdots \ri
 \sigma^{-1}X^n \ri 0 \ri
$$
must be contractible. There are maps $\sigma^{-1}X^i\ri 
\sigma^{-1}X^{i-1}$
so that, for each $i$, the sum of the two composites
$$
\CD
\sigma^{-1}X^i @>>> \sigma^{-1}X^{i+1} \\
@VVV @VVV \\
\sigma^{-1}X^{i-1} @>>> \sigma^{-1}X^i
\endCD
$$
is the identity on
$\sigma^{-1}X^i$. 
By Proposition~\ref{P5.2}, the contracting homotopy may
be lifted to the complex
$$
\ri 0\ri \pi X^0 \stackrel\partial\ri \pi X^1 
\stackrel\partial\ri \cdots \stackrel\partial\ri
\pi X^{n-1} \stackrel\partial\ri \pi X^n \ri 0 \ri
$$
For each $i$ there are maps $\pi X^i\ri \pi X^{i-1}$, so that the
two composites
$$
\CD
\pi X^i @>\partial>> \pi X^{i+1} \\
@V{D}VV @VV{D}V \\
\pi X^{i-1} @>\partial>> \pi X^i
\endCD
$$
add to the identity on $\pi X^i$.

Now let $Y^i\in\cb^c$ be the complex
$$
\ri 0\ri X^0 \ri X^1 \ri \cdots \ri
X^{i-1} \ri X^i \ri 0 \ri
$$
For each $i$, there is a triangle
$$
\CD
\Sigma^{-i-1}X^{i+1} @>>> Y^{i+1} @>>> Y^i @>>> \Sigma^{-i}X^{i+1}.
\endCD
$$
The functor $\pi$ is triangulated, and hence for each $i$
we deduce a triangle
$$
\CD
\Sigma^{-i-1}\pi X^{i+1} @>>> \pi Y^{i+1} @>>> \pi Y^i 
@>\rho_i^{}>> \Sigma^{-i}\pi X^{i+1}~.
\endCD
$$
We shall prove, by induction on $i$, that 
\begin{roenumerate}
\item The map 
$$
\CD
\pi Y^i 
@>\rho_i^{}>> \Sigma^{-i}\pi X^{i+1}
\endCD
$$
is a split monomorphism in $\cb=D(R)$.
\item
For each $i$ we shall produce an explicit splitting; that is,
we shall produce a map
$$
\CD
\Sigma^{-i}\pi X^{i+1}
@>\theta_i^{}>> \pi Y^i
\endCD
$$
so that $\theta_i^{}\rho_i^{}$ is the identity on $\pi Y^i$.
\item
$1-\rho_i^{}\theta_i^{}$ is an endomorphism of $\Sigma^{-i}\pi
X^{i+1}$.
We shall show it to be the composite
$$
\CD
\Sigma^{-i}\pi X^{i+1}
@>\Sigma \partial>> \Sigma^{-i}\pi X^{i+2}
@>\Sigma D>> \Sigma^{-i}\pi X^{i+1}
\endCD
$$
with $\partial$ and $D$ as above, satisfying
$1=D\partial+\partial D$.
\end{roenumerate} 
Note that for $i<-1$, $X^{i+1}=Y^{i}=0$, and there is nothing 
to do. We may assume that (i)-(iii) hold for
some $i$. We only need to show the induction step; that
is, if it holds for $i$ then it holds also for $i+1$.

It is easy to compute, in the derived category $D(R)$, the
composite $\alpha\beta$, with $\alpha$ and $\beta$ the morphisms
in the triangles below
$$
\CD
\Sigma^{-i-1} X^{i+1} @>\beta>>  Y^{i+1} @>>>  Y^i 
@>>> \Sigma^{-i} X^{i+1}\\
\Sigma^{-i-2} X^{i+2} @>>>  Y^{i+2} @>>>  Y^{i+1} 
@>\alpha>> \Sigma^{-i-1} X^{i+2}
\endCD
$$
The morphism $\alpha\beta$ is just $\Sigma^{-i-1}$
applied to the differential $X^{i+1}\ri X^{i+2}$. Applying the functor
$\pi$ we conclude the following. By the part (iii) of the induction
hypothesis, the composite
$$
\CD
\Sigma^{-i-1}\pi X^{i+1}
@>\partial>> \Sigma^{-i-1}\pi X^{i+2}
@>D>> \Sigma^{-i-1}\pi X^{i+1}
\endCD
$$
is equal to $1-\Sigma^{-1}(\rho_i^{}\theta_i^{})$. 
By the above, it factors
further as
$$
\CD
\Sigma^{-i-1}\pi X^{i+1}@>\pi\beta>>  \pi Y^{i+1}@>\pi\alpha>>
\Sigma^{-i-1}\pi X^{i+2}
@>D>> \Sigma^{-i-1}\pi X^{i+1}
\endCD
$$
Now look at the longer composite
$$
\CD
\Sigma^{-i-1}\pi X^{i+1}@>\pi\beta>>  \pi Y^{i+1}
@>D\circ(\pi\alpha)>> \Sigma^{-i-1}\pi X^{i+1} 
@>\pi\beta>>  \pi Y^{i+1}
\endCD
$$
It is equal to $(\pi\beta)\big[1-\Sigma^{-1}(\rho_i^{}\theta_i^{})\big]$.
The distinguished triangle
$$
\CD
\Sigma^{-i-1}\pi X^{i+1} @>\pi\beta>>  \pi Y^{i+1} @>>>  \pi Y^i 
@>\rho_i^{}>> \Sigma^{-i} \pi X^{i+1}
\endCD
$$
coupled with the fact that $\rho_i^{}$ is a split monomorphism,
guarantees that the triangle is really a split exact sequence
in $\cc$
$$
\CD
0 @>>> \Sigma^{-1}\pi Y^i 
@>\Sigma^{-1}\rho_i^{}>> \Sigma^{-i-1}\pi X^{i+1}
@>\pi\beta>> \pi \Sigma Y^{i+1} @>>> 0~.
\endCD
$$
But then $\pi\beta$ is a split epimorphism, and its
composite with $\Sigma^{-1}\rho_i^{}$ vanishes.
>From the vanishing of $\{\pi\beta\}\{\Sigma^{-1}\rho \}$
it follows that 
$$
(\pi\beta)\big[1-\Sigma^{-1}(\rho_i^{}\theta_i^{})\big]
=\pi\beta,
$$
and hence that
$$
\big[1-(\pi\beta)\circ D\circ(\pi\alpha)\big](\pi\beta)~=~0~.
$$
Since $\pi\beta$ is a split epimorphism, we conclude that
$$(\pi\beta)\circ D\circ(\pi\alpha)~=~1~.$$
But $\pi\alpha:\pi Y^{i+1} \ri
\Sigma^{-i-1}\pi X^{i+2}$ is nothing other than the map
$\rho_{i+1}^{}$, and if we put $\theta_{i+1}^{}=(\pi\beta)\circ D$,
then we have proved parts (i) and (ii) for $i+1$.

It only remains to establish (iii). But by construction,
$\Sigma\rho_{i+1}^{}\Sigma\theta_{i+1}^{}$ is given by the composite
$$
\CD
\Sigma^{-i}\pi X^{i+2}
@>D>> \Sigma^{-i}\pi X^{i+1} @>\pi\alpha>>
\pi \Sigma Y^{i+1} 
@>\Sigma\rho_{i+1}^{}>> \Sigma^{-i}\pi X^{i+2},
\endCD
$$
which is nothing other than a suspension of $\partial D$.
Hence this equals $1-D\partial$.

This completes the induction. Now choose $i>n$. The complex $Y^i$
is nothing other than $X\in\cb$, and by (i) we conclude that $\pi X$
is a direct summand of $\pi X^{i+1}=0$. 
It follows that $\pi X=0$.
This forces $X\in\ca$, but we know that $X\in\cb^c$. Hence
$X\in\ca\cap\cb^c=\ca^c$.
\eprf

\pro{P5.4}
Suppose 
$x$ is an object in $\cc^c$, and suppose $Tx=0$,
where 
$$T\pi~:~\cb~=~D(R)\ri \cd~=~D(\rs)~;~x \mapsto \sigma^{-1}x$$ 
is the functor induced by tensor with $\rs$, as in Proposition~\ref{P5.1}. 
Then $x=0$.
\epro

\prf
By Proposition~\ref{P5.3} we know that if $x\in\cb^c$, and if
$$
T\pi x~=~\br{^L\otimes_R^{}}x~=~0~,
$$
then $x\in\ca^c$, in other words $\pi x=0$. The Proposition
is therefore true for all objects $\pi x\in\cc^c$, with 
$x\in\cb^c$.

By \ref{T3.8.4}, the map  $\cb^c/\ca^c\ri\cc^c$ is fully
faithful, and $\cc^c$ is the smallest thick subcategory
containing
$\cb^c/\ca^c\subset\cc^c$.
By Corollary~4.5.12 of~\cite{Neeman99} we conclude the following.
For any object $t\in\cc^c$ there exists
an object $x\in\cb^c$ with 
$t\oplus \Sigma t\simeq\pi x$. 
Then $0=Tt\oplus\Sigma Tt\simeq T\pi x$,
and by the above this means $\pi x=0$. But $t$ is a direct
summand of $\pi x$; hence $t=0$.
\eprf

\section{Chain complex lifting}
\label{lift2}

We consider the  chain complex lifting problem of deciding if
a bounded complex $D$ of f.g. projective $\sigma^{-1}R$-modules
is chain equivalent to $\sigma^{-1}C$ for a bounded complex $C$
of f.g. projective $R$-modules. In terms of the functor 
$$T\pi~:~\cb~=~D(R)\ri \cd~=~D(\rs)~;~X \mapsto \br\oti X$$ 
the problem is to decide if a compact object $D$ in $\cd$
lifts to a compact object $C$ in $\cb$.
\medskip

To have any chance, we must start with a complex of induced
f.g. projective $\sigma^{-1}R$-modules, of the form
$$\CD
D~:~\cdots @>>> \sigma^{-1} x^{i-1} @>>> \sigma^{-1} x^{i} @>>>
\sigma^{-1} x^{i+1} @>>> \cdots
\endCD$$
with the $x^i$'s f.g. projective $R$-modules. We can write $D$ as
$$\CD
D~:~\cdots @>>> T\pi x^{i-1} @>>> T\pi x^{i} @>>>T\pi x^{i+1} @>>> \cdots~.
\endCD$$
By Proposition~\ref{P5.2}, this may be lifted uniquely to a chain
complex of objects in $\cc^c$ 
$$
\CD
\widetilde{D}~:~\cdots @>>> \pi x^{i-1} @>>>
 \pi x^{i} @>>>
\pi x^{i+1} @>>> \cdots
\endCD
$$
with $\cc= \cb/\ca= D(R)/D(R,\sigma)$.  Because all the objects lie in
the image of $\pi:\cb^c\ri\cc^c$ and because the functor
$\cb^c/\ca^c\ri\cc^c$ is fully faithful (see~\ref{T3.8.4}), we may view
the chain complex $\widetilde{D}$, uniquely, as lying in $\cb^c/\ca^c$. 
The next results discuss lifting this to $\cb^c$.

\lem{L103.6}
Given any diagram in $\cb^c/\ca^c$ of the form
$$
\CD
\pi x^0_{}@>>>\pi x^1_{} \\
@V|V\wr V           @. \\
\pi y^0_{} @.
\endCD
$$
where the vertical map is an isomorphism, we may complete to a
commutative square
$$
\CD
\pi x^0_{}@>>>\pi x^1_{} \\
@V|V\wr V           @V|V\wr V \\
\pi y^0_{} @>\pi f>> \pi y^1_{} 
\endCD
$$
where both vertical maps are isomorphisms in $\cb^c/\ca^c$, and $\pi
f:\pi y^0_{}\ri \pi y^1_{}$ is obtained by applying the functor $\pi$
to some map $f: y^0_{}\ri y^1_{}$.  
\elem

\prf
In $\cb^c/\ca^c$, we have a map $\pi y^0_{}\ri \pi x^1_{} $.
Such maps are equivalence classes of diagrams in $\cb^c$
$$
\CD
@. x^1_{} \\
@.           @VV\alpha V \\
y^0_{} @> f>>  y^1_{} 
\endCD
$$
with $\pi\alpha$ an isomorphism. Taking $\pi$ of this,
we get our commutative square.
\eprf

After these preliminaries, we  return to the problem of lifting
a bounded chain complex of induced f.g. projective $\rs$ modules 
$$
\CD
D~:~\cdots @>>> \s_{}^{-1} x_{}^{i-1} @>>>
 \s_{}^{-1} x_{}^{i} @>>>
\s_{}^{-1} x_{}^{i+1} @>>> \cdots
\endCD
$$
to a bounded chain complex of f.g. projective $R$-modules $C$
with a chain equivalence $\sigma^{-1}C \simeq D$.
We shall only treat the special case of a complex of length
3 in detail, but the general case is similar to this one. 
\medskip

We begin with
$$
\CD
D~: @>>>0 @>>>\s_{}^{-1} x_{}^{0} @>>>
 \s_{}^{-1} x_{}^{1} @>>>
\s_{}^{-1} x_{}^{2} @>>>\s_{}^{-1} x_{}^{3} @>>>0 @>>> 
\endCD
$$
In the derived category $\cd^c=D^c(\rs)$, we have several distinguished
triangles. We wish to consider three of them. They are
given by the mapping cones
$$
\CD
@>>>0 @>>> 0 @>>>\s_{}^{-1} x_{}^{0} @>>>
0 @>>>  \\
@.  @VVV @VVV @VVV @VVV @. \\
@>>>0 @>>> 0 @>>>\s_{}^{-1} x_{}^{1} @>>>
0 @>>>\\
@.  @VVV @VVV @VVV @VVV @. \\
@>>>0 @>>> \s_{}^{-1} x_{}^{0} @>>>\s_{}^{-1} x_{}^{1} @>>>
0 @>>>
\endCD
$$
and 
$$
\CD
@>>>0 @>>> 0 @>>>\s_{}^{-1} x_{}^{0} @>>>
\s_{}^{-1} x_{}^{1} @>>>0 @>>>  \\
@.  @VVV @VVV @VVV @VVV @VVV @. \\
@>>>0 @>>> 0 @>>>0 @>>>\s_{}^{-1} x_{}^{2} @>>>
0 @>>>\\
@.  @VVV @VVV @VVV @VVV @VVV @. \\
@>>>0 @>>> \s_{}^{-1} x_{}^{0} @>>>\s_{}^{-1} x_{}^{1} @>>>\s_{}^{-1} x_{}^{2} @>>>
0 @>>>
\endCD
$$
and
$$
\CD
@>>>0 @>>> 0 @>>>\s_{}^{-1} x_{}^{0} @>>>
\s_{}^{-1} x_{}^{1} @>>>\s_{}^{-1} x_{}^{2} @>>>0 @>>>  \\
@.  @VVV @VVV @VVV @VVV @VVV @VVV @. \\
@>>>0 @>>> 0 @>>>0 @>>>0 @>>>\s_{}^{-1} x_{}^{3} @>>>
0 @>>>\\
@.  @VVV @VVV @VVV @VVV @VVV @VVV @. \\
@>>>0 @>>> \s_{}^{-1} x_{}^{0} @>>>\s_{}^{-1} x_{}^{1} @>>>\s_{}^{-1} x_{}^{2} @>>>
\s_{}^{-1} x_{}^{3} @>>>
0 @>>>
\endCD
$$
The more abstract way of stating this is as follows.
In the derived category $D^c(\rs)$, the map $\s_{}^{-1} x_{}^{0}\ri
\s_{}^{-1} x_{}^{1}$ may be completed to a triangle
$$
\CD
\s_{}^{-1} x_{}^{0}@>>>\s_{}^{-1} x_{}^{1}@>>> X_1 @>>>\Sigma\s_{}^{-1} x_{}^{0}~.
\endCD
$$
This is the first of our three distinguished triangles above.
Because the composite $\s_{}^{-1} x_{}^{0}\ri
\s_{}^{-1} x_{}^{1}\ri
\s_{}^{-1} x_{}^{2}$ vanishes, we may factor the map $\s_{}^{-1} x_{}^{1}\ri
\s_{}^{-1} x_{}^{2}$ as $\s_{}^{-1} x_{}^{1}\ri X_1\ri
\s_{}^{-1} x_{}^{2}$. The factorization $X_1\ri \s_{}^{-1} x_{}^{2}$ is
unique, since its ambiguity is up to a map
$\Sigma\s_{}^{-1} x_{}^{0}\ri \s_{}^{-1}x_{}^{2}$, which must vanish
because it is a map from an object in ${D(\rs)}_{}^{\leq-1}$ to
an object in ${D(\rs)}_{}^{\geq0}$. The composite
$X_1\ri \s_{}^{-1}x_{}^{2}\ri \s_{}^{-1}x_{}^{3}$ must be zero, 
because it is the unique
factorization of the zero map
$\s_{}^{-1}x_{}^{1}\ri \s_{}^{-1}x_{}^{2}\ri \s_{}^{-1}x_{}^{3}$
through $\s_{}^{-1}x_{}^{1}\ri X_1$.
Next complete $X_1\ri \s_{}^{-1}x_{}^{2}$ to a triangle
 $$
\CD
X_1@>>>\s_{}^{-1} x_{}^{2}@>>> X_2 @>>>\Sigma X_1~.
\endCD
$$
This is the second of our three triangles above.
The vanishing of the composite $X_1\ri \s_{}^{-1}x_{}^{2}\ri \s_{}^{-1}x_{}^{3}$
tells us that the map $\s_{}^{-1}x_{}^{2}\ri \s_{}^{-1}x_{}^{3}$ must factor
as $\s_{}^{-1}x_{}^{2}\ri X_2\ri \s_{}^{-1}x_{}^{3}$. The factorization
is unique up to a morphism $\Sigma_{}^{-1}X_1\ri \s_{}^{-1}x_{}^{3}$,
and all such maps vanish because $\Sigma_{}^{-1}X_1\in  {D(\rs)}_{}^{\leq-1}$
while $\s_{}^{-1}x_{}^{3}\in {D(\rs)}_{}^{\geq0}$. Finally, we
may complete $X_2\ri\s_{}^{-1}x_{}^{3}$ to a triangle
$$\CD
X_2@>>>\s_{}^{-1} x_{}^{3}@>>> X_3 @>>>\Sigma X_2~.
\endCD$$
This gives the third triangle above.  The question is whether this
construction can be lifted to $D^c(R)$.  Since in $D^c(\rs)$ the choices of
factorizations were all unique, any lifting of the triangles and
factorizations to $D^c(R)$ will map, under the functor $T\pi$, to the
above.  The question is only whether the diagram of distinguished
triangles just constructed exists in $D^c(R)$.  We treat first the
problem of lifting by the functor $\pi$.  The obstructions are the
well-known Toda brackets (= Massey products) -- see Chapter IV.3 of
Gelfand and Manin \cite{GelfandManin}.  We give a detailed treatment of
this only in the simplest case, of a 3-dimensional complex.

\thm{TToda}
Let $D$ be a 3-dimensional chain complex of induced f.g. projective $\rs$-modules
$$
\CD
D~:~\s_{}^{-1} x_{}^{0} @>>> \s_{}^{-1} x_{}^{1} @>>>
\s_{}^{-1} x_{}^{2} @>>>\s_{}^{-1} x_{}^{3}~, 
\endCD$$
which we rewrite as
$$\CD
D~:~T\pi x_{}^{0} @>>> T\pi x_{}^{1} @>>> T\pi x_{}^{2} @>>>T\pi x_{}^{3}~.
\endCD$$
By Proposition~\ref{P5.2}, $D$ may be lifted uniquely to a chain complex
in $\cc^c$
$$\CD
\widetilde{D}~:~\pi x_{}^{0} @>f>> \pi x_{}^{1} @>g>>
\pi x_{}^{2} @>h>>\pi x_{}^{3}~.
\endCD$$
There is defined an element
$$\theta(D) \in \frac{\cc^c(\Sigma \pi x_{}^{0},\pi x_{}^{3})}
{\text{\rm Im}\{\cc^c(\Sigma \pi x_{}^{0},\pi x_{}^{2})
\oplus \cc^c(\Sigma \pi x_{}^{1},\pi x_{}^{3})\}}$$ 
such that the following conditions are equivalent :
\begin{itemize}
\item[(i)] $\theta(D)=0$.
\item[(ii)] There exist three triangles in $\cc^c$
$$\CD
\pi x_{}^{0} @>f>>
 \pi x_{}^{1} @>\alpha>>
X_1 @>>>\Sigma \pi x_{}^{0}\\
X_1 @>\beta>>
 \pi x_{}^{2} @>\gamma>>
X_2 @>>>\Sigma X_1\\
X_2 @>\delta>>
 \pi x_{}^{3} @>>>
X_3 @>>>\Sigma X_2
\endCD
$$
such that $g=\beta\alpha$ and $h=\delta\gamma$.
\item[(iii)] There exist three triangles in $\cb^c=D^c(R)$
$$
\CD
y_{}^{0} @>\tilde f>>
 y_{}^{1} @>\alpha>>
Y_1 @>>>\Sigma y_{}^{0}\\
Y_1 @>\beta>>
 y_{}^{2} @>\gamma>>
Y_2 @>>>\Sigma Y_1\\
Y_2 @>\delta>>
 y_{}^{3} @>>>
Y_3 @>>>\Sigma Y_2
\endCD
$$
and an isomorphism of chain complexes
in $\cc^c$
$$
\CD
\pi x_{}^{0} @>f>>
 \pi x_{}^{1} @>g>>
\pi x_{}^{2} @>h>>\pi x_{}^{3}\\
@V|V\wr V @V|V\wr V @V|V\wr V @V|V\wr V\\
\pi y_{}^{0} @>\pi\tilde f>>
\pi  y_{}^{1} @>\pi(\beta\alpha)>>
\pi y_{}^{2} @>\pi(\delta\gamma)>>\pi y_{}^{3}
\endCD
$$
\end{itemize}
In particular, $Y_3$ is a bounded f.g. projective $R$-module chain
complex such that $\sigma^{-1}Y_3 \simeq D$, solving the lifting problem.
\ethm
\prf In the first instance, we define $\theta(D)$.
We may always complete $f$ to a triangle
$$
\CD
\pi x_{}^{0} @>f>> \pi x_{}^{1} @>\alpha>> X_1 @>>>\Sigma \pi x_{}^{0}~.
\endCD
$$
The fact that $gf=0$ permits us to factor $g$ as
$$
\CD 
\pi x_{}^{1} @>\alpha>> X_1 @>\beta>> \pi x_{}^{2}~.
\endCD
$$
But $\beta$ is not well-defined; we may change our
choice by any element $\phi\in\cc^c(\Sigma \pi x_{}^{0},\pi x_{}^{2})$.

Now we may study the maps
$$
\CD
\pi x_{}^{1} @>\alpha>> X_1 @>\beta>> \pi x_{}^{2} @>h>> \pi x_{}^{3}~.
\endCD
$$
The composite $h\beta\alpha=hg=0$. We cannot be certain
that $h\beta$ vanishes, but we know that $h\beta$ composes
with $\alpha$ to give zero. From the triangle above,
it follows that $h\beta$ factors through 
$$
\CD
X_1 @>>> \T \pi x_{}^{0}
@>\theta>>\pi x_{}^{3}~.
\endCD
$$
The composite $h\beta$ will vanish if and only if
the map $\theta$ factors further as
$$
\CD
\T \pi x_{}^{0} @>\T f>>\T \pi x_{}^{1}
@>\psi>> \pi x_{}^{3}~.
\endCD
$$
We deduce that there is an obstruction to continuing the
process, given by $\theta\in\cc^c(\Sigma \pi x_{}^{0},\pi x_{}^{3})$.
And this $\theta$ is well defined up to adding
a $\phi\in\cc^c(\Sigma \pi x_{}^{0},\pi x_{}^{2})$
and a $\psi\in\cc^c(\Sigma \pi x_{}^{1},\pi x_{}^{3})$.
The element defined by
$$\theta(D)~=~[\theta] \in \frac{\cc^c(\Sigma \pi x_{}^{0},\pi x_{}^{3})}
{\text{\rm Im}\{\cc^c(\Sigma \pi x_{}^{0},\pi x_{}^{2})
\oplus \cc^c(\Sigma \pi x_{}^{1},\pi x_{}^{3})\}}
$$
is such that $\theta(D)=0$ if and only if we may choose
$\beta$ so that $h\beta=0$.
\smallskip

We now prove (i) $\Longleftrightarrow$ (ii) $\Longleftrightarrow$ (iii).
\smallskip

\noindent
(iii) $\Longrightarrow$ (ii) $\Longrightarrow$ (i) Obvious.\\
(i)  $\Longrightarrow$ (ii) 
As above, factor $g$ as $\beta\alpha$ 
so that $h\beta=0$.
Complete $\beta$ to a triangle
$$
\CD
X_1 @>\beta>>
 \pi x_{}^{2} @>\gamma>>
X_2 @>>>\Sigma X_1~.
\endCD
$$
Because $h\beta=0$, we may factor $h$ as
$$
\CD
 \pi x_{}^{2} @>\gamma>>
X_2 @>\delta>>\pi x_{}^{3}~.
\endCD
$$
Complete $\delta$ to a triangle
$$
\CD
X_2 @>\delta>>
 \pi x_{}^{3} @>>>
X_3 @>>>\Sigma X_2,
\endCD
$$
and we are done.\\
(ii) $\Longrightarrow$ (iii) We may assume that
we are given three triangles in $\cb^c/\ca^c$
$$
\CD
\pi x_{}^{0} @>f>>
 \pi x_{}^{1} @>\alpha'>>
X_1 @>>>\Sigma \pi x_{}^{0}\\
X_1 @>\beta'>>
 \pi x_{}^{2} @>\gamma'>>
X_2 @>>>\Sigma X_1\\
X_2 @>\delta'>>
 \pi x_{}^{3} @>>>
X_3 @>>>\Sigma X_2
\endCD
$$
such that 
$g=\beta'\alpha'$ and $h=\delta'\gamma'$.
Put  $y_{}^{0}=x_{}^{0}$.
Applying Lemma~\ref{L103.6} to the diagram
$$
\CD
\pi x_{}^{0} @>f>>
 \pi x_{}^{1} \\
@V1VV  @. \\
\pi y_{}^{0} @.
\endCD
$$
we may complete to a commutative square
$$
\CD
\pi x_{}^{0} @>f>>
 \pi x_{}^{1} \\
@V1VV  @V|V\wr V \\
\pi y_{}^{0} @>\pi\tilde f>> \pi y_{}^{1}~. 
\endCD
$$
Form the triangle
$$
\CD
y_{}^{0} @>\tilde f>> y_{}^{1} @>\alpha>> Y_1
@>>> \T y_{}^{0}~.
\endCD
$$
Applying the functor $\pi$, we have a commutative diagram
$$
\CD
\pi x_{}^{0} @>f>>
 \pi x_{}^{1} @>\alpha'>>
X_1 @>>>\Sigma \pi x_{}^{0}\\
@V1VV @V|V\wr V @. @V|V\wr V\\
\pi y_{}^{0} @>\pi \tilde f>> \pi y_{}^{1} @>\pi\alpha>> \pi Y_1
@>>> \T \pi y_{}^{0}
\endCD
$$
which may be extended to an isomorphism of triangles.
We have a commutative diagram
$$
\CD
 \pi x_{}^{1} @>\alpha'>>
X_1 @>\beta'>>\Sigma \pi x_{}^{2}\\
@V|V\wr V  @V|V\wr V @.\\
\pi y_{}^{1} @>\pi\alpha>> \pi Y_1
@.
\endCD
$$
which, again by Lemma~\ref{L103.6}, we may extend to 
$$
\CD
 \pi x_{}^{1} @>\alpha'>>
X_1 @>\beta'>> \pi x_{}^{2}\\
@V|V\wr V  @V|V\wr V @V|V\wr V\\
\pi y_{}^{1} @>\pi\alpha>> \pi Y_1
@>\pi\beta>> \pi y_{}^{2}~.
\endCD
$$
Now complete $\beta$ to a triangle
$$
\CD
 Y_1
@>\beta>> y_{}^{2} @>\gamma>> Y_2 @>>> \T Y_1~.
\endCD
$$
Again, we have a commutative diagram
$$
\CD
X_1 @>\beta'>>
 \pi x_{}^{2} @>\gamma'>>
X_2 @>>>\Sigma X_1\\
@V1VV @V|V\wr V @. @V|V\wr V\\
\pi Y_1 @>\pi \beta>> \pi y_{}^{2} @>\pi\gamma>> \pi Y_2
@>>> \T \pi Y_1
\endCD
$$
which we extend to an isomorphism of triangles.
Lemma~\ref{L103.6} allows us to extend the commutative diagram
$$
\CD
 \pi x_{}^{2} @>\gamma'>>
X_2 @>\delta'>>\pi x_{}^{3} \\
@V1VV @V|V\wr V @. \\
\pi y_{}^{2} @>\pi\gamma>> \pi Y_2
@.
\endCD
$$
to the diagram
$$
\CD
 \pi x_{}^{2} @>\gamma'>>
X_2 @>\delta'>>\pi x_{}^{3} \\
@V1VV @V|V\wr V @V|V\wr V \\
\pi y_{}^{2} @>\pi\gamma>> \pi Y_2 @>\pi\delta>>
\pi y_{}^{3}
\endCD
$$
Finally, we form the triangle
$$
\CD
 Y_2
@>\delta>> y_{}^{3} @>>> Y_3 @>>> \T Y_2
\endCD
$$
\eprf

\lem{L103.12}
Let $M$ and $N$ be f.g. projective $R$-modules.
There is a natural isomorphism
$$
\cc^c(\T \pi M,\pi N)~\cong~\text{\rm Tor}^R_2(\s^{-1}M^*,\s^{-1}N).
$$
In this formula, $M^*=\Hom_R(M,R)$ is the dual of $M$.
\elem

\prf
By Corollary~\ref{C8.4}, $\text{\rm Tor}^R_1(\s^{-1}R,\s^{-1}R)=0$.
The case $n=1$ of
Corollary~\ref{C8.6} then tell us that
\begin{eqnarray*}
\text{\rm Tor}^R_2(\s^{-1}R,\s^{-1}R)& = &H^{-1}(G\pi R) \\
&=& \cb(\T R,G\pi R)\\
&=&\cc^c(\T\pi R,\pi R).
\end{eqnarray*}
All we are doing is extending this isomorphism first
to free modules, then to their direct summands.
\eprf

\rmk{R103.13}
Lemma~\ref{L103.12} permits us to write the obstruction
class $\theta(D)$ of Theorem~\ref{TToda} as lying in the group
$$
\frac
{\text{\rm Tor}^R_2(\s^{-1}{\{x^0_{}\}}^*,\s^{-1}x^3_{})}
{\text{\rm Im}\{
\text{\rm Tor}^R_2(\s^{-1}{\{x^0_{}\}}^*,\s^{-1}x^2_{})
\oplus\text{\rm Tor}^R_2(\s^{-1}{\{x^1_{}\}}^*,\s^{-1}x^0_{})
\}}~.$$
Note that if ${\rm Tor}^R_2(\rs,\rs)=0$ this group is 0.
\ermk

\rmk{RTodageneral}
It is easy to generalize this to longer complexes. Given a bounded
f.g. projective $R$-module chain complex
$$D~:~\cdots\ri\s^{-1}x^{i-1}_{}\ri\s^{-1}x^{i}_{}\ri
\s^{-1}x^{i+1}_{}\ri\cdots$$
there is a series of obstructions to lifting all the associated
triangles to $D^c(R)$
$$\theta_{i,j}(D) \in \frac{\cc^c(\Sigma^{j-2} \pi x_{}^{i},\pi x_{}^{i+j})}
{\text{\rm Im}\{\cc^c(\Sigma^{j-2} \pi x_{}^{i},\pi x_{}^{i+j-1})
\oplus \cc^c(\Sigma^{j-2} \pi x_{}^{i+1},\pi x_{}^{i+j})\}}$$
for $j \geq 3$.
As in Remark \ref{R103.13} these are related to Tor-groups
by a spectral sequence,
and are 0 if ${\rm Tor}^R_*(\rs,\rs)=0$ for $* \geq 2$.
\ermk

\section{Waldhausen's approximation and localization theorems}
\label{SWald}

We have been studying noncommutative localization using derived
categories techniques. Next we want to apply our results
to deduce $K$-theoretic consequences. In order to do so,
we briefly review some results of Waldhausen's. 
\medskip

Let ${\mathbf C}$ be a category with cofibrations and weak equivalences.
Out of ${\mathbf C}$ Waldhausen constructs a spectrum, denoted 
$K({\mathbf C})$. In Thomason's \cite{ThomTro}, the category ${\mathbf C}$
is assumed to be a full subcategory of the category of chain
complexes over some abelian category, the cofibrations
are maps of complexes which are split monomorphisms in each 
degree, and the weak equivalences are the quasi-isomorphisms.
We shall call such categories {\em permissible Waldhausen categories.}

\rmk{BiWald}
Thomason's term for them is {\em complicial biWaldhausen categories.}
\ermk

\nin
Given a permissible Waldhausen category ${\mathbf C}$, one can form its derived
category; just invert the weak equivalences. We denote this
derived category by $D({\mathbf C})$. We have two major theorems
here, both of which are special cases of more general theorems
of Waldhausen.

\thm{Approx}{\bf (Waldhausen's Approximation Theorem).}\ \ 
Let $F:{\mathbf C}\ri{\mathbf D}$ be an exact functor of essentially small
permissible
Waldhausen categories (categories of chain complexes, as
above). Suppose that the induced map of derived categories
\[
D(F)~:~D({\mathbf C})\ri D({\mathbf D})
\]
is an equivalence of categories. Then the induced map of
spectra
\[
K(F)~:~K({\mathbf C})\ri K({\mathbf D})
\]
is a homotopy equivalence.
\ethm

\nin
In this sense, Waldhausen's $K$-theory is almost an invariant
of the derived categories. To construct it, one needs to
have a great deal more structure. One must begin with
a permissible category with cofibrations and weak equivalences.
But the Approximation Theorem asserts that the dependence
on the added structure is not strong. 

\thm{Local}{\bf (Waldhausen's Localization Theorem).}\ \ 
Let $\mathbf A$, $\mathbf B$ and $\mathbf C$ be essentially small
permissible Waldhausen categories.  Suppose
\[
{\mathbf A}\ri {\mathbf B}\ri{\mathbf C}
\] 
are exact functors of permissible
Waldhausen categories. Suppose further that the induced
triangulated functors of derived categories
\[
D({\mathbf A})\ri D({\mathbf B})\ri D({\mathbf C})
\]
compose to zero, and that the natural map
\[
\CD
D({\mathbf B})/D({\mathbf A}) @>>> D({\mathbf C})
\endCD
\] 
is an equivalence of categories. Then the sequence of spectra
 \[
K({\mathbf A})\ri K({\mathbf B})\ri K({\mathbf C})
\]
is a homotopy fibration.\hfill$\Box$
\ethm

To obtain a homotopy fibration using Waldhausen's localization theorem,
we need to produce three permissible Waldhausen categories, and a
sequence
\[
{\mathbf A}\ri {\mathbf B}\ri{\mathbf C}
\]
so that 
\[
\CD
D({\mathbf B})/D({\mathbf A}) @>>> D({\mathbf C})
\endCD
\] 
is an equivalence of categories.  In particular, we want to find
triangulated categories $\ca^c=D({\mathbf A})$, $\cb^c=D({\mathbf B})$
and $\cc^c=D({\mathbf C})$ so that $\cc^c=\cb^c/\ca^c$.  Of course, it
is not enough to just find the triangulated categories $\ca^c$, $\cb^c$
and $\cc^c$; to apply the localization theorem, we must also find the
permissible Waldhausen categories $\mathbf A$, $\mathbf B$ and
$\mathbf C$,
and the exact functors
\[
{\mathbf A}\ri {\mathbf B}\ri{\mathbf C}~.
\]
In Theorem~\ref{Approx}, we learned that the $K$-theory is largely
independent of the choices of $\mathbf A$, $\mathbf B$ and
$\mathbf C$.
In this article, we shall allow ourselves some latitude.  Thomason is
careful to check, in \cite{ThomTro}, that the choices of permissible
Waldhausen categories can be made; we shall consider this a technical
point, and explain only how to produce $\ca^c=D({\mathbf A})$,
$\cb^c=D({\mathbf B})$ and $\cc^c=D({\mathbf C})$.  We shall also commit
the notational sin of writing $K(\ca^c)$ for $K({\mathbf A})$, where
$\ca^c=D({\mathbf A})$, and similarly $K(\cb^c)$ for $K({\mathbf B})$,
and $K(\cc^c)$ for $K({\mathbf C})$.
\medskip

As the notation of the previous paragraph was designed to suggest, we
want to apply the results $\ca^c\subset\cb^c$ and $\cc^c$, with $\ca$,
$\cb$ and $\cc$ as we have seen them in the previous sections.  That
is, $\cb=D(R)$ is the derived category of a ring $R$, $\ca$ is
generated by a set $\s$ of morphisms in $\cb^c$, and $\cc=\cb/\ca$.  By
the discussion above, we have a fibration in \kth\
\[
\CD
K(\ca^c)@>>>  K(\cb^c)@>>> K(\cb^c/\ca^c)~.
\endCD
\]
In Theorem~\ref{T3.8} we learned that 
the natural map $\cb^c/\ca^c\ri\cc^c$ is fully faithful,
and that up to splitting idempotents it is an equivalence.
Grayson's cofinality theorem then tells us that the map
\[
\CD
K(\cb^c/\ca^c) @>>> K(\cc^c)
\endCD
\]
induces an isomorphism $K_i(\cb^c/\ca^c)\ri K_i(\cc^c)$
when $i>0$, while $K_0(\cb^c/\ca^c)\ri K_0(\cc^c)$
is injective. 
We conclude that, up to the failure of surjectivity in $\pi_0$,
\[
\CD
K(\ca^c)@>>>  K(\cb^c)@>>> K(\cc^c)
\endCD
\]
is a homotopy fibration. 
\medskip

We know also that $\cb^c=D^c(R)$, and in Proposition~\ref{P5.2} we
produced a functor $T:\cc^c\ri \cd^c=D^c(\rs)$.  For any ring $S$, we
have $K(S)=K(D^c(S))$; Waldhausen's \kth\ of the derived category agrees
with Quillen's \kth\ of $S$.  Applying this to the rings $R$ and $\rs$,
we have

\thm{Tlovelier}\label{lovelier}
In the diagram
\[
\CD
K(\ca^c)@>>>  K(\cb^c)@>>> K(\cc^c) \\
@. @V|V\wr V    @VK(T)VV \\
@. K(R)  @. K(\rs)
\endCD
\]
not only is the top row a fibration up to the failure of
surjectivity on $\pi_0^{}$, but $K(\cb^c)$ agrees with 
Quillen's $K(R)$, and $T:\cc^c \ri \cd^c=D^c(\sigma^{-1}R)$
induces a morphism $K(T):K(\cc^c)\ri K(\rs)$.\hfill$\Box$
\ethm

In particular, Theorem \ref{lovelier} gives a long exact sequence
$$\dots \ri K_n(R) \ri K_n(\cc^c) \ri K_n(R,\sigma) \ri K_{n-1}(R) \ri \dots~,$$
with $K_*(R,\sigma)=K_{*-1}(\ca^c)$.  In the coming sections, we shall
study the range in which the map $K(T):K_*(\cc^c) \ri
K_*(\sigma^{-1}R)$ induces an isomorphism in homotopy.

\section[$K_0$-isomorphism]
{$T$ induces a $K_0$-isomorphism}
\label{S8}

Let the notation be as in Theorem~\ref{Tlovelier}.
In this section, we shall prove that the functor
$$\CD
T~:~\cc^c @>>> \cd^c~=~D^c(\rs)
\endCD$$
induces an isomorphism in $K_0$.  We shall do it through a sequence of
lemmas.  We remind the reader that $\cb=D(R)$ has a standard {\it
t}-structure, and that the functor $G\pi$ behaves well with respect to
it.  See Remark~\ref{t-structure} and Lemma~\ref{L6.1}.

\lem{L6.2}
Let $n$ be an integer.
Let $X\in\cb^c$ be an object of $\cb^{\leq n}$, and let $P$ be a 
f.g. projective $R$-module. Then the functor $T:\cc^c\ri \cd^c$
of Proposition~\ref{P5.1} gives a homomorphism
\[
\CD
\cc^c(\pi\Sigma^{-n} P,\pi X) @>>> \cd^c(T\pi\Sigma^{-n}P,T\pi X)~.
\endCD
\]
We assert that this map is an isomorphism.
\elem

\prf
By translation, we may assume $n=0$.
We need to prove the map injective and surjective. Let us prove
surjectivity first. Recall that $X\in\cb^{\leq 0}$ is isomorphic
to a chain complex of f.g. projectives
\[
\ri X^m\ri X^{m+1}\ri \cdots\ri
X^{-1}\ri X^0\ri 0\ri 0\ri 
\]
This makes $T\pi X$ the chain complex
\[
\ri \s^{-1}_{}X^m\ri \s^{-1}_{}X^{m+1}\ri \cdots\ri
\s^{-1}_{}X^{-1}\ri \s^{-1}_{}X^0\ri 0\ri 0\ri
\]
Let $P$ be a f.g. projective $R$-module,
concentrated in degree 0.
Now the complex of $\rs$-modules $T\pi P$
is a single projective module $\sigma^{-1} P$, concentrated in degree $0$.
Any map in the derived category, 
from the bounded above complex of projectives 
$\s^{-1}_{}P=\br\otimes_R^{} P$ to the complex
$\br{^L\otimes_R^{}}X$,
can be represented by a chain map. There is a 
map $\s^{-1}_{} P\ri \s^{-1}_{}X_0$ inducing it.
By Proposition~\ref{P5.2}, this comes from a map 
$\pi P\ri \pi X^0$. But then the 
composite
\[
\CD
\pi P @>>> \pi X^0
@>>> \pi X
\endCD
\]
gives a map $\pi P\ri 
\pi X$ in $\cc^c$, inducing $T\pi P\ri
T\pi X$.

This proved the surjectivity. For the injectivity, note
that there is a short exact sequence of chain complexes
\[
\begin{array}{ccccccccccccccc}
\ri &0&\ri& 0&\ri &\cdots&\ri&
0&\ri &X^0&\ri& 0&\ri& 0&\ri \\
&\downarrow& &\downarrow & & & & 
\downarrow& &\downarrow & &\downarrow& &\downarrow &\\
\ri &X^m&\ri& X^{m+1}&\ri &\cdots&\ri&
X^{-1}&\ri &X^0&\ri& 0&\ri& 0&\ri \\
&\downarrow& &\downarrow & & & & 
\downarrow& &\downarrow & &\downarrow& &\downarrow &\\
\ri &X^m&\ri& X^{m+1}&\ri &\cdots&\ri&
X^{-1}&\ri &0 &\ri& 0&\ri& 0&\ri 
\end{array} 
\]
Write the corresponding triangle as
\[
\CD
X^0 @>>> X @>>> Y @>>> \Sigma X^0.
\endCD
\]
We have a triangle 
\[
\CD
\pi X^0 @>>>\pi X @>>>\pi Y @>>>\pi \Sigma X^0.
\endCD
\]
Let $P$ be a f.g. projective $R$-module,
concentrated in degree 0.
Suppose we are given a map $\pi P\ri \pi X$.
Composing to $Y$, we deduce a
map 
\[
\CD
\pi P @>>>\pi X @>>>\pi Y~.
\endCD
\]
By adjunction, this corresponds to a map
\[
\CD
P @>>>G\pi Y~,
\endCD
\]
which must vanish. After all, $Y\in\cb^{\leq -1}$,
and by Lemma~\ref{L6.1} it follows that $G\pi Y$ is
also in $\cb^{\leq -1}$. The map from a projective
object $P$ in degree $0$ to the complex $G\pi Y\in\cb^{\leq -1}$
must vanish.
\medskip

It follows that the map $\pi P \ri\pi X$ must factor
as
\[
\CD
\pi P @>>>\pi X^0  @>>>\pi X .
\endCD
\]
Now assume that the map vanishes in $\cd^c=D^c(\rs)$. 
That is, the composite
\[
\CD
\sigma^{-1} P @>>>\sigma^{-1} X^0 
@>>>\br{^L\otimes_R^{}} X 
\endCD
\]
vanishes in $\cd^c$. Then it must be null homotopic. The map 
$\sigma^{-1}P \ri \sigma^{-1}X^0$ must factor as
\[
\CD
\sigma^{-1}P @>>>\sigma^{-1}X^{-1} @>>>\sigma^{-1}X^0.
\endCD
\]
By Proposition~\ref{P5.2}, this tells us that the map
$\pi P\ri \pi X^0$ must factor as
\[
\CD
\pi P @>>>\pi X^{-1}  @>>>\pi X^0  
\endCD
\]
and hence the map
\[
\CD
\pi P @>>>\pi X^{-1}  @>>>\pi X^0 @>>> \pi X 
\endCD
\]
must vanish. 
\eprf

\lem{L6.3}
Let $n$ be an integer.
Let $Z$ be an object of $\cc^c$, and suppose for all
$r\geq n$, $H^r(T Z)=0$. Then there is an
object $X\in\cb^c$, that is a bounded complex of projective
$R$-modules
\[
\ri 0\ri X^m\ri X^{m+1}\ri \cdots\ri
X^{\ell-1}\ri X^\ell\ri 0\ri  
\]
so that $Z$ is a direct summand of $\pi X$, and $\ell\leq n$.
\elem

\prf
By suspending, we may assume $n=0$. By \ref{T3.8.4}, we may 
certainly find an $X$ with $Z$ a direct summand of $\pi X$.
What is not clear is that we may choose $X$ to be a complex
\[
\ri 0\ri X^m\ri X^{m+1}\ri \cdots\ri
X^{\ell-1}\ri X^\ell\ri 0\ri  
\]
with $\ell\leq0$. Assume therefore that $\ell>0$, and we shall
show that we may reduce $\ell$ by $1$.

We recall the short exact sequence of chain complexes
\[
\begin{array}{ccccccccccccccc}
\ri&0 &\ri &0&\ri& 0&\ri &\cdots&\ri&
0&\ri &X^\ell&\ri& 0&\ri  \\
&\downarrow& &\downarrow& &\downarrow & & & & 
\downarrow& &\downarrow & &\downarrow& \\
\ri&0 &\ri &X^m&\ri& X^{m+1}&\ri &\cdots&\ri&
X^{\ell-1}&\ri &X^\ell&\ri& 0&\ri\\
&\downarrow& &\downarrow& &\downarrow & & & & 
\downarrow& &\downarrow & &\downarrow& \\
\ri&0&\ri &X^m&\ri& X^{m+1}&\ri &\cdots&\ri&
X^{\ell-1}&\ri &0 &\ri& 0&\ri
\end{array} 
\]
It gives a triangle which we write as
\[
\CD
\Sigma^{-\ell}X^\ell @>>> X @>a>> Y @>>>\Sigma^{-\ell+1}X^\ell.
\endCD
\]
We also have that $Z$ is a direct summand of $\pi X$. That is,
there are maps
\[
\CD
\pi X @>b>> Z @>c>> \pi X
\endCD
\]
so that $bc=1_Z^{}$. Now we wish to consider the
composite
\[
\CD
\pi\Sigma^{-\ell}X^\ell @>>> \pi X @>b>> Z @>c>> \pi X~.
\endCD
\]
We know that 
$X^\ell$ is a f.g. projective $R$-module,
and $X\in\cb$ lies  in $\cb^{\leq\ell}$. The conditions are
as in Lemma~\ref{L6.2}. In order to prove that the composite
vanishes, it suffices to prove that $T$ of it vanishes, in
$\cd^c=D^c(\rs)$.\\
\indent But in $\cd^c$ the map becomes the composite
\[
\CD
T\pi\Sigma^{-\ell}X^\ell @>>> T\pi X  @>>> T Z @>>> T\pi X~.
\endCD
\]
We assert that already the shorter composite,
$T\pi\Sigma^{-\ell}X^\ell \ri T\pi X 
\ri T Z$ must vanish. After all, it is a map
\[
\CD
T\pi\T^{-\ell} X^\ell  
@>>> T Z  
\endCD
\]
By hypothesis, $TZ$ vanishes above degree 0.
It is quasi-isomorphic to a complex of $\rs$-modules in degree
$\leq0$. And $T\pi\T^{-\ell} X^\ell=  
\Sigma^{-\ell}\sigma^{-1}X^\ell$ is a
single projective $\rs$-module, concentrated in degree $\ell>0$.
Hence the vanishing.
The composite
\[
\CD
\pi\Sigma^{-\ell}X^\ell @>>> \pi X @>b>> Z @>c>> \pi X
\endCD
\]
must therefore vanish.
Since $c$ is a split monomorphism, we deduce that the composite
\[
\CD
\pi\Sigma^{-\ell}X^\ell @>>> \pi X @>b>> Z 
\endCD
\]
also vanishes.\\
\indent But now the triangle
\[
\CD
\pi\Sigma^{-\ell}X^\ell @>>> \pi X @>a>> \pi Y @>>>\pi\Sigma^{-\ell+1}X^\ell
\endCD
\]
tells us that the map $b:\pi X\ri Z$ must factor as
\[
\CD
\pi X @>a>> \pi Y @>\beta>> Z.
\endCD
\]
The composite 
\[
\CD
Z @>c>>\pi X @>a>> \pi Y @>\beta>> Z
\endCD
\]
is the identity, and hence $Z$ is a direct summand of $\pi Y$, with
$Y$ the complex
\[
\begin{array}{ccccccccccccccc}
\ri&0&\ri &X^m&\ri& X^{m+1}&\ri &\cdots&\ri&
X^{\ell-1}&\ri &0 &\ri& 0&\ri
\end{array} 
\]
\eprf

\lem{L6.4}
Let $n$ be an integer.
Let $Z$ be an object of $\cc^c$, and suppose for all
$r\geq n$, $H^r(T Z)=0$. Given any f.g. projective $R$-module $P$,
and any map
\[
\CD
T\pi P~=~\sigma^{-1}P @>a>> H^n(T Z)~,
\endCD
\]
there is a map in $\cc^c$
\[
\CD
\pi \Sigma^{-n}P @>\mu>> Z
\endCD
\]
so that $H^n(T\mu)=a$.
\elem

\prf
By translating, we may assume $n=0$.
Let $Z$ be an object of $\cc^c$, and suppose for all
$r\geq 0$, $H^r(T Z)=0$. By Lemma~\ref{L6.3},
there exists a complex $X\in D^c(R)$
\[
\ri 0\ri X^m\ri X^{m+1}\ri \cdots\ri X^{-1}\ri X^0\ri 0\ri  
\]
so that $Z$ is a direct summand of $\pi X$.
We have two maps
\[
\CD
\pi X @>b>> Z @>c>> \pi X
\endCD
\]
so that $bc=1_Z^{}$. This gives us two maps
\[
\CD
T\pi X @>Tb>>T Z @>Tc>> T\pi X
\endCD
\]
with $(Tb)(Tc)=1$. Given any map 
\[
\CD
\sigma^{-1}P @>a>> H^0(T Z),
\endCD
\]
we can form the composite
\[
\CD
\sigma^{-1}P @>a>> H^0(TZ)
@>H^0(Tc)>> H^0(T\pi X)~.
\endCD
\]
Of course, $T\pi X$ is just the chain complex
\[
 \cdots\ri
\sigma^{-1}X^{-1}\ri \sigma^{-1}X^0\ri 0\ri  
\]
and any map from a projective $\sigma^{-1}P$ to
$H^0(T\pi X)$ lifts to a map
\[
\CD
\sigma^{-1}P @>>> T\pi X~.
\endCD
\]
By Lemma~\ref{L6.2}, the above map is $T\gamma$, for a (unique)
map  
\[
\CD
\pi P @>\gamma>> \pi X~.
\endCD
\]
Now let $\mu$ be the composite
\[
\CD
\pi P @>\gamma>> \pi X @>b>> Z~.
\endCD
\]
Applying the functor $H^0\circ T$, we compute  
$H^0(T\mu)$ to be the composite
\[
\CD
\sigma^{-1}P @>a>> H^0(TZ)
@>H^0(Tc)>> H^0(T\pi X) @>H^0(Tb)>> 
H^0(TZ),
\endCD
\]
which is nothing other than the map $a$.
\eprf

\lem{L6.5}
For any f.g. projective $\rs$-module $M$, there is a 
canonically unique object $\widetilde M\in\cc^c$ so that
\sthm{L6.5.1}
$H^0(T\widetilde M)=0$ for $n\neq0$.
\esthm
\sthm{L6.5.2}
$H^0(T\widetilde M)=M$.
\esthm
\nin
The functor $H^0(T-)$ is an equivalence
of categories between the full subcategory of objects
$\widetilde M\in\cc^c$ and f.g. projective 
$\rs$-modules.
\elem

\prf
Let us first prove existence. Let $M$ be a f.g.
projective $\rs$-module. There exists a $\rs$-module
$N$, so that $M\oplus N\cong \br^r$. There is an idempotent
$\br^r\ri \br^r$ which is the map
\[
\CD
M\oplus N @>1^{}_M\oplus 0^{}_N>> M\oplus N~.
\endCD
\]
Write this map as
$1^{}_M\oplus 0^{}_N:T\pi R^r\ri T\pi R^r$.
By Proposition~\ref{P5.2}, there is a unique lifting
$e:\pi R^r\ri \pi R^r$. The uniqueness of the lifting
allows us to easily show that $e^2=e$. But idempotents
split in $\cc$, by Proposition~1.6.8 of \cite{Neeman99}.
Define $\widetilde M$ by splitting the idempotent $e$.

Then $H^n(T\widetilde M)$ is computed
by splitting the idempotent $H^n(Te)$ on $H^n(\rs^r)$;
this gives us zero when $n\neq0$, and $M$ when $n=0$.
We have proved the existence of a $\widetilde M$
satisfying \ref{L6.5.1} and \ref{L6.5.2}.

Now suppose $X$ is an object of $\cc^c$, and that
\begin{roenumerate}
\item $H^n(T X)=0$ for $n\neq0$,
\item $H^0(TX)=M$.
\end{roenumerate}
We wish to produce an isomorphism $\widetilde M\ri X$.
In any case, we have a map
\[
\CD
\sigma^{-1}R^r @>>> M~=~H^0(TX)~,
\endCD
\]
namely the projection to the direct summand.
By Lemma~\ref{L6.4}, there is a map
\[
\CD
\pi R^r @>>> X~,
\endCD
\]
which induces the projection. We may form the composite
\[
\CD
\widetilde M @>>> \pi R^r @>>> X,
\endCD
\]
and it is very easy to check that the map
 \[
\CD
T\widetilde M  @>>> T X
\endCD
\]
is a homology isomorphism, hence an isomorphism in
$D^c(\rs)$.
If we complete $\widetilde M  \ri X$ to a triangle
in $\cc^c$
\[
\CD
\widetilde M  @>>> X @>>> Y @>>> \Sigma\widetilde M, 
\endCD
\]
then $TY=0$.  But by Proposition~\ref{P5.4} it then follows that $Y=0$,
and $\widetilde M \ri X$ is an isomorphism.\\

Finally it remains to check that $\cc(\widetilde M,\widetilde N)
=\Hom_{\rs}^{}(M,N)$. By the construction of $\widetilde M$
and $\widetilde N$ as direct summands of $\pi R^r$ and $\pi R^s$,
this reduces to knowing that 
$$\cc(\pi R^r,\pi R^s)~=~\cd(T\pi R^r,T\pi R^s)~.$$ 
But we know this from Proposition~\ref{P5.2}.
\eprf

\thm{T6.6}
The map $T:\cc^c\ri \cd^c$ of Proposition~\ref{P5.1} induces
a $K_0$-isomorphism.
\ethm

\prf
We have maps of categories
\[
\CD{\script P}(\rs) @>a>> \cc^c @>T>> \cd^c
\endCD
\]
with ${\script P}(\rs)$ the category of f.g. projective $\sigma^{-1}R$-modules.
The map $T$ is given by Proposition~\ref{P5.1}; the map $a$
takes a f.g. projective $\rs$-module $M$
to $a(M)=\widetilde M$. In $K$-theory, the composite
\[
K_0(\rs)\ri K_0(\cc^c)\ri K_0(\cd^c)
\]
is clearly an isomorphism. To prove that both maps are
isomorphisms, it suffices to show that the
map $K_0(a):K_0(\rs)\ri K_0(\cc^c)$ is onto. This is what we shall
do.

Let $Z$ be an object of $\cc^c$.  We want to show that its class
$[Z]\in K_0(\cc^c)$ lies in the image of $K_0(a)$.  We shall prove this
by induction on the length of $TZ$.  In any case, $TZ$ is an object of
$\cd^c=D^c(\rs)$; it is a bounded complex of f.g.  projective
$\rs$-modules.

Suppose the length of $TZ$ is zero.  Replacing $Z$ by a suspension,
this means that $H^n(TZ)=0$ unless $n=0$.  But then $H^0(TZ)=M$ must be
a f.g.  projective $\rs$-module, and by Lemma~\ref{L6.5} we know that
$Z$ is (canonically) isomorphic to $\widetilde M$.  Thus $Z$ is in the
image of $a$.

Suppose now that we know the induction hypothesis.  We are given
$n\geq0$.  We know that if $Z$ is an object of $\cc^c$ so that the
length of $TZ$ is $\leq n$, then the class $[Z]\in K_0(\cc^c)$ lies in
the image of $K_0(a)$.  Let $Z$ be a complex of length $n+1\geq1$. 
Replacing $Z$ by a suspension, this means that $H^r(TZ)=0$ unless
$-n-1\leq r\leq 0$.  Now $H^0(TZ)$ is a finitely presented
$\rs$-module; we may choose a f.g.  free $R$-module $F$, and a
surjection $ \sigma^{-1}F\ri H^0(TZ)$.  By Lemma~\ref{L6.4}, there is a
map
\[
\CD
\pi F @>>> Z
\endCD
\]
lifting this surjection. Form the triangle in $\cc^c$
\[
\CD
\pi F @>>> Z @>>> Y @>>> \Sigma \pi F.
\endCD
\]
It is easily computed that the length of $TY$ is $\leq n$, so by
induction $[Y]$ lies in the image of $K_0(a):K_0(\rs)\ri K_0(\cc^c)$. 
Clearly $[\pi F]= [\widetilde{\sigma^{-1}F}]$ also lies in the image of
$K_0(a)$, and the triangle tells us that $[Z]=[Y]+[\pi F]$.  
\eprf

\section[$K_1$-isomorphism]{$T$ induces a $K_1$-isomorphism}

\label{S9}

In Proposition~\ref{P5.1} we produced a triangulated functor of triangulated 
categories
$$T~:~\cc^c~\ri \cd^c~=~D^c(\rs)~.$$
In Section~\ref{S8} we proved that the induced map $K_0(T):K_0(\cc^c)\ri
K_0(\cd^c)$ is an isomorphism. The main result
of this section is that so is $K_1(T):K_1(\cc^c)\ri K_1(\cd^c)$.
First we must address a point concerning Waldhausen \kth. 
\medskip

Let $\bB$ be the category whose objects are all bounded chain complexes
of f.g. projective $R$-modules. The morphisms
in $\bB$ are the chain maps. The cofibrations are the maps
which are split monomorphisms in each degree. The weak 
equivalences are the homology isomorphisms. Clearly, $\bB$
is a model for the triangulated category $\cb^c=D^c(R)$.
\medskip

Let $\bA$ be the full subcategory of all objects in $\bB$
whose image in $\cb^c$ is contained in $\ca$. Then $\bA$ is
a model for the triangulated category $\ca^c$, and the
inclusion $\bA\ri\bB$ is a model for the map $\ca^c\ri\cb^c$.
Let $\bC$ be the same category as $\bB$, with the same cofibrations,
but different weak equivalences. The weak equivalences in $\bC$
are the maps in $\bC=\bB$ whose mapping cone lies in $\bA\subset\bB$.
We have exact functors of Waldhausen categories
$$\bA\ri\bB\ri\bC$$
which induce the maps of triangulated categories
$$\ca^c\ri\cb^c\ri\cb^c/\ca^c~.$$
Let $\bD$ be the Waldhausen category of all bounded chain
complexes of f.g. projective $\sigma^{-1}R$-modules.
The cofibrations are the maps which are split monomorphisms
in each degree. The weak equivalences are the homology
isomorphisms. The functor $X\mapsto \sigma^{-1}X$ is
clearly an exact functor of Waldhausen categories 
$\bB\ri\bD$. In Proposition~\ref{P5.2} we showed that it
factors as
$$\bB\ri\bC\ri\bD~.$$
Waldhausen's localization theorem tells us that there is
a homotopy fibration in Waldhausen $K$-theory
$$K(\bA)\ri K(\bB)\ri K(\bC)~;$$
somewhat loosely, we have been referring to this homotopy
fibration as
$$K(\ca^c)\ri K(\cb^c)\ri K(\cb^c/\ca^c)~.$$
The map we have been calling $K(\cb^c/\ca^c)\ri K(\cc^c)$
is an isomorphism except possibly on $K_0$. Since this
section deals only with $K_1$, there is no point in 
explaining what this map is, on the level of Waldhausen
models.
There is also a map $K(\bC)\ri K(\bD)$. This is 
the map we have loosely been referring to
as $K( \cb^c/\ca^c)\ri K(\cd^c)$.
\medskip

In the proof of Theorem~\ref{T6.6}, we also introduced
a functor ${\script P}(\rs) \ri \cc^c$. We do not know a
Waldhausen model for this map. But in the following
discussion we shall show that there is an induced map
on $K_1$. The group $K_1(\rs)$ is generated by determinants
of automorphisms of free (or projective) modules. Certainly,
the following generate: given any projective $R$-module $P$,
and an automorphism $\phi:\sigma^{-1}P\ri \sigma^{-1}P$,
the determinant of $\phi$ is an element of $K_1(\rs)$, and the
collection of all determinants of all $\phi$'s generates $K_1(\rs)$.
We want to produce a map $K_1(\rs)\ri K_1(\bC)$; to define the
map, it suffices to say what it does on all $\phi$'s
as above.
\medskip

To define what the map does to $\phi$, let us remind ourselves that the
zero-space of the spectrum $K(\bC)$ has a Gillet-Grayson model (see
\cite{Gillet-Grayson87}) , which we denote $GG(\bC)$.  That is, there
is a homotopy equivalence
$$GG(\bC)~\simeq~ \Omega^\infty K(\bC)~.$$
The space $GG(\bC)$ is an $H$-space, and hence
$$K_1(\bC)~=~\pi_1^{}K(\bC)~=~\pi_1^{}GG(\bC)~=~H_1 GG(\bC)~.$$
Starting with an automorphism $\phi:\sigma^{-1}P\ri \sigma^{-1}P$,
we need to produce a class in the first homology group $H_1 GG(\bC)$.
\medskip

We note that, by Proposition~\ref{P5.2}, $\phi:\sigma^{-1}P\ri
\sigma^{-1}P$ corresponds to a unique automorphism
$$\varphi~:~\pi P\ri \pi P~.$$
This is an automorphism defined in $\cb^c/\ca^c$, and $\bC$ is a 
Waldhausen model for $\cb^c/\ca^c$. It follows that there
exist weak equivalences $a:Q\ri P$ and $b:Q\ri P$, with
$\varphi=ab^{-1}$. But then $P$ and $Q$ are 0-cells in
the Gillet-Grayson model $GG(\bC)$. The weak equivalences
$a:Q\ri P$ and $b:Q\ri P$ are 1-cells. Now $[a]-[b]$ is a
cycle, that is an element in
$H_1\big(GG(\bC)\big)=\pi_1^{}\big(GG(\bC)\big)=K_1(\bC)$. 
We leave it to the reader to check that the map sending $\phi$ to 
$[a]-[b]$ extends to a well-defined homomorphism $K_1(\rs)\ri K_1(\bC)$.
\medskip

The composite $K_1(\rs)\ri K_1(\bC)\ri K_1(\bD)$ is easily 
seen to be an isomorphism; in our looser notation, it is the map
$$K_1(\rs)\ri K_1(\cb^c/\ca^c)~=~K_1(\cc^c)\ri K_1(\cd^c)~.$$
To prove that both maps are isomorphisms it suffices
therefore to check that $K_1(\rs)\ri K_1(\bC)$ is epi.
We have a localization exact sequence
$$\CD
K_1(\bB) @>>> K_1(\bC) @>>> K_0(\bA) @>>> K_0(\bB)~.
\endCD$$
Note that $K(\bB)=K\big(D^c(R)\big)=K(R)$.
The composite $K_1(R)\ri K_1(\rs)\ri K_1(\bC)$ is easily
computed to agree with the natural $K_1(R)= K_1(\bB)\ri K_1(\bC)$;
we deduce a commutative diagram where the bottom row is exact
$$\CD
K_1(R) @>>> K_1(\rs) @. @. \\
@V1VV  @V\psi VV @. @. \\
K_1(R) @>>> K_1(\bC) @>>> K_0(\bA) @>>> K_0(R)~.
\endCD
$$
To prove $\psi$ epi, it suffices to show that the 
composite 
$$
\CD
K_1(\rs) @.  \\
@V\psi VV @.  \\
K_1(\bC) @>>> K_0(\bA) 
\endCD
$$
surjects to the kernel of $K_0(\bA) \ri K_0(R)$. But the composite
is easy to compute. Take an automorphism $\phi:\rs\otimes P\ri 
\rs\otimes P$ as above, which corresponds as above to an automorphism
$$\varphi~:~\pi P\ri \pi P~.$$
Choose weak equivalences $a:Q\ri P$ and $b:Q\ri P$, with
$\varphi=ab^{-1}$. Then $\phi$ gets sent to $[A]-[B]$, where
$$\begin{array}{l}
A~:~\dots \ri 0\ri Q \stackrel a\ri P\ri0\ri \dots~,\\[1ex]
B~:~\dots \ri 0\ri Q \stackrel b\ri P\ri0\ri~.
\end{array}$$
It will therefore suffice to show that every element in the
kernel of the map $K_0(\bA) \ri K_0(R)$ can be expressed
as a difference $[A]-[B]$, as above. We shall prove something
stronger.

\thm{T7.1}
Every element in $K_0(\bA)=K_0(\ca^c)$ is a linear
combination of complexes of length $\leq1$. That is,
it may be written as $\sum\pm[A_i]$, with $A_i\in\ca^c$
being complexes of f.g. projective $R$-modules of the form
$$\dots \ri 0\ri X  \ri Y\ri0\ri \dots~.\eqno{\Box}$$
\ethm

Note that with $K_0$, it makes no difference whether we
compute $K_0$ of a triangulated category, or $K_0$ of
its model; from now on we can forget all about models.
Before the proof of Theorem~\ref{T7.1}, let us state the main corollary

\cor{C7.2}
Every object in the kernel of the map $K_0(\ca^c)\ri
K_0(\cb^c)=K_0(R)$ is of the form $[A]-[B]$, where $A$ is a complex
$$\dots \ri 0\ri Q \stackrel a\ri P\ri0\ri \dots
$$
and $B$ is a complex
$$
\dots \ri 0\ri Q \stackrel b\ri P\ri0\ri \dots
$$
By the discussion preceding Theorem~\ref{T7.1}, this means
that the map $T:\cc^c\ri \cd^c$ induces a $K_1$-isomorphism.
\ecor

\nin
{\em Proof that Corollary~\ref{C7.2} follows from Theorem~\ref{T7.1}.}
\ \ 
Suppose we have an element of the kernel of the map
$K_0(\ca^c)\ri K_0(\cb^c)=K_0(R)$.
By Theorem~\ref{T7.1}, just by virtue of being an element of
$K_0(\ca^c)$, it has an expression as $\sum\pm[A_i]$, with $A_i\in\ca^c$
being complexes of f.g. projective $R$-modules of the form
$$
\dots \ri 0\ri X_i  \ri Y_i\ri0\ri \dots
$$
Recalling that $[\Sigma A_i]=-[A_i]$, up to changing signs
in the sum we may assume that all the $X_i$ are in degree $-1$,
all the $Y_i$ in degree $0$. Collecting together
all the terms of equal sign, we may rewrite the sum as
$$
[\oplus A_i]-[\oplus B_j].
$$
That is, we have an element $[A]-[B]$ in the kernel
of $K_0(\ca^c)\ri K_0(\cb^c)=K_0(R)$, where $A,B$ are complexes of the form
$$\begin{array}{l}
A~:~ \dots \ri 0\ri A^{-1} \stackrel a\ri A^0\ri0\ri \dots~,\\[1ex]
B~:~ \dots \ri 0\ri B^{-1} \stackrel b\ri B^0\ri0\ri \dots~.
\end{array}$$
The fact that $[A]-[B]$ lies in the kernel of the map $K_0(\ca^c)\ri
K_0(\cb^c)=K_0(R)$ tells us that, in $K_0(R)$, there is an identity
$$[A^{-1}]+[B^0]~=~[B^{-1}]+[A^0]~.$$
This in turn says that there is a projective $R$-module $X$, 
and an isomorphism
$$A^{-1}\oplus B^0\oplus X~\cong~B^{-1}\oplus A^0\oplus X~.$$
The object $[A]\in\ca^c$ is isomorphic to the complex
$$
\CD
\dots @>>> 0 @>>> A^{-1}\oplus B^0\oplus X 
@>a\oplus 1_{B^0}^{}\oplus 1_X^{}>> A^{0}\oplus B^0\oplus X 
@>>> 0 @>>>\dots
\endCD
$$
while the object $[B]\in\ca^c$ is isomorphic to the complex
$$
\CD
\dots @>>> 0 @>>> B^{-1}\oplus A^0\oplus X 
@>b\oplus 1_{A^0}^{}\oplus 1_X^{}>> B^{0}\oplus A^0\oplus X 
@>>> 0 @>>>\dots
\endCD
$$
Put $Q=A^{-1}\oplus B^0\oplus X\,\cong\,B^{-1}\oplus A^0\oplus X$,
and $P=A^{0}\oplus B^0\oplus X$. Then $A$ is isomorphic in 
$\ca^c$ to a complex
$$\dots \ri 0\ri Q \stackrel \alpha\ri P\ri0\ri \dots $$
and $B$ is isomorphic in 
$\ca^c$ to a complex
$$\dots\ri 0\ri Q \stackrel \beta\ri P\ri0\ri \dots$$
as required. \hfill{$\Box$}

\bigskip

\nin
{\em Proof of Theorem~\ref{T7.1}.}\ \ 
It remains to prove Theorem~\ref{T7.1}. Let $X$ be an object
of $\ca^c$. We need to show that
the class $[X]$ in $K_0(\ca^c)$ can be written as a linear
combination of classes of objects of length $\leq1$.\\
\indent Because $X\in\ca^c\subset\cb^c$,
we have that $X$ is a bounded complex of f.g. projective
$R$-modules. Suspending suitably, we may assume it has the form
$$
\ri 0\ri X^{-m}\ri X^{-m+1}\ri \cdots\ri
X^{-1}\ri X^0\ri 0\ri  
$$
If $m\leq1$ we are done; the complex has length $\leq 1$.  The proof is
by induction.  Assume we are given an integer $n\geq1$.  Assume further
that, for every $X\in\ca^c$ of length $m\leq n$, $[X]$ is equal in
$K_0(\ca^c)$ to a linear combination of complexes of length $\leq1$. 
Take a complex $X$ as above, with $m=n+1\geq2$.  We need to show that
it can also be expressed as a linear combination of complexes of length
$\leq1$.

By Lemma~\ref{Lbound.7}, we have that every object of $\ca^c$ is
isomorphic to a direct summand of an object in $\cs$, with $\cs$
defined as in Definition~\ref{Dbound.4}.  Choose a chain complex $Y\in
\cs$ and maps $X\ri Y\ri X$ composing to the identity on $X$.  Clearly,
the map $H^0(Y)\ri H^0(X)$ must be surjective.

By Lemma~\ref{Lbound.6}, there exist exists a triangle
$$
\CD
U @>>> Y @>>> V @>>> \Sigma U
\endCD
$$
with $U\in\cs[1,\infty)$ and $V\in\cs(-\infty,1]$.
The composite
$$
\CD
U @>>> Y @>>> X 
\endCD
$$
is a map from $U\in\cs[1,\infty)$ to $X\in\cb^{\leq0}$,
which must vanish. It follow that $Y\ri X$ factors as
$Y\ri V\ri X$. And since $H^0(Y)\ri H^0(X)$ is epi and
factors through $H^0(V)$, we deduce that $H^0(V)\ri H^0(X)$
must be epi. Replacing $Y$ by $V$, we may assume
$Y\in\cs(-\infty,1]$.

Next we apply Lemma~\ref{Lbound.6} again, this time to deduce
that $Y\in\cs(-\infty,1]$
can be expressed as the mapping cone on a map $U\ri V$,
with $U\in\cs(-\infty,0]$ and $V\in\cs[-1,1]$. There is a triangle
$$
\CD
U @>>> V @>>> Y @>>> \Sigma U~,
\endCD
$$
hence an exact sequence
$$
\CD
H^0( V) @>>>H^0(  Y) @>>> H^1( U)~=~0~.
\endCD
$$
The map $H^0( V) \ri H^0(  Y) \ri H^0(X)$
is the composite of two epis, hence
is epi. Replacing $Y\ri X$ by the composite
$V\ri Y\ri X$, we may assume $Y\in\cs[-1,1]$.

The last time we apply Lemma~\ref{Lbound.6} is to express
$Y$ as the mapping cone of a map $U\ri Z$, with
$U\in\cs[0,1]$ and $Z\in\cs[0,1]$. The
only observations we wish to make is that $H^1(Z)\ri H^1(Y)$
is epi, and that in $K_0(\ca^c)$, $[Z]$ and $[Y]=[Z]-[U]$ are
both linear combinations of objects in $\ca$ of length $\leq1$. 
Let us summarize: we have constructed maps 
$Z\ri Y\ri X$, with 
\begin{roenumerate}
\item
$Z\in\cs[0,1]$, $Y\in\cs[-1,1]$,
\item 
$H^0(Y)\ri H^0(X)$ epi,
\item 
$H^1(Z)\ri H^1(Y)$ epi,
\item 
both $[Y]$ and $[Z]$ are linear combinations of objects in $\ca$ 
of length $\leq1$. 
\end{roenumerate}

Form the mapping cone on the map $Y\ri X$, to obtain a triangle
$$
\CD
Y @>>> X @>>> X' @>>> Y~.
\endCD
$$
Since $Y\in\cs[-1,1]$ while $X\in\ca^c$ is supported
on the interval $[-m,0]$ with $m\geq2$, the mapping cone $X'$
is an object of $\ca^c$ supported in $[-m,0]$.
The long exact sequence in homology gives
$$
\CD
H^0(Y) @>>>  H^0(X) @>>> H^0(X') @>>> H^1(Y
) @>>> H^1(X)~.
\endCD
$$
We have $H^1(X)=0$, while $H^0(Y) \ri  H^0(X)$ is an epimorphism.
Hence $H^0(X')=H^1(Y)$. But we know that
the map $H^1(Z)\ri H^1(Y)$ is an epimorphism,
by (iii). And $Z$ is a complex of the form
$$\dots \ri 0\ri Z^0\ri Z^1\ri 0\ri \dots $$
that is a complex of length $\leq1$. It follows that we can extend
the epimorphism
$\beta:H^1(Z)\ri H^0(X')$ to a map from the 
presentation $Z$; there is a map $\Sigma Z\ri X'$,
inducing $\beta$ in $H^0$.
We may  form the mapping cone, obtaining a triangle
$$
\CD
\Sigma Z @>>> X' @>>> X'' @>>> \Sigma^2 Z~.
\endCD
$$
Since $\Sigma Z\in\cs[-1,0]$ and $X'$ is supported on $[-m,0]$,
we conclude that $X''$ is supported on $[-m,0]$
But now the long exact homology sequence
$$
\CD
H^0(\Sigma Z) @>\alpha>>  H^0(X') @>>> H^0(X'') @>>> H^1(\Sigma Z) 
\endCD
$$
has $\alpha$ surjective, while $H^1(\Sigma Z)=0$.
We conclude that $H^0(X'')=0$. The complex $X''$ is
supported in the interval $[m,-1]$. By induction,
its class in $K_0(\ca^c)$ is a linear combination
of complexes of length $\leq1$.

But now the triangles above give the identities
$$
[X]=[X']+[Y],\qquad\qquad [X']=[X'']-[Z]
$$ 
and hence $[X]=[X'']+[Y]-[Z]$, and all the
terms on the right may be expressed as linear
combinations of complexes of length $\leq1$.
\hfill{$\Box$}

\section{$T$ is an equivalence if and only if the Tor-groups vanish}
\label{Sequivalence}

In Proposition~\ref{P5.1} we constructed a functor 
$$T~:~\cc^c~=~(D(R)/D(R,\sigma))^c \ri \cd^c~=~D^c(\rs)~.$$
In Sections~\ref{S8} and \ref{S9} we showed that $T$ induces
an isomorphism in $K_0$ and $K_1$. For higher $K$-theory,
the useful result we have is a necessary and sufficient
conditions for the functor $T$ to be an equivalence
of categories. An equivalence of categories trivially induces
an isomorphism in \kth. It is very easy to see a necessary
condition:

\lem{L104.1}
If $T:\cc^c \ri \cd^c$ is an equivalence
then, for all $n\neq0$, $H^n(G\pi R)=0$.
\elem

\prf
Suppose $T$ is an equivalence. We have isomorphisms
\begin{eqnarray*}
H^n(G\pi R) &=& \cb^c( R, \T^nG\pi R)\\
&=&\cc^c(\pi R,\T^n\pi R)\\
&=&\cd(T\pi R,\T^nT\pi R)~.
\end{eqnarray*}
But $T\pi R=\rs$,
and in the category $\cd=D(\rs)$, the maps $\rs\ri\T^n\rs$
all vanish, whenever $n\neq0$.
\eprf

\rmk{R104.2}
The main result of this section is that the converse 
of Lemma~\ref{L104.1} holds.
The condition is also sufficient.
\ermk

\pro{P104.3}
The groups $\{H^n(G\pi R),n\neq0\}$ all vanish if and
only if the groups $\{\text{\rm Tor}^R_n(\rs,\rs),n\neq0\}$
all vanish. 
\epro

\prf
By Lemma~\ref{L6.1}, the groups $H^n(G\pi R)$ vanish when $n>0$,
while the groups $\text{\rm Tor}^R_n(\rs,\rs)$ clearly vanish
when $n<0$. By Corollary~\ref{C8.4}, 
$\text{\rm Tor}^R_1(\rs,\rs)=0$.
By Corollary~\ref{C8.6}, the smallest 
 non-zero integer $n$ for which $H^{-n}(G\pi R)\neq0$
is equals the smallest integer $n\neq-1$ for
which  $\text{\rm Tor}^R_{n+1}(\rs,\rs)\neq0$.
\eprf

\cor{C3.1.9}
Suppose $H^n(G\pi R)=0$ for all $n\neq0$.
Put $S=\rs$. We assert
\sthm{C3.1.9.1}
$S\otimes_R^{}S=S$ via the multiplication map, and 
\esthm
\sthm{C3.1.9.2}
$\hbox{\rm Tor}^R_i(S,S)=0$ for all $i>0$.
\esthm
\ecor

\prf
\ref{C3.1.9.2} is immediate from Proposition~\ref{P104.3},
while \ref{C3.1.9.1} follows from Corollary~\ref{C8.4}.
\eprf

Next we shall formally study the consequences of \ref{C3.1.9.1} and
\ref{C3.1.9.2}.  Let us begin with an easy general observation, about
natural transformations of functors respecting coproducts.

\lem{L3.1.6.5}
Let $\cs$ and $\ct$ be triangulated categories, let $F$ and $G$
be triangulated functors $\cs\ri\ct$, and let $\rho:F\Longrightarrow
G$
be a natural transformation commuting with the suspension. Define
the full subcategory $\ci\subset\cs$ by the formula
$$\text{\rm Ob}(\ci)~=~
\{x\in {\rm Ob}(\cs)\mid \rho_x^{}\text{ is an isomorphism}\}~.$$
Then the category $\ci$ is triangulated. If furthermore both
$F$ and $G$ commute with coproducts, then $\ci$ contains
all coproducts in $\cs$ of its objects.
\elem

\prf
Left to the reader.
\eprf

\lem{L3.1.10}
Let $R \ri S$ be a ring homomorphism such that
$R$ and $S$ satisfy the hypotheses \ref{C3.1.9.1} and
\ref{C3.1.9.2}; we remind the reader
\bd
\item[\ref{C3.1.9.1}]
$S\otimes_R^{}S=S$ via the multiplication map, and 
\item[\ref{C3.1.9.2}] $\hbox{\rm Tor}^R_i(S,S)=0$ for all $i>0$.
\ed
Define the functor
$$\Theta~:~D(S) \ri D(S)~;~X\mapsto S{^L\otimes_R}X~.$$
Multiplication defines a natural
transformation $\mu:\Theta\Longrightarrow 1$. We assert that $\mu$
is an isomorphism.
\elem

\prf
The functors $\Theta$ and $1$ are triangulated functors respecting
coproducts, and the natural transformation $\mu$ commutes with
the suspension. By Lemma~\ref{L3.1.6.5}, the full subcategory
of all $x\in D(S)$ for which $\mu$ is an isomorphism
is triangulated category closed under coproducts. By
\ref{C3.1.9.1} and \ref{C3.1.9.2}, it also contains the object
$S\in D(S)$. From Lemma~\ref{L3.1.7}, it follows that $\mu_x^{}$
is an isomorphism for all $x\in D(S)$.
\eprf

\thm{T3.1.11}
Let $\alpha:R\ri S$ be a ring homomorphism such that $R$ and $S$ satisfy
the hypotheses \ref{C3.1.9.1} and \ref{C3.1.9.2}; we remind the reader
\bd
\item[\ref{C3.1.9.1}]
$S\otimes_R^{}S=S$ via the multiplication map, and 
\item[\ref{C3.1.9.2}]
$\hbox{\rm Tor}^R_i(S,S)=0$ for all $i>0$.
\ed
Put $\cb=D(R)$, $\cd=D(S)$. There is a functor
$\pi:\cb\ri\cd$ taking $X\in D(R)$ to $S{^L\otimes_R}X\in D(S)$.
The functor $\pi$ has a fully faithful right adjoint $G:\cd\ri\cb$.
The unit of adjunction $\eta:1\Longrightarrow G\pi$ is identified 
by saying that the two maps below are naturally isomorphic
$$\CD
X@>\eta_X^{}>> G\pi X\quad,\quad
R{^L\otimes_R}X @>\alpha\otimes X>> S{^L\otimes_R}X~.
\endCD
$$
The functor $G$ respects coproducts. Finally, there exists
a thick subcategory $\ca\subset\cb$ so that the pair of
adjoints $\pi,G$ are a Bousfield localization for the pair
$\ca\subset\cb$.
\ethm

\prf
The functor 
$$\pi~:~\cb~=~D(R)\ri \cd~=~D(S)~;~X\mapsto S{^L\otimes_R}X$$
is left adjoint to the forgetful
functor $G:\cd \ri \cb$ taking a complex of $S$-modules
and just viewing it as a complex of $R$-modules. This is just
the derived category version of the
classical fact that for any $R$-module $M$ and $S$-module $N$
$$\Hom_R^{}(M,N)~=~\Hom_S^{}(S\otimes_RM,N)~.$$
The functor $G$ clearly respects coproducts, and the unit of
adjunction clearly identifies as
$$\CD
R{^L\otimes_R}X @>\alpha\otimes X>> S{^L\otimes_R}X.
\endCD$$

The first statement that requires proof is that $G$ is
fully faithful. Suppose therefore that $X$ and $Y$ are
complexes of $S$-modules. Then
$$
\begin{array}{rclcl}
\Hom_{D(R)}^{}(X,Y) &=&\Hom_{D(S)}^{}(S{^L\otimes_R}X,Y)
&\qquad&\text{because $\pi$ is left adjoint to $G$}\\[1ex]
&=&\Hom_{D(S)}^{}(X,Y)
&\qquad&\text{by Lemma~\ref{L3.1.10}, $S{^L\otimes_R}X=X$.}
\end{array}
$$
Since the homomorphisms are the same in $D(S)=\cd$ as
in $D(R)=\cb$, the inclusion $G$ is fully faithful.

This means that the fully faithful inclusion $G:\cd\ri\cb$
has a left adjoint, and hence, by Proposition~9.1.18
in \cite{Neeman99}, a Bousfield colocalization
exists for the pair $\cd\subset\cb$. Put 
$\ca=\cd^\perp$, in the notation of Chapter~9
loc. cit. By Corollary~9.1.14
loc. cit. a Bousfield localization functor exists
for the pair $\ca\subset\cb$, with $\cd={^\perp\ca}$,
in other words $G:\cd\ri\cb$ is the right adjoint of
the quotient map. The quotient map must be the left
adjoint of $G$, that is $\pi:D(R)\ri D(S)$, as
above.
\eprf

It remains to apply Theorem~\ref{T3.1.11} to the 
special case of Corollary~\ref{C3.1.9}, that is where
$S=\rs$.

\thm{Tfinal}
Let $R$ be a ring, $\s$ a set of morphisms of f.g. projective $R$-modules.
Set
$$\ca~=~D(R,\sigma)~,~\cb~=~D(R)~,~\cc~=~D(R)/D(R,\sigma)~,~
\cd~=~D(\sigma^{-1}R)$$
with $\ca \subset \cb$ the subcategory generated by $\s$. 
Let $\pi:\cb\ri\cc=\cb/\ca$ be the projection,
$G:\cc\ri\cb$ its right adjoint. The following conditions
are equivalent~:
\begin{itemize}
\item[(i)] $H^n(G\pi R)=0$ for all $n\neq0$,
\item[(ii)] ${\rm Tor}^R_i(\rs,\rs)=0$ for $i \geq 1$,
\item[(iii)] the natural functor $T:\cc\ri \cd$,
of Proposition~\ref{P5.1}, is an equivalence of categories,
\item[(iv)] the natural functor $T:\cc^c\ri \cd^c=D^c(\rs)$
is an equivalence of categories.
\end{itemize}
\ethm
\prf
The equivalence of (i) and (ii) is given by Proposition \ref{P104.3}.\\
\indent The implication (iv) $\Longrightarrow$ (i) is Lemma \ref{L104.1}.\\
\indent The implication (iii) $\Longrightarrow$ (iv) is obvious.\\
\indent It remains to prove (ii) $\Longrightarrow$ (iii).
The composite $T\pi$ is the functor $X\mapsto \br\oti X$,
and by (ii) and Theorem~\ref{T3.1.11} we know that it is a
projection $\cb\ri\cb/{^\perp\cd}$. Suppose we could
show that ${^\perp\cd}=\ca$. Then both $T\pi$ and $\pi$
would be identified as the projection $\cb\ri\cb/\ca=\cb/{^\perp\cd}$,
and by the universal property of the Verdier quotient,
$T$ would be an equivalence.
We are reduced to showing ${^\perp\cd}=\ca$.

The category ${^\perp\cd}$ is very explicit: it is the kernel
of the functor $T\pi$. That is, ${^\perp\cd}$ is the collection
of all $X\in\cb$ with $\br\oti X=0$. 
By Lemma~\ref{L100.3}, $\ca\subset{^\perp\cd}$. We must
prove the reverse inclusion.
Suppose therefore that $\br\oti X=0$; we need to
show that $X\in\ca$.

We may form a triangle
$$
\CD
a @>>> X @>\eta_X^{}>> G\pi X @>>> \T a.
\endCD
$$
By \ref{T3.1.5.5} we have that $a\in\ca$, and hence
$\br\oti a=0$. We are assuming $\br\oti X=0$;
the triangle tells us that 
$\br\oti G\pi X=0$. We want to prove that $X\in\ca$,
in other words we want to prove that $G\pi X=0$.
It clearly suffices to prove that 
$G\pi X$ is isomorphic to $\br\oti G\pi X$. We shall 
prove, more precisely, that
 the natural map
$$
\CD
R\oti G\pi X @>f\oti 1_{G\pi}^{}>> \br\oti G\pi X
\endCD
$$
is an isomorphism, where $f:R\ri \rs$ is the natural map
$\eta_R^{}:R\ri H^0(G\pi R)=G\pi R$.

The map $f\oti 1_{G\pi}^{}$ is a natural transformation between
triangulated functors respecting coproducts, and $f\oti 1_{G\pi}^{}$
commutes with the suspension.  Form the full subcategory $\ci$ given by
$$\text{\rm Ob}(\ci)~=~\{x\in {\rm Ob}(\cd)\mid 
f\oti 1_{G\pi x}^{}\text{ is an isomorphism}\}~.$$
We wish to show that $\ci$ is all of $\cb=D(R)$.
By Lemma~\ref{L3.1.6.5}, we know that $\ci$ is a triangulated 
subcategory of $\cb$, closed under coproducts.
By Lemma~\ref{L3.1.7},
it suffices to show that $R\in\ci$.
But $f\oti 1_{G\pi R}^{}$ is the map 
$f\oti 1:\rs\ri\br\oti\br$. Now (ii) tells us that the multiplication
map $\mu:\br\oti\br\ri\rs$ is a homology isomorphism,
hence an isomorphism in the derived category. It is
clear that the composite
$$
\CD
\rs @>f\oti 1>> \br\oti \br @>\mu>> \rs
\endCD
$$
is the identity, forcing $f\oti 1$ to be the two-sided
inverse of the invertible map $\mu$.
\eprf

As in the Introduction~:

\dfn{stablyflat} A noncommutative localization $\sigma^{-1}R$ is {\it
stably flat over $R$} if it satisfies the equivalent conditions of Theorem
\ref{Tfinal}.  
\edfn

\thm{Plifting} Suppose we are given a chain complex $D\in D^c(\rs)$,
that is a bounded chain complex of f.g.  projective $\rs$ modules. 
Suppose every module in the chain complex is induced; that is, they are
all of the form $\s^{-1}P$, with $P$ a finitely generated projective
$R$-module.  If $\rs$ is stably flat over $R$ then there exists a
complex $C\in D^c(R)$ and a homotopy equivalence $D\cong\s^{-1}C$. 
\ethm

\prf
We are assuming that $\rs$ is stably flat over $R$. 
By Theorem~\ref{Tfinal} this means that the functor
$T:\cc^c\ri \cd^c$ is an equivalence. By~\ref{T3.8.4}, the map
$\cb^c/\ca^c\ri\cc^c=\cd^c$ is fully faithful, and $\cc^c=\cd^c$
is the smallest thick subcategory containing
$\cb^c/\ca^c$. Now note that the objects
of the Verdier quotient $\cb^c/\ca^c$ are the same
as the objects of $\cb^c$, and the projection
map $\pi:\cb^c\ri\cb^c/\ca^c$ induces a surjective map
in $K_0$. Combining this with Proposition~4.5.11 of~\cite{Neeman99},
an object $D\in\cd^c=D^c(\rs)$ is in the image of the
functor $T\pi:\cb^c\ri\cd^c$ if and only if the class of $[D]$ in
$K_0(\cd^c)$ lies in the image of the map 
$K_0(T\pi):K_0(\cb^c)\ri K_0(\cd^c)$.
The functor $T\pi$ is by definition the functor taking
a complex $C\in\cb^c=D^c(R)$ to $\s^{-1}C\in\cd^c=D^c(\rs)$.
We chose our complex $D$ to be a complex of induced modules.
The class of $[D]$ most certainly lies in
the image of $K_0(T\pi):K_0(\cb^c)\ri K_0(\cd^c)$. We conclude that
there exists a $C$ with $D$ isomorphic in $D^c(\rs)$ to $\s^{-1}C$.
\eprf

\thm{Kexact} If $\sigma^{-1}R$ is stably flat over $R$ then the functor 
$T:\cc^c \ri \cd^c=D^c(\sigma^{-1}R)$ induces isomorphisms
$$T~:~K_*(\cc^c)~=~K_*(D^c(R)/D^c(R,\sigma)) \ri K_*(\cd^c)~=~K_*(\sigma^{-1}R)$$
and there is a localization exact sequence in algebraic $K$-theory 
$$\dots \ri K_n(R) \ri K_n(\sigma^{-1}R) \ri K_n(R,\sigma) \ri 
K_{n-1}(R) \ri \dots $$
\ethm
\prf Combine Theorems \ref{Tlovelier} and \ref{Tfinal}.
\eprf

In Theorem~\ref{Plifting} we saw that, if $\rs$ is stably flat over
$R$ and $D\in D^c(\rs)$ is a complex of induced modules, then
there exists a complex $C\in D^c(R)$ and a homotopy equivalence
$\s^{-1}C\cong D$. Now we want a more precise version of this.

\pro{Prelifting}
If $C\in D^c(R)$ is a bounded complex of f.g. projective
$R$-modules, and if $\s^{-1}C$ is homotopy equivalent to
a complex $D\in D^c(\rs)$ vanishing outside an interval $[0,n]$
with $n\geq1$, then there exists a $B\in D^c(R)$, vanishing
outside $[0,n]$, with $\s^{-1}B$ homotopy equivalent to
$D\cong  \s^{-1}C$.
\epro

\prf
We need to show that $C$ can be shortened. Suppose therefore that
$C$ is the complex
\[
\CD
\cdots @>>> 0 @>>> C^{-1} @>>> C^0 @>>> \cdots @>>> C^n@>>> 0 @>>> \cdots
\endCD
\]
and assume that there is a homotopy equivalence
of $\s^{-1}C$ with a shorter complex, that is a 
commutative diagram
\[
\CD
 @>>> 0 @>>> \s^{-1}C^{-1} @>\partial>> \s^{-1}C^0 
@>>> \cdots @>>> \s^{-1}C^n @>>> 0 @>>> \\
@. @VVV @VVV @VVV @. @VVV @VVV @. \\
 @>>> 0 @>>> 0 @>>> D^0 @>>> \cdots @>>> D^n @>>> 0 @>>> \\
@. @VVV @VVV @VVV @. @VVV @VVV @. \\
 @>>> 0 @>>> \s^{-1}C^{-1} @>\partial>> \s^{-1}C^0 
@>>> \cdots @>>> \s^{-1}C^n @>>> 0 @>>> 
\endCD
\]
so that the composite is homotopic to the identity.
In particular, there is a map $d:\s^{-1}C^0\ri \s^{-1}C^{-1}$
so that $d\partial:\s^{-1}C^{-1}\ri \s^{-1}C^{-1}$ is the
identity.

By Proposition~\ref{P5.2}, the map $d:\s^{-1}C^0\ri \s^{-1}C^{-1}$
lifts uniquely to a map $d':\pi C^0\ri \pi C^{-1}$.
By Proposition~\ref{Pbound.8}, the map $d'$ can be represented
as $\alpha^{-1}\beta$, where $\alpha$ and $\beta$ are, respectively,
the chain maps
\[
\CD
 @>>> 0 @>>> 0 @>>> C^{-1} @>>> 0 @>>> \\
@. @VVV @VVV @VVV  @VVV @. \\
 @>>> 0 @>>> X @>r>> Y @>>> 0 @>>>
\endCD
\]
and
\[
\CD
 @>>> 0 @>>> 0 @>>> C^{0} @>>> 0 @>>> \\
@. @VVV @VVV @VgVV  @VVV @. \\
 @>>> 0 @>>> X @>r>> Y @>>> 0 @>>>
\endCD
\]
The fact that $\s^{-1}\alpha$ is an equivalence tells
us that the map $\s^{-1}r:\s^{-1}X\ri\s^{-1}Y$ is injective,
with cokernel $ \s^{-1}C^{-1}$. The fact that $\alpha^{-1}\beta$
agrees with $d'$ means that the composite
\[
\CD
\s^{-1}C^{0} @>\s^{-1}g>> \s^{-1}Y @>>> \text{\rm Coker}(\s^{-1}r) 
\endCD
\]
is just the map $d:\s^{-1}C^{0}\ri \s^{-1}C^{-1}$. 
Let $B$ be the chain complex
\[
\CD
@>>> 0  @>>> C^0 \oplus X@>{\begin{pmatrix} 
\partial & 0 \\ g & r
\end{pmatrix}}>> C^{1}\oplus Y @>>>
\cdots @>>> C^n @>>> 0 @>>>
\endCD
\]
There is a natural map $f:B\ri C$, and $\s^{-1}f$ is a homology 
isomorphism of bounded complexes of projectives, hence a
homotopy equivalence. Thus
$\s^{-1}B$ is homotopy equivalent to $\s^{-1}C\cong D$.

This permits us to shorten on the left. Shortening the complex on the
right is dual.
\eprf

\section{Torsion modules}
\label{Storsion}

Until now all our theorems were general, in the sense that we did not
impose any restrictions on the ring $R$ or on the set of maps $\s$.  

\hyp{inj1} In this section, we assume that all the
morphisms in $\s$ are injections.
\ehyp

The main theorem of this section is that, under the above
restriction, the higher Waldhausen \kth\ of the triangulated category
$\ca^c=D^c(R,\sigma)$ agrees with the higher Quillen \kth\ of the exact
category $\ce=H(R,\sigma)$ of $\sigma$-torsion $R$-modules with
projective dimension $\leq 1$.
\medskip

\pro{inj2} If $R \ri \sigma^{-1}R$ is an injection then every $s:P \to
Q$ in $\sigma$ is an injection, i.e.  Hypothesis \ref{inj1} is
satisfied.
\epro
\prf Since $R \ri\rs$ is a monomorphism and $P$ is projective
and therefore flat, we deduce that
\[
\CD
R\otimes_R^{}P @>>> \br \otimes_R^{}P
\endCD
\]
is a monomorphism. In other words, the map $P\ri\s^{-1}P$ is mono.
Consider the commutative diagram
$$\xymatrix@C+10pt@R+10pt{P \ar[r]^{\displaystyle{s}} \ar[d] & Q \ar[d] \\
\sigma^{-1}P \ar[r]^-{\s^{-1}s} & \sigma^{-1}Q}$$
By the above, $P \ri \sigma^{-1}P$ is an injection. Now $\s^{-1}s:
\sigma^{-1}P \ri \sigma^{-1}Q$ is an isomorphism. From the
commutativity of the square we deduce that $s:P \to Q$
is an injection.
\eprf

\exm{inj3} The converse of Proposition~\ref{inj2}
does not hold in general. 
The set $\sigma=\{0 \ri R\}$ satisfies Hypothesis \ref{inj1},
but $R \ri \sigma^{-1}R=0$ is not injective.
\eexm

\dfn{Dtorsion.1.6}  An {\em $(R,\sigma)$-module} $T$ is an
$R$-module which admits a f.g.  projective $R$-module resolution
of length 1
$$\xymatrix{0 \ar[r]& P \ar[r]^{\displaystyle{s}} & Q \ar[r] &
T\ar[r]&0}$$ 
with $\s^{-1}s:\sigma^{-1}P \ri \sigma^{-1}Q$ a $\sigma^{-1}R$-module 
isomorphism. Let $\ce=H(R,\s)$ be the full subcategory of the category of
$R$-modules with objects the $(R,\sigma)$-modules.  \edfn

\rmk{cell2} 
An $R$-module $T$ is an $(R,\sigma)$-module if and only if~:
\begin{roenumerate}
\item ${\rm Tor}^R_i(\sigma^{-1}R,T)~=~0~$ for all $i\in\zz$.
In particular,
$\sigma^{-1}T~=~{\rm Tor}_0^R(\sigma^{-1}R,T)~=~0~.$
\item $T$ has
projective dimension $\leq 1$.
\item $T$ is finitely presented.
\end{roenumerate}
\ermk

\lem{Ltorsion.1.7}
The category $\ce$ of $(R,\sigma)$-modules, as in
Definition~\ref{Dtorsion.1.6}, is closed under extensions and kernels. 
Furthermore, it is idempotent complete; concretely, any direct summand
of an object in $\ce$ lies in $\ce$.  
\elem

\prf
Suppose we are given a short exact sequence of $R$-modules
\[
\CD
0 @>>> T' @>>> T @>>> T'' @>>> 0~.
\endCD
\]
The long exact sequence for $\Tor$ tells us that
if two of the terms lie in $\ce$, then for all
$i\in\zz$
\[
\Tor^R_i(\rs,T')~=~\Tor^R_i(\rs,T)~=~\Tor^R_i(\rs,T'')~=~0~.
\]
It is clear that
if $T''$ and $T$ are finitely presented $R$-modules of projective
dimension $\leq1$ then so is $T'$, and that if $T''$ and $T'$ are
finitely presented $R$-modules of projective dimension $\leq1$ then so
is $T$.  Hence $\ce$ is closed under kernels and extensions.\\
\indent
Suppose now that $T$ is an object of $\ce$, and that as
$R$-modules, $T=A\oplus B$. Since
\[
0~=~\Tor^R_i(\rs,T)~=~\Tor^R_i(\rs,A)\oplus \Tor^R_i(\rs,B)~,
\]
we deduce that, for all $i\in\zz$,
$\Tor^R_i(\rs,A)= \Tor^R_i(\rs,B)=0$. The projective dimensions
of $A$ and $B$ are bounded above by the projective dimension
of $T$, which is $\leq1$. Furthermore, we have an exact
sequence
\[\CD
T @>e>> T @>>> A @>>> 0
\endCD\]
where $e$ is an idempotent map on $T$. This expresses $A$ as
a quotient of two finitely presented modules. Hence
$A$ is finitely presented. Thus $A$ lies in $\ce$.
\eprf

\dfn{Dtorsion.3}
The bounded derived category of the exact category $\ce$, denoted
$D^b(\ce)$, is defined as follows.  The objects are bounded chain
complexes of objects of $\ce$.  The morphisms are obtained from the
chain maps by formally inverting that maps whose mapping cones are
acyclic (as complexes of $R$-modules).  There is an obvious functor
$i:D^b(\ce)\ri D(R)$.  
\edfn

\lem{Ltorsion.4}
The functor $i:D^b(\ce)\ri D(R)$ is fully faithful.
\elem

\prf
Let us begin by showing that, for any objects $T,T'\in\ce$ and any $n\in\zz$,
$$\{D^b(\ce)\}(T,\Sigma^n T') ~=~ \{D(R)\}(T,\Sigma^n T')~.$$
Take a map in $D^b(\ce)$ of the form $T\ri\T^n T'$. There
exists a bounded complex in $\ce$
which we call $X$, a quasi--isomorphism
$g:X\ri T$, and a map of complexes $f:X\ri \T^n T'$ so that
our map is $fg^{-1}_{}$.  That is, we have a complex
\[
\CD
\cdots @>>> X^{-2} @>>> X^{-1} @>>> X^{0} 
@>\partial_0^{}>> X^1 @>>> \cdots
\endCD
\]
There is a quasi--isomorphism $X\ri T$; in particular
$H^0(X)=T$. We have an exact sequence
\[
\CD
@>>> X^{-1} @>>> \text{\rm ker}(\partial_0^{}) @>>>  T @>>>0~.
\endCD
\]
But $T\in\ce$ means that $T$ is of projective
dimension $\leq1$. There is an exact sequence
\[
\CD
0 @>>> P @>s>> Q @>>> T @>>> 0
\endCD
\]
with  $P$ and $Q$ f.g. projective.
Since $P$ and $Q$ are projective $R$-modules,
there exists a map 
\[
\CD
 P @>>> Q @>>> T @>>> 0 \\
@VVV @VVV @V1VV \\
X^{-1} @>>> \text{\rm ker}(\partial_0^{}) @>>>  T @>>>0
\endCD
\]
Let $Z$ be given by the pushout square
\[
\CD
 P @>>> Q  \\
@VVV @VVV \\
X^{-1} @>>> Z .
\endCD
\]
The short exact sequence
\[
\CD
0 @>>> X^{-1} @>>> Z @>>> T @>>> 0
\endCD
\]
establishes that $Z\in\ce$ and that the complex
$X^{-1} \ri Z$ is quasi-isomorphic to $T$. We deduce
a quasi--isomorphism $h:X'\ri X$
of complexes, given below:
\[
\CD
\cdots @>>> 0 @>>> X^{-1} @>>> Z @>>> 0@>>>\cdots\\
@. @VVV @VVV @VVV @VVV @VVV \\
\cdots @>>> X^{-2} @>>> X^{-1} @>>> X^{0} @>>> X^1 @>>>\cdots
\endCD
\]
It follows that the map $fg^{-1}_{}:T\ri \T^n T'$ 
is equal to the map $\{fh\}{\{gh\}}^{-1}_{}$.
Since $X'$ is concentrated in degrees $0$ and $1$, it
follows that $fh$ vanishes unless $n=0$ or $1$.
Unless $n=0$ or $1$, we have proved that
$\{D^b(\ce)\}(T,\Sigma^n T')$ vanishes. As for
$$\{D(R)\}(T,\Sigma^n T')~=~\Ext^n_R(T,T')~,$$
it must vanish since the projective dimension of $T$ is
$\leq1$. In other words, for $n\neq0,1$ the equality
$$\{D^b(\ce)\}(T,\Sigma^n T')~=~\{D(R)\}(T,\Sigma^n T')$$
is just because both sides vanish.

We leave to the reader to check that the two sides
are equal also when $n=0$ or $1$. For $n=0$ both sides
identify as $\ce(T,T')$, while for $n=1$ both
sides identify as ${\rm Ext}^1_R(T,T')$. 

Let $T$ be an object of $\ce$.
Consider next the full subcategory $\ct\subset D^b(\ce)$
defined by
\[
\text{\rm Ob}(\ct)~=~
\left\{\begin{array}{l}$\,$ \\
Y\in\text{\rm Ob}(D^b(\ce)) \\
$\,$\end{array}\quad\right|\left.\quad \begin{array}{l}
\forall n\in\zz,\\
\{D^b(\ce)\}(T,\Sigma^n Y) \ri \{D(R)\}(T,\Sigma^n Y)\\
\text{\rm is an isomorphism}\end{array}\right\}.
\]
By the above, $\ct$ contains $\ce$, and clearly it is
triangulated. Hence $\ct$ contains all of $D^b(\ce)$.
Next, take any $Y$ in $D^b(\ce)$, and
 consider the full subcategory $\car\subset D^b(\ce)$
given by
\[
\text{\rm Ob}(\car)~=~
\left\{\begin{array}{l}$\,$ \\
X\in\text{\rm Ob}(D^b(\ce)) \\
$\,$\end{array}\quad\right|\left.\quad \begin{array}{l}
\forall n\in\zz,\\
\{D^b(\ce)\}(X,\Sigma^n Y) \ri \{D(R)\}(X,\Sigma^n Y)\\
\text{\rm is an isomorphism}\end{array}\right\}.
\]
By the above, $\ce\subset\car$, and $\car$ is clearly 
triangulated. Hence $\car$ contains $D^b(\ce)$.
\eprf

\lem{Lequivtorsion}
Assume that maps in $\s$ are all injections.
The natural map  $D^b(\ce)\ri D(R)$ 
factors through $\ca^c=D^c(R,\sigma)\subset D(R)$, and the
induced map $D^b(\ce)\ri \ca^c$
is an equivalence of categories.
\elem

\prf
Every object of $\ce$ is quasi--isomorphic to a
complex $0\ri P\ri Q\ri 0$ of f.g. projectives; that
is, $\ce\subset\cb^c$. Furthermore, every object
$e\in\ce$ satisfies $\br\oti e=0$, and by
Proposition~\ref{P5.3} this means $e\in\ca$. Therefore
$e\in\ca\cap\cb^c=\ca^c$. 

By Lemma~\ref{Ltorsion.4}, $D^b(\ce)$ is a full, triangulated
subcategory of $\cb=D(R)$. Clearly, it is the smallest
full, triangulated subcategory containing $\ce$.
Since $\ca^c$ contains $\ce$, it follows that $D^b(\ce)\subset\ca^c$.

We also know that the maps in $\s$ are injections.
If $s:P\ri Q$ lies in $\s$, then its cokernel lies in $\ce$,
and hence $D^b(\ce)\subset\ca^c$ contains all $s\in\s$. 
By \ref{T3.8.3}, $\ca^c$ is the smallest thick subcategory 
of $\cb$ containing $\s$. If we could prove that $D^b(\ce)$
is thick, it would follow that $\ca^c\subset D^b(\ce)$;
we are reduced to proving that $D^b(\ce)$ is thick.
But Lemma~\ref{Ltorsion.1.7} tells us that $\ce$ is
idempotent complete, and 
Theorem~2.8 of Balmer and Schlichting's~\cite{Balmer-Schlichting}
allows us to deduce that $D^b(\ce)$ is idempotent complete,
hence thick.
\eprf

\thm{Ttorsion.5}\label{torsionK}
Suppose every morphism in $\s$ is injective.
Then the algebraic \kth\ of the Waldhausen category $\ca^c=D^c(R,\sigma)$
is isomorphic to the algebraic \kth\ of the exact category $\ce=H(R,\sigma)$
$$K_*(R,\sigma)~=~K_{*-1}(D^c(R,\sigma))~=~K_{*-1}(H(R,\sigma))~.$$
\ethm

\prf
By Lemma~\ref{Lequivtorsion}, the natural map $D^b(\ce)\ri D(R)$
induces a triangulated equivalence of $D^b(\ce)$ with $\ca^c$. 
Hence the induced map in \kth\ is an isomorphism.
But Waldhausen's $K_i(D^b(\ce))$ agree with Quillen's
$K_i(\ce)$.
\eprf

\section{Algebraic $L$-theory}
\label{Ltheory}

We now extend our results to the algebraic $L$-theory of rings with
involution. We refer to Ranicki \cite{Ranicki1981},
\cite{Ranicki1998} more detailed expositions of algebraic $L$-theory.
\medskip

An {\it involution} on a ring $R$ is an anti-automorphism
$$R \ri R~;~r \mapsto \overline{r}~.$$
The involution is used to regard a left $R$-module $M$ as a right 
$R$-module by
$$M \times R \ri M ~;~(x,r) \mapsto \overline{r} x~.$$
The {\it dual} of a (left) $R$-module $M$ is the $R$-module
$$M^*~=~{\rm Hom}_R(M,R)~,~R \times  M^* \ri M^*~;~
(r,f) \mapsto (x \mapsto f(x)\overline{r})~.$$
The {\it dual} of an $R$-module morphism $s:P \ri Q$ is the $R$-module
morphism
$$s^*~:~Q^* \ri P^*~;~f \mapsto (x \mapsto f(s(x)))~.$$
If $M$ is f.g. projective then so is $M^*$, and 
$$M \ri M^{**}~;~x \mapsto (f \mapsto \overline{f(x)})$$
is an isomorphism which is used to identify $M^{**}=M$. 

\begin{hypothesis} \label{linj1}
In this section, we assume that 
\begin{itemize}
\item[(i)] $R$ is a ring with involution, 
\item[(ii)] the duals of morphisms $s:P \ri Q$ in $\sigma$
are morphisms $s^*:Q^* \ri P^*$ in $\sigma$, 
\item[(iii)] $\epsilon \in R$ is a central unit such that
$\overline{\epsilon}=\epsilon^{-1}$ {\rm (}e.g. $\epsilon= \pm 1${\rm )}.
\end{itemize}
The noncommutative localization $\sigma^{-1}R$ 
is then also a ring with involution, with $\epsilon \in \sigma^{-1}R$
a central unit such that $\overline{\epsilon}=\epsilon^{-1}$.
\hfill$\Box$
\end{hypothesis}

We review briefly the chain complex construction of the f.g. projective
$\epsilon$-quadratic $L$-groups $L_*(R,\epsilon)$
and the $\epsilon$-symmetric $L$-groups $L^*(R,\epsilon)$.
Given an $R$-module chain complex $C$ let
the generator $T \in \bZ_2$ act on the $\bZ$-module
chain complex $C\otimes_RC$ by the $\epsilon$-transposition duality
$$T_{\epsilon}~:~C_p \otimes_R C_q \ri C_q \otimes_R C_p~:~x \otimes y
\mapsto (-1)^{pq}\epsilon y \otimes x~.$$
Let $W$ be the standard free $\bZ[\bZ_2]$-module resolution of 
$\mathbb Z$
$$W~:~\dots \ri \bZ[\bZ_2] 
\xrightarrow[]{1-T} \bZ[\bZ_2] \xrightarrow[]{1+T} \bZ[\bZ_2]
\xrightarrow[]{1-T} \bZ[\bZ_2]~.$$
The {\it $\epsilon$-symmetric} (resp. {\it $\epsilon$-quadratic})
{\it $Q$-groups} of $C$ are the 
$\bZ_2$-hypercohomology (resp. $\bZ_2$-hyperhomology)
groups of $C \otimes_RC$ 
$$\begin{array}{l}
Q^n(C,\epsilon)~=~H^n(\bZ_2;C\otimes_RC)~=~
H_n({\rm Hom}_{\bZ[\bZ_2]}(W,C\otimes_RC))~,\\[1ex]
Q_n(C,\epsilon)~=~H_n(\bZ_2;C\otimes_RC)~=~
H_n(W\otimes_{\bZ[\bZ_2]}(C\otimes_RC))~.
\end{array}$$
The $Q$-groups are chain homotopy invariants of $C$.
There are defined forgetful maps
$$\begin{array}{l}
1+T_{\epsilon}~:~Q_n(C,\epsilon) \ri Q^n(C,\epsilon)~;~ 
\psi \mapsto (1+T_{\epsilon})\psi~,\\[1ex]
Q^n(C,\epsilon) \ri H_n(C \otimes_R C)~;~\phi
\mapsto \phi_0~.
\end{array}$$
For f.g. projective $C$ the function
$$C\otimes_RC \ri {\rm Hom}_R(C^*,C)~;~
x \otimes y \mapsto (f \mapsto \overline{f(x)}y)$$
is an isomorphism of $\bZ[\bZ_2]$-module chain complexes, with
$T \in \bZ_2$ acting on ${\rm Hom}_R(C^*,C)$ by $\theta \mapsto \epsilon \theta^*$.
The element $\phi_0 \in H_n(C\otimes_RC)=H_n({\rm Hom}_R(C^*,C))$ is
a chain homotopy class of $R$-module chain maps $\phi_0:C^{n-*} \ri C$.
\medskip

An {\it $n$-dimensional $\epsilon$-symmetric complex over $R$} $(C,\phi)$
is a bounded f.g. projective $R$-module chain complex $C$ together with
an element $\phi \in Q^n(C,\epsilon)$. The complex $(C,\phi)$ 
is {\it Poincar\'e} if the $R$-module chain map $\phi_0:C^{n-*} \ri C$ 
is a chain equivalence. 

\exm{form} 
A 0-dimensional $\epsilon$-symmetric Poincar\'e complex $(C,\phi)$ over $R$
is essentially the same as a nonsingular $\epsilon$-symmetric form 
$(M,\lambda)$ over $(R,\sigma)$, with $M=(C_0)^*$ a f.g. projective
$R$-module and 
$$\lambda~=~\phi_0~:~M \times M \ri R$$
a sesquilinear pairing such that the adjoint 
$$M \ri M^*~;~x \mapsto (y \mapsto \lambda(x,y))$$
is an $R$-module isomorphism.
\eexm

See pp.\ 210--211 of \cite{Ranicki1998} for the notion of an
{\it $\epsilon$-symmetric (Poincar\'e) pair}.
The {\it boundary} of an $n$-dimensional 
$\epsilon$-symmetric complex $(C,\phi)$ is the $(n-1)$-dimensional 
$\epsilon$-symmetric Poincar\'e complex 
$$\partial (C,\phi)~=~(\partial C,\partial \phi)$$
with $\partial C=C(\phi_0:C^{n-*} \ri C)_{*+1}$ and $\partial \phi$
as defined on p.\ 218 of \cite{Ranicki1998}.
The {\it $n$-dimensional $\epsilon$-symmetric $L$-group}
$L^n(R,\epsilon)$ is the cobordism group of $n$-dimensional
$\epsilon$-symmetric Poincar\'e complexes $(C,\phi)$ over $R$
with $C$ $n$-dimensional.
In particular, $L^0(R,\epsilon)$ is the Witt group of nonsingular
$\epsilon$-symmetric forms over $R$.
\medskip

An $n$-dimensional $\epsilon$-symmetric complex $(C,\phi)$ over $R$
is {\it $\sigma^{-1}R$-Poincar\'e} if the $\sigma^{-1}R$-module chain map 
$\sigma^{-1}\phi_0:\sigma^{-1}C^{n-*} \ri \sigma^{-1}C$
is a chain equivalence, in which case $\sigma^{-1}(C,\phi)$ is an
$n$-dimensional $\epsilon$-symmetric Poincar\'e complex over
$\sigma^{-1}R$.  
\medskip

The {\it $n$-dimensional $\epsilon$-symmetric $\Gamma$-group}
$\Gamma^n(R\ri \sigma^{-1}R,\epsilon)$ is the cobordism group of
$n$-dimensional $\epsilon$-symmetric $\sigma^{-1}R$-Poincar\'e
complexes $(C,\phi)$ over $R$ such that $\sigma^{-1}C$ is chain
equivalent to an $n$-dimensional induced f.g.  projective
$\sigma^{-1}R$-module chain complex.  The {\it $n$-dimensional
$\epsilon$-symmetric $L$-group} $L^n(R,\sigma,\epsilon)$ is the
cobordism group of $(n-1)$-dimensional $\epsilon$-symmetric Poincar\'e
complexes over $R$ $(C,\phi)$ such that $C$ is
$\sigma^{-1}R$-contractible, i.e.  $\sigma^{-1}C \simeq 0$.  \medskip

Similarly in the $\epsilon$-quadratic case, with groups
$L_n(R,\epsilon)$, $\Gamma_n(R \ri \sigma^{-1}R,\epsilon)$,
$L_n(R,\sigma,\epsilon)$. The $\epsilon$-quadratic $L$- and $\Gamma$-groups
are 4-periodic
$$\begin{array}{l}
L_n(R,\epsilon)~=~L_{n+2}(R,-\epsilon)~=~L_{n+4}(R,\epsilon)~,\\[1ex]
\Gamma_n(R\ri \sigma^{-1}R,\epsilon)~=~
\Gamma_{n+2}(R\ri \sigma^{-1}R,-\epsilon)~=~\Gamma_{n+4}(R\ri \sigma^{-1}R,\epsilon)~,\\[1ex]
L_n(R,\sigma,\epsilon)~=~L_{n+2}(R,\sigma,-\epsilon)~=~L_{n+4}(R,\sigma,\epsilon)~.
\end{array}$$

\pro{L1} For any ring with involution $R$ and noncommutative localization
$\sigma^{-1}R$ there is defined a localization exact sequence
of $\epsilon$-symmetric $L$-groups
$$\xymatrix{
\dots \ar[r] & L^n(R,\epsilon) \ar[r]&
\Gamma^n(R \ri \sigma^{-1}R,\epsilon) \ar[r]^-{\partial} &
L^n(R,\sigma,\epsilon)\ar[r] & L^{n-1}(R,\epsilon) \ar[r] & \dots~.}$$
Similarly in the $\epsilon$-quadratic case, with an exact sequence
$$\xymatrix{\dots \ar[r] & L_n(R,\epsilon) \ar[r]&
\Gamma_n(R \ri \sigma^{-1}R,\epsilon) \ar[r]^-{\partial} &
L_n(R,\sigma,\epsilon)\ar[r] & L_{n-1}(R,\epsilon) \ar[r] & \dots~.}$$
\epro
\prf The relative group of $L^n(R,\epsilon) \ri \Gamma^n(R \ri
\sigma^{-1}R,\epsilon)$ is the cobordism group of
$n$-dimensional $\epsilon$-symmetric $\sigma^{-1}R$-Poincar\'e
pairs over $R$ $(f:C \ri D,(\delta\phi,\phi))$ with $(C,\phi)$ Poincar\'e.
The effect of algebraic surgery on $(C,\phi)$ using this pair is a
cobordant $(n-1)$-dimensional $\epsilon$-symmetric Poincar\'e complex
$(C',\phi')$ with $C'$ $\sigma^{-1}R$-contractible. The function
$(f:C \ri D,(\delta\phi,\phi)) \mapsto (C',\phi')$ defines an isomorphism
between the relative group and $L^n(R,\sigma,\epsilon)$.
\eprf

Define
$$I~=~{\rm im}(K_0(R)\ri K_0(\sigma^{-1}R))~,$$
the subgroup of $K_0(\sigma^{-1}R)$ consisting of the projective classes 
of the f.g.  projective $\sigma^{-1}R$-modules induced from f.g.  projective 
$R$-modules. By definition, $L^n_I(\sigma^{-1}R,\epsilon)$ is the cobordism 
group of $n$-dimensional $\epsilon$-symmetric Poincar\'e complexes over 
$\sigma^{-1}R$
$(B,\theta)$ such that $[B] \in I$. 
There are evident morphisms of $\Gamma$- and $L$-groups
$$\begin{array}{l}
\sigma^{-1}\Gamma^*~:~\Gamma^n(R \ri \sigma^{-1}R,\epsilon) \ri 
L_I^n(\sigma^{-1}R,\epsilon)~;~(C,\phi) \mapsto \sigma^{-1}(C,\phi)~,\\[1ex]
\sigma^{-1}\Gamma_*~:~\Gamma_n(R \ri \sigma^{-1}R,\epsilon) \ri 
L^I_n(\sigma^{-1}R,\epsilon)~;~(C,\psi) \mapsto \sigma^{-1}(C,\psi)~.
\end{array}$$
In general, the morphisms $\sigma^{-1}\Gamma^*,\sigma^{-1}\Gamma_*$
need not be isomorphisms, since a bounded f.g.  projective
$\sigma^{-1}R$-module chain complex $D$ with $[D] \in I$ need not be
chain equivalent to $\sigma^{-1}C$ for a bounded f.g.  projective
$R$-module chain complex $C$ (cf.  the chain complex lifting problem
considered in sections  \ref{lift2} and \ref{Sequivalence}).  
\medskip

It was proved in Chapter 3 of Ranicki \cite{Ranicki1981} that 
if $R \ri \sigma^{-1}R$ is an injective Ore localization then the morphisms
$\sigma^{-1}Q^*,\sigma^{-1}Q_*$, 
$\sigma^{-1}\Gamma^*,\sigma^{-1}\Gamma_*$ are isomorphisms, 
so that there are defined localization exact sequences 
for both the $\epsilon$-symmetric and the $\epsilon$-quadratic $L$-groups
$$\xymatrix@R-23pt{
\dots \ar[r] & L^n(R,\epsilon) \ar[r]&
L^n_I(\sigma^{-1}R,\epsilon) \ar[r]^-{\partial} &
L^n(R,\sigma,\epsilon)\ar[r] & L^{n-1}(R,\epsilon) \ar[r] & \dots~,\\
\dots \ar[r] & L_n(R,\epsilon) \ar[r]&
L^I_n(\sigma^{-1}R,\epsilon) \ar[r]^-{\partial} &
L_n(R,\sigma,\epsilon)\ar[r] & L_{n-1}(R,\epsilon) \ar[r] & \dots~.}$$
Special cases of these sequences were obtained by Milnor-Husemoller,
Karoubi, Pardon, Smith, Carlsson-Milgram.
\medskip

For any bounded f.g. projective $R$-module chain complex $C$
the natural $R$-module chain map
$$\mathop{\varinjlim}\limits_{(B,\beta)} B ~=~G\pi(C) \ri \sigma^{-1}C$$
induces morphisms
$$\begin{array}{l}
\sigma^{-1}Q^*~:~\mathop{\varinjlim}\limits_{(B,\beta)}
Q^n(B,\epsilon)~=~Q^n(G\pi(C),\epsilon) \ri Q^n(\sigma^{-1}C,\epsilon)~,\\[1ex]
\sigma^{-1}Q_*~:~\mathop{\varinjlim}\limits_{(B,\beta)} 
Q_n(B,\epsilon)~=~Q_n(G\pi(C),\epsilon) \ri Q_n(\sigma^{-1}C,\epsilon)
\end{array} $$
with the direct limits taken over all the bounded f.g.  projective
$R$-module chain complexes $B$ with a chain map $\beta:C\ri B$ such
that $\sigma^{-1}\beta:\sigma^{-1}C \ri \sigma^{-1}B$ is a
$\sigma^{-1}R$-module chain equivalence.  The natural projection
$D\otimes_RD \ri D\otimes_{\sigma^{-1}R}D$ is an isomorphism for any
bounded f.g.  projective $\sigma^{-1}R$-module chain complex $D$ (since
this is already the case for $D=\sigma^{-1}R$), so the $Q$-groups of
$\sigma^{-1}C$ are the same whether $\sigma^{-1}C$ is regarded as an
$R$-module or $\sigma^{-1}R$-module chain complex.

\thm{L2}  {\rm (Vogel \cite{Vogel1982}, Theorem 8.4)}
For any ring with involution $R$ and noncommutative localization
$\sigma^{-1}R$ the morphisms
$$\sigma^{-1}\Gamma_*~:~\Gamma_n(R \ri \sigma^{-1}R,\epsilon) \ri 
L^I_n(\sigma^{-1}R,\epsilon)~;~(C,\psi) \mapsto \sigma^{-1}(C,\psi)$$
are isomorphisms, and there is a localization exact sequence
of $\epsilon$-quadratic $L$-groups
$$\xymatrix{
\dots \ar[r] & L_n(R,\epsilon) \ar[r]&
L^I_n(\sigma^{-1}R,\epsilon) \ar[r]^-{\partial} &
L_n(R,\sigma,\epsilon)\ar[r] & L_{n-1}(R,\epsilon) \ar[r] & \dots~.}$$
\ethm
\prf By algebraic surgery below the middle dimension it suffices
to consider only the special cases $n=0,1$. In effect, it was proved
in \cite{Vogel1982} that $\sigma^{-1}Q_*$ is an isomorphism for
0- and 1-dimensional $C$. 
\eprf

It was claimed in Proposition 25.4 of Ranicki \cite{Ranicki1998}
that $\sigma^{-1}\Gamma^*$ is also an isomorphism, assuming (incorrectly)
that the chain complex lifting problem can always be solved.
However, we do have :

\thm{L3} 
If $\sigma^{-1}R$ is a noncommutative localization of a ring with
involution $R$ which is stably flat over $R$, there is a localization
exact sequence of $\epsilon$-symmetric $L$-groups
$$\xymatrix{
\dots \ar[r] & L^n(R,\epsilon) \ar[r]&
L^n_I(\sigma^{-1}R,\epsilon) \ar[r]^-{\partial} &
L^n(R,\sigma,\epsilon)\ar[r] & L^{n-1}(R,\epsilon) \ar[r] & \dots~.}$$
\ethm
\prf By Theorem \ref{Tfinal} for any bounded f.g. projective
$R$-module chain complex $C$ the natural $R$-module chain map $G\pi(C)
\ri \sigma^{-1}C$ induces isomorphisms in homology
$$H_*(G\pi(C))~\cong~H_*(\sigma^{-1}C)~.$$
Thus the natural $\bZ[\bZ_2]$-module chain map
$$G\pi(C)\otimes_RG\pi(C) \ri \sigma^{-1}C\otimes_R\sigma^{-1}C~=~
\sigma^{-1}C\otimes_{\sigma^{-1}R}\sigma^{-1}C$$
induces isomorphisms of $\epsilon$-symmetric $Q$-groups
$$\sigma^{-1}Q^*~:~\mathop{\varinjlim}\limits_{(B,\beta)} 
Q^n(B,\epsilon) \ri Q^n(\sigma^{-1}C,\epsilon)$$
(and also isomorphisms $\sigma^{-1}Q_*$ of $\epsilon$-quadratic $Q$-groups).
By Proposition \ref{Prelifting} every $n$-dimensional induced f.g. projective
$\sigma^{-1}R$-module chain complex $D$ is chain equivalent to $\sigma^{-1}C$
for an $n$-dimensional f.g. projective $R$-module chain complex $C$, with
$$Q^n(D,\epsilon)~=~Q^n(\sigma^{-1}C,\epsilon)~=~
\mathop{\varinjlim}\limits_{(B,\beta)}Q^n(B,\epsilon) ~.$$
It follows that the morphisms of $\epsilon$-symmetric $\Gamma$- and $L$-groups
$$\sigma^{-1}\Gamma^*~:~\Gamma^n(R \ri \sigma^{-1}R,\epsilon) \ri 
L^n_I(\sigma^{-1}R,\epsilon)~;~(C,\phi) \mapsto \sigma^{-1}(C,\phi)$$ 
are also isomorphisms, and the localization exact sequence is given
by Proposition \ref{L1}.
\eprf

\begin{hypothesis} \label{linj2}
For the remainder of this section, we assume Hypothesis \ref{linj1}
and also that $R \ri \sigma^{-1}R$ is an injection.\hfill$\Box$
\end{hypothesis}
As in Proposition  \ref{inj2} it follows that all the morphisms in $\s$ are injections.
\medskip

We shall now obtain the $L$-theoretic analogue of the algebraic
$K$-theory identification
$K_*(R,\sigma)=K_{*-1}(H(R,\sigma))$ obtained in section
\ref{Storsion}, with $H(R,\sigma)$ the exact category of 
$(R,\sigma)$-modules.
We generalize the results of Ranicki
\cite{Ranicki1981} and Vogel \cite{Vogel1980} to prove that under
Hypotheses \ref{linj1},\ref{linj2} the relative $L$-groups
$L^*(R,\sigma,\epsilon)$, $L_*(R,\sigma,\epsilon)$ in the $L$-theory
localization exact sequences are the $L$-groups of $H(R,\sigma)$ with
respect to the following duality involution.
\medskip

Define the {\it torsion dual} 
of an $(R,\sigma)$-module $M$ to be the $(R,\s)$-module
$$M\widehat{~}~=~{\rm Ext}^1_R(M,R)~,$$
using the involution on $R$ to define the left $R$-module
structure. If $M$ has f.g. projective $R$-module resolution
$$0 \ri P_1 \xrightarrow[]{s} P_0 \ri M \ri 0$$
with $s \in \sigma$ the torsion dual $M\widehat{~}$ has the dual f.g. 
projective $R$-module resolution
$$0 \ri P_0^* \xrightarrow[]{s^*} P_1^* \ri M\widehat{~} \ri 0$$
with $s^* \in \sigma$. 
\medskip

\pro{iso} Let $M={\rm coker}(s:P_1 \ri P_0)$,
$N={\rm coker}(t:Q_1 \ri Q_0)$ be $(R,\sigma)$-modules.\\
{\rm (i)} The adjoint of the pairing
$$M \times M\widehat{~} \ri \sigma^{-1}R/R~;~(g \in P_0,f\in P_1^*)
\mapsto fs^{-1}g$$
defines a natural $R$-module isomorphism
$$M\widehat{~} \ri {\rm Hom}_R(M,\sigma^{-1}R/R)~;~f \mapsto (g \mapsto fs^{-1}g)~.$$
{\rm (ii)} The natural $R$-module morphism
$$M \ri M\widehat{~}\widehat{~}~;~x \mapsto (f \mapsto \overline{f(x)})$$
is an isomorphism.\\
{\rm (iii)} There are natural identifications
$$\begin{array}{l}
M \otimes_R N~=~{\rm Tor}^R_0(M,N)~=~{\rm Ext}^1_R(M\widehat{~},N)~=~
H_0(P\otimes_RQ)~,\\[1ex]
{\rm Hom}_R(M\widehat{~},N)~=~{\rm Tor}^R_1(M,N)~=~
{\rm Ext}^0_R(M\widehat{~},N)~=~H_1(P\otimes_RQ)~.
\end{array}$$
The functions
$$\begin{array}{l}
M \otimes_R N \ri N \otimes_R M~;~x \otimes y \mapsto y \otimes x~,\\[1ex]
{\rm Hom}_R(M\widehat{~},N) \ri {\rm Hom}_R(N\widehat{~},M)~;~
f \mapsto f\widehat{~}
\end{array}$$
determine transposition isomorphisms
$$T~:~{\rm Tor}_i^R(M,N) \ri {\rm Tor}_i^R(N,M)~~(i=0,1)~.$$
{\rm (iv)} For any finite subset
$V=\{v_1,v_2,\dots,v_k\} \subset M\otimes_RN$ there exists an exact sequence of
$(R,\sigma)$-modules 
$$0 \ri N \ri L  \ri \oplus_k M\widehat{~} \ri 0$$
such that $V \subset {\rm ker}(M\otimes_RN \ri M\otimes_RL)$.
\epro
\prf (i)  Apply the snake lemma to the morphism of short exact sequences
$$\xymatrix{
0 \ar[r] &{\rm Hom}_R(P_0,R) \ar[r] \ar[d]^-{\displaystyle{s^*_{}}} &
{\rm Hom}_R(P_0,\sigma^{-1}R) \ar[r] \ar[d]^-{\displaystyle{s^*_1}} &
{\rm Hom}_R(P_0,\sigma^{-1}R/R) \ar[r] \ar[d]^-{\displaystyle{s^*_2}}& 0 \\
0 \ar[r] &{\rm Hom}_R(P_1,R) \ar[r] &
{\rm Hom}_R(P_1,\sigma^{-1}R) \ar[r] &
{\rm Hom}_R(P_1,\sigma^{-1}R/R) \ar[r] & 0}$$
with $s^*$ injective, $s^*_1$ an isomorphism and $s^*_2$ surjective, 
to verify that the $R$-module morphism
$$M\widehat{~}~=~{\rm coker}(s^*) \ri
{\rm Hom}_R(M,\sigma^{-1}R/R)~=~{\rm ker}(s^*_2)$$
is an isomorphism.\\
(ii) Immediate from the identification
$$s^{**}~=~s~:~(P_0)^{**}~=~P_0 \ri (P_1)^{**}~=~P_1~.$$
(iii) Exercise for the reader.\\
(iv) Lift each $v_i \in M\otimes_RN$ to an element
$$v_i \in P_0 \otimes_RQ_0~=~{\rm Hom}_R(P_0^*,Q_0)~~(1 \leq i \leq k)~.$$
The $R$-module morphism defined by
$$u~=~\begin{pmatrix}s^* & 0 & 0 & \dots & 0 \\
0 & s^* & 0 & \dots & 0 \\
0 & 0 & s^* & \dots & 0 \\
\vdots & \vdots & \vdots &\ddots & \vdots \\
v_1 & v_2 & v_3 & \dots & t \end{pmatrix}~:~
U_1~=~(\oplus_k P_0^*) \oplus 
Q_1 \ri U_0~=~(\oplus_k P_1^*)\oplus Q_0$$
is in $\sigma$, so that $L={\rm coker}(u)$ is an $(R,\sigma)$-module with
a f.g. projective $R$-module resolution
$$0 \ri U_1 \xrightarrow[]{u} U_0 \ri L \ri 0~.$$
The short exact sequence of 1-dimensional f.g. projective $R$-module chain
complexes
$$0 \ri Q \ri U \ri \oplus_k P^{1-*} \ri 0$$
is a resolution of a short exact sequence of $(R,\sigma)$-modules
$$0 \ri N \ri L  \ri \oplus_k M\widehat{~} \ri 0~.$$
The first morphism in the exact sequence
$${\rm Tor}^R_1(M,\oplus_k M\widehat{~}) \ri M\otimes_RN \ri
M\otimes_RL \ri M \otimes_R(\oplus_k M\widehat{~}) \ri 0$$
sends $1_i \in {\rm Tor}^R_1(M,\oplus_kM\widehat{~})=
\oplus_k {\rm Hom}_R(M\widehat{~},M\widehat{~})$
to $v_i \in {\rm ker}(M\otimes_RN \ri M\otimes_RL)$.
\eprf

Given an  $(R,\sigma)$-module chain complex $C$ define the {\it
$\epsilon$-symmetric} (resp.  {\it $\epsilon$-quadratic}) {\it torsion
$Q$-groups} of $C$ to be the $\bZ_2$-hypercohomology (resp. 
$\bZ_2$-hyperhomology) groups of the $\epsilon$-transposition
involution $T_{\epsilon}=\epsilon T$ on the $\bZ$-module chain complex
${\rm Tor}_1^R(C,C)={\rm Hom}_R(C\widehat{~},C)$
$$\begin{array}{l}
Q_{\rm tor}^n(C,\epsilon)~=~H^n(\bZ_2;{\rm Tor}^R_1(C,C))~=~
H_n({\rm Hom}_{\bZ[\bZ_2]}(W,{\rm Tor}^R_1(C,C)))~,\\[1ex]
Q^{\rm tor}_n(C,\epsilon)~=~H_n(\bZ_2;{\rm Tor}^R_1(C,C))~=~
H_n(W\otimes_{\bZ[\bZ_2]}({\rm Tor}^R_1(C,C)))~.
\end{array}$$
There are defined forgetful maps
$$\begin{array}{l}
1+T_{\epsilon}~:~Q^{\rm tor}_n(C,\epsilon) \ri Q^n_{\rm tor}(C,\epsilon)~;~ 
\psi \mapsto (1+T_{\epsilon})\psi~,\\[1ex]
Q_{\rm tor}^n(C,\epsilon) \ri H_n({\rm Tor}^R_1(C,C))~;~\phi
\mapsto \phi_0~.
\end{array}$$
The element $\phi_0 \in H_n({\rm Tor}^R_1(C,C))$ is
a chain homotopy class of $R$-module chain maps 
$\phi_0:C^{n-}\widehat{~} \ri C$.
\medskip

An {\it $n$-dimensional $\epsilon$-symmetric complex over $(R,\sigma)$} 
$(C,\phi)$ is a bounded $(R,\sigma)$-module chain complex $C$ together with
an element $\phi \in Q_{\rm tor}^n(C,\epsilon)$. The complex $(C,\phi)$ 
is {\it Poincar\'e} if the $R$-module chain maps 
$\phi_0:C^{n-}\widehat{} \ri C$ are chain equivalences. 

\exm{linking form} 
A 0-dimensional $\epsilon$-symmetric Poincar\'e complex $(C,\phi)$ over $(R,\sigma)$
is essentially the same as a nonsingular $\epsilon$-symmetric linking form 
$(M,\lambda)$ over $(R,\sigma)$, with $M=(C_0)\widehat{~}$ an 
$(R,\s)$-module and 
$$\lambda~=~\phi_0~:~M \times M \ri \sigma^{-1}R/R$$
a sesquilinear pairing such that the adjoint 
$$M \ri M\widehat{~}~;~x \mapsto (y \mapsto \lambda(x,y))$$
is an $R$-module isomorphism.
\eexm

The {\it $n$-dimensional torsion $\epsilon$-symmetric $L$-group}
$L_{\rm tor}^n(R,\sigma,\epsilon)$ is the cobordism group of
$n$-dimensional $\epsilon$-symmetric Poincar\'e complexes $(C,\phi)$ over
$(R,\sigma)$, with $C$ $n$-dimensional.  In particular, 
$L^0_{\rm tor}(R,\sigma,\epsilon)$ is the Witt group of nonsingular
$\epsilon$-symmetric linking forms over $(R,\sigma)$.  
\medskip

Similarly in the $\epsilon$-quadratic case, with torsion $L$-groups 
$L^{\rm tor}_n(R,\sigma,\epsilon)$.  The $\epsilon$-quadratic torsion
$L$-groups are 4-periodic
$$L^{\rm tor}_n(R,\sigma,\epsilon)~=~
L^{\rm tor}_{n+2}(R,\sigma,-\epsilon)~=~L^{\rm tor}_{n+4}(R,\sigma,\epsilon)~.$$

\thm{Lfin} 
If $R \ri \sigma^{-1}R$ is injective the relative $L$-groups in the localization exact 
sequences of Proposition \ref{L1} 
$$\xymatrix@R-23pt{
\dots \ar[r] & L^n(R,\epsilon) \ar[r]&
\Gamma^n(R \ri \sigma^{-1}R,\epsilon) \ar[r]^-{\partial} &
L^n(R,\sigma,\epsilon)\ar[r] & L^{n-1}(R,\epsilon) \ar[r] & \dots\\
\dots \ar[r] & L_n(R,\epsilon) \ar[r]&
\Gamma_n(R \ri \sigma^{-1}R,\epsilon) \ar[r]^-{\partial} &
L_n(R,\sigma,\epsilon)\ar[r] & L_{n-1}(R,\epsilon) \ar[r] & \dots}
$$
are the torsion $L$-groups
$$\begin{array}{l}
L^*(R,\sigma,\epsilon)~=~L_{\rm tor}^*(R,\sigma,\epsilon)~,\\[1ex]
L_*(R,\sigma,\epsilon)~=~L^{\rm tor}_*(R,\sigma,\epsilon)~.
\end{array}$$
\ethm
\prf For any bounded $(R,\sigma)$-module chain complex $T$
there exists a bounded f.g. projective $R$-module chain complex $C$
with a homology equivalence $C \ri T$. Working as in \cite{Vogel1980}
there is defined a distinguished triangle of $\bZ[\bZ_2]$-module chain complexes
$$\Sigma {\rm Tor}^R_1(T,T) \ri
C\otimes_RC \ri T\otimes_RT \ri \Sigma^2{\rm Tor}^R_1(T,T)$$
with $\bZ_2$ acting by the $\epsilon$-transposition $T_{\epsilon}$
on the $\bZ$-module chain complex ${\rm Tor}^R_1(T,T)$ 
and by the $(-\epsilon)$-transpositions $T_{-\epsilon}$
on $C\otimes_RC$ and $T\otimes_RT$, inducing long exact sequences
$$\xymatrix@R-23pt{
\dots \ar[r] & Q^n_{\rm tor}(T,\epsilon) \ar[r]&
Q^{n+1}(C,-\epsilon) \ar[r] &
Q^{n+1}(T,-\epsilon) \ar[r] & Q^{n-1}_{\rm tor}(T,\epsilon) \ar[r] & \dots\\
\dots \ar[r] & Q_n^{\rm tor}(T,\epsilon) \ar[r]&
Q_{n+1}(C,-\epsilon) \ar[r] &
Q_{n+1}(T,-\epsilon) \ar[r] & Q_{n-1}^{\rm tor}(T,\epsilon) \ar[r] & \dots~.}
$$
Passing to the direct limits over all the bounded $(R,\sigma)$-module
chain complexes $U$ with a homology equivalence $\beta:T \ri U$ use
Proposition \ref{iso} (iv) to obtain
$$\begin{array}{l}
\mathop{\varinjlim}\limits_{(U,\beta)} Q^{n+1}(U,-\epsilon)~=~0~,\\[1ex]
\mathop{\varinjlim}\limits_{(U,\beta)} Q_{n+1}(U,-\epsilon)~=~0
\end{array}$$
and hence
$$\begin{array}{l}
\mathop{\varinjlim}\limits_{(U,\beta)} 
Q^n_{\rm tor}(U,\epsilon)~=~Q^{n+1}(C,-\epsilon)~,\\[1ex]
\mathop{\varinjlim}\limits_{(U,\beta)} 
Q_n^{\rm tor}(U,\epsilon)~=~Q_{n+1}(C,-\epsilon)~.
\end{array}$$
\eprf

\rmk{final} The identification $L_*(R,\sigma,\epsilon)=L^{\rm
tor}_*(R,\sigma,\epsilon)$ for noncommutative $\sigma^{-1}R$ was first
obtained by Vogel \cite{Vogel1980}.  
\ermk

\ifx\undefined\bysame
\newcommand{\bysame}{\leavevmode\hbox to3em{\hrulefill}\,}
\fi

\end{document}